\documentstyle[12pt]{article}
\makeatletter

\def\URA#1{\begin{picture}(3,2)\put(0,0){\vector(3,2){3}}%
        \put(1.5,.5){${\scriptstyle{#1}}$}\end{picture}}
\def\ULA#1{\begin{picture}(3,2)\put(3,0){\vector(-3,2){3}}%
        \put(1.5,1.2){${\scriptstyle{#1}}$}\end{picture}}
\def\DRA#1{\begin{picture}(3,2)\put(0,2){\vector(3,-2){3}}%
        \put(1.5,1.2){${\scriptstyle{#1}}$}\end{picture}}

\def\LRA#1#2{\@tempdimb=\c@enumiv\@tempdima%
   \vcenter{\offinterlineskip\halign{##\cr%
   \hfil${\scriptstyle{#1}}$\hfil\crcr%
   \hbox to \@tempdimb{\rightarrowfill}\cr%
   \noalign{\kern-1ex}%
   \hbox to \@tempdimb{\leftarrowfill}\cr%
   \hfil${\scriptstyle{#2}}$\hfil\crcr}}}
\def\RA#1{\@tempdimb=\c@enumiv\@tempdima\vbox{\offinterlineskip%
   \halign{##\cr\hfil${\scriptstyle {#1}}$\hfil\crcr%
   \hbox to \@tempdimb{\rightarrowfill}\cr}}}
\def\LA#1{\@tempdimb=\c@enumiv\@tempdima\vbox{\offinterlineskip%
   \halign{##\cr\hfil${\scriptstyle {#1}}$\hfil\crcr%
   \hbox to \@tempdimb{\leftarrowfill}\cr}}}
\def\UA#1{\strut\begin{picture}(1,2)\put(.5,0){\vector(0,1){2}}%
        \put(.7,.8){${\scriptstyle{#1}}$}\end{picture}}
\def\DA#1{\strut\begin{picture}(1,2)\put(.5,2){\vector(0,-1){2}}%
        \put(.7,.8){${\scriptstyle{#1}}$}\end{picture}}
\def\diag{\leavevmode\bgroup\setcounter{enumiv}{1}%
   \unitlength1em \@tempdima3em \def\\{\crcr&}\vbox\bgroup%
   \def\multicolumn##1##2{\multispan##1\setcounter{enumiv}{##1}%
   \hfil{##2}\hfil\setcounter{enumiv}{1}}
   \offinterlineskip\halign\bgroup\vrule height.8em depth.7em %
   width0pt##&&\hfil${\displaystyle{##}}$\hfil\cr&}
\def\enddiag{\crcr\egroup\egroup\egroup}
\makeatother
\if@twoside \oddsidemargin -20pt \evensidemargin 20pt \marginparwidth 85pt
\else \oddsidemargin 0pt \evensidemargin 0pt
\fi
\advance\textheight by \topskip
\font\symbolfont=msbm10
\font\p=msbm10 at 12pt
\font\teneu=eufm10
\font\egteu=eufm8
\def\dn#1{\mathchoice{\hbox{\teneu #1}}{\hbox{\teneu #1}}%
   {\hbox{\egteu #1}}{\hbox{\egteu #1}}}
\font\symbolfont=msbm10

\def\fff{\hbox{\symbolfont F}}
\def\nnn{\hbox{\symbolfont N}}
\def\ppp{\hbox{\symbolfont P}}

\def\ggg{\hbox{\symbolfont G}}
\def\zzz{\hbox{\symbolfont Z}}
\def\cc{\mathop{\cal C}\nolimits}
\def\dd{\mathop{\cal D}\nolimits}
\def\ff{\mathop{\cal F}\nolimits}
\def\ii{\mathop{\cal I}\nolimits}

\def\tt{\mathop{\cal T}\nolimits}

\def\oo{\mathop{\cal O}\nolimits}
\def\pp{\mathop{\cal P}\nolimits}
\def\rr{\mathop{\cal R}\nolimits}

\def\ma{\mathop{\rm M}\nolimits}
\def\ind{\mathop{\rm ind}\nolimits}
\def\ob{\mathop{\rm ob}\nolimits}

\def\lat{\mathop{\rm lat}\nolimits}
\def\mod{\mathop{\rm mod}\nolimits}
\def\Mod{\mathop{\rm Mod}\nolimits}

\def\sub{\mathop{\rm sub}\nolimits}
\def\length{\mathop{\rm length}\nolimits}

\def\pr{\mathop{\rm pr}\nolimits}
\def\rin{\mathop{\rm in}\nolimits}
\def\pd{\mathop{\rm pd}\nolimits}
\def\hom{\mathop{\rm Hom}\nolimits}
\def\ext{\mathop{\rm Ext}\nolimits}

\def\endm{\mathop{\rm End}\nolimits}
\def\Cok{\mathop{\rm Cok}\nolimits}
\def\Ker{\mathop{\rm Ker}\nolimits}
\def\Im{\mathop{\rm Im}\nolimits}
\def\soc{\mathop{\rm soc}\nolimits}
\def\add{\mathop{\rm add}\nolimits}

\def\resdim#1#2{\mathop{{#1}\mbox{-{\rm resol.dim}}\strut\kern2pt {#2}}\nolimits}
\def\plim{\mathop{\raise-.2em\hbox{$\def\arraystretch{0}\begin{array}{c}\lim\\ {\scriptstyle\stackrel{\longleftarrow}{}}\end{array}$}}}
\font\p=msbm10 at 12pt

\def\subsetneq{\mathop{\mbox{ {\p\char'050} }}\nolimits}

\def\XZA{0.1}
\def\XZB{0.2}
\def\XZC{0.3}
\def\XA{1}
\def\XAA{1.1}
\def\XAAA{1.1.1}
\def\XAAB{1.1.2}
\def\XAAC{1.1.3}
\def\XAAD{1.1.4}
\def\XAB{1.2}
\def\XABA{1.2.1}
\def\XAC{1.3}
\def\XACA{1.3.1}
\def\XAD{1.4}
\def\XADA{1.4.1}
\def\XADB{1.4.2}
\def\XAE{1.5}
\def\XAEA{1.5.1}
\def\XAEB{1.5.2}
\def\XAEC{1.5.3}
\def\XAF{1.6}
\def\XAFA{1.6.1}
\def\XAFB{1.6.2}
\def\XAFC{1.6.3}
\def\XAFD{1.6.4}
\def\XAFE{1.6.5}
\def\XAG{1.7}
\def\XB{2}
\def\XBA{2.1}
\def\XBB{2.2}
\def\XBBA{2.2.1}
\def\XBBB{2.2.2}
\def\XBBC{2.2.3}
\def\XBC{2.3}
\def\XBCA{2.3.1}
\def\XBCB{2.3.2}
\def\XBCC{2.3.3}
\def\XBCD{2.3.4}
\def\XBD{2.4}
\def\XBDA{2.4.1}
\def\XBDB{2.4.2}
\def\XBE{2.5}
\def\XBEA{2.5.1}
\def\XBEB{2.5.2}
\def\XC{3}
\def\XCA{3.1}
\def\XCB{3.2}
\def\XCBA{3.2.1}
\def\XCBB{3.2.2}
\def\XCBC{3.2.3}
\def\XCC{3.3}
\def\XCCA{3.3.1}
\def\XCD{3.4}
\def\XCDA{3.4.1}
\def\XCE{3.5}
\def\XCEA{3.5.1}
\def\XCEB{3.5.2}
\def\XCEC{3.5.3}
\def\XCED{3.5.4}
\def\XCF{3.6}
\def\XCFA{3.6.1}
\def\XCFB{3.6.2}
\def\XCG{3.7}
\def\XD{4}
\def\XDA{4.1}
\def\XDAA{4.1.1}
\def\XDAB{4.1.2}
\def\XDAC{4.1.3}
\def\XDAD{4.1.4}
\def\XDB{4.2}
\def\XDC{4.3}
\def\XDD{4.4}
\def\XDDA{4.4.1}
\def\XDDB{4.4.2}
\def\XDDC{4.4.3}
\def\XDE{4.5}
\def\XDEA{4.5.1}
\def\XDF{4.6}
\def\XDFA{4.6.1}
\def\XDFB{4.6.2}
\def\XDFC{4.6.3}
\def\XDFD{4.6.4}
\def\XDG{4.7}
\def\XDGA{4.7.1}
\def\XDGB{4.7.2}
\def\XDH{4.8}
\def\XDHA{4.8.1}
\def\XDHB{4.8.2}
\def\XDI{4.9}
\def\XE{5}
\def\XEA{5.1}
\def\XEAA{5.1.1}
\def\XEAB{5.1.2}
\def\XEAC{5.1.3}
\def\XEB{5.2}
\def\XEBA{5.2.1}
\def\XEBB{5.2.2}
\def\XEC{5.3}
\def\XECA{5.3.1}
\def\XECB{5.3.2}
\def\XED{5.4}
\def\XEDA{5.4.1}
\def\XEE{5.5}
\def\XEEA{5.5.1}
\def\XEEB{5.5.2}
\def\XEEC{5.5.3}
\def\XEF{5.6}

\def\xx{\mathop{\cal X}\nolimits}
\def\yy{\mathop{\cal Y}\nolimits}

\def\tr{\mathop{\rm Tr}\nolimits}
\def\pd{\mathop{\rm pd}\nolimits}
\def\id{\mathop{\rm id}\nolimits}
\def\gl{\mathop{\rm gl.dim}\nolimits}
\def\expdim{\mathop{\rm exp.dim}\nolimits}
\def\wresdim{\mathop{\rm wresol.dim}\nolimits}
\def\rdim{\mathop{\rm rep.dim}\nolimits}
\def\rwrdim{\mathop{\rm rwrep.dim}\nolimits}
\def\fdim{\mathop{\rm fin.dim}\nolimits}

\def\LL{\mathop{\rm LL}\nolimits}

\def\pr{\mathop{\rm pr}\nolimits}
\begin{document}
\begin{center}
\vspace*{1cm}{\large Rejective subcategories of artin algebras and orders}
\footnote{2000 {\it Mathematics Subject Classification.}
Primary 16E10; Secondary 16G10, 16G30}
\vskip1em{\sc Osamu Iyama}
\end{center}

{\footnotesize
{\sc Abstract. }We will study the resolution dimension of functorially finite subcategories. The subcategories with the resolution dimension zero correspond to ring epimorphisms, and rejective subcategories correspond to surjective ring morphisms. We will study a chain of rejective subcategories to construct modules with endomorphisms rings of finite global dimension. We apply these result to study a function $r_\Lambda:\mod\Lambda\to\nnn_{\ge0}$ which is a natural extension of Auslander's representation dimension.}

\vskip1em
In the representation theory of artin algebras and orders, it often plays an important role to study certain classes of subcategories of the module category. Typical examples are given by subcategories induced by morphisms of rings (\S\XAC), and subcategories induced by cotilting modules (\S\XAD). These subcategories are functorially finite in the sense of Auslander-Smalo [AS1]. One object of this paper is to study functorially finite subcategories from the viewpoint of its resolution dimension (\S\XAA). The subcategories of resolution dimension zero is often called bireflective [St], and we shall show that they correspond to ring epimorphisms (\S\XAFA). We shall introduce a special class of bireflective subcategories called {\it rejective subcategories} (\S\XAE), which was well-known in the representation theory of orders and recently played a crucial role in the study of representation-finite orders [I1,2][Ru1,2]. They correspond to factor algebras of artin algebras, and overrings of orders (\S\XAFA) which are non-commutative analogy of the normalization in the commutative ring theory.

Another object of this paper is to study certain chains of rejective subcategories called {\it rejective chains} (\S\XBB), which give a method to construct rings of finite global dimension (\S\XBBB). Recently, rejective chains were applied to give positive answers to two open problems in [I3,4]. One is Solomon's conjecture on zeta functions of orders [S1,2], and another is the finiteness problem of the representation dimension of artin algebras [A1][Xi1] (see \S\XDAA). We shall formulate these construction of rejective chains by using a certain functor $\fff_{\cc}$ (\S\XBC). Typical examples are given by preprojective partition of Auslander-Smalo [AS2,3], and Bass chains of Drozd-Kirichenko-Roiter [DKR] and Hijikata-Nishida [HN] (\S\XBCD). It was first observed by Dlab-Ringel [DR4] that certain chains of subcategories are related to quasi-hereditary algebras, introduced by Cline-Parshall-Scott in the representation theory of Lie algebras and algebraic groups [CPS1,2]. In \S\XC, we shall study the relationship between rejective chains and quasi-hereditary algebras (\S\XCEA), and calculate the global dimension of rings with rejective chains (\S\XCC). We shall also relate neat algebras of Agoston-Dlab-Wakamatsu [ADW].

In \S\XD, we shall generalize the concept of the representation dimension of artin algebras to orders over complete discrete valuation rings. We can apply the results in previous sections to study the representation dimension of artin algebras and orders, since it equals to the resolution dimension of finite subcategories (\S\XDAB). The representation dimension of artin algebras was introduced by M. Auslander [A1] as a homological invariant to measure how far an artin algebra is from being representation-finite. We shall introduce a function $r_\Lambda$ for an artin algebra or order $\Lambda$, whose value at $\Lambda\oplus D\Lambda$ equals to the representation dimension (\S\XDA). Our $r_\Lambda$ would give us much more information, and would be quite natural concept in Auslander's philosophy, the homological approach to the representation theory. Although $\rdim\Lambda$ does not distinguish tame hereditary algebras and wild hereditary algebras (\S\XDAD), we shall show that the supremum of $r_\Lambda$ determines the representation type of hereditary algebras (\S\XDFC) as an application of Rouquier's result on exterior algebras (\S\XDF). Moreover, we shall show that the value of $r_\Lambda(\Lambda)$ is closely related to the reflexive-finiteness of $\Lambda$ (\S\XDGB). In \S\XE, we shall study $r_\Lambda$ under finite equivalences. We shall generalize the recent result of Xiangqian [X], any stable equivalence preserves the representation dimension, by using the relative Auslander-Reiten theory introduced by Auslander-Solberg [ASo].

Some results in this paper was announced in [I5] without proof.

\vskip.5em{\bf\XZA\ Notations }
In this paper, any module is assumed to be a left module. For a ring $\Lambda$, we denote by $J_\Lambda$ the Jacobson radical of $\Lambda$, and by $\Mod\Lambda$ (resp. $\mod\Lambda$, $\pr\Lambda$) the category of (resp. finitely generated, finitely generated projective) $\Lambda$-modules. Mainly we shall treat two kinds (1) and (2) of rings below, where it is well-known that their representation theory have many common aspect [A2][Y].

(1) Let $R$ be a commutative local artinian ring and $E$ the injective full of the simple $R$-module. An $R$-algebra $\Lambda$ is called an {\it artin $R$-algebra} if it is a finitely generated $R$-module. We have the duality $D:=\hom_R(\ ,E):\mod\Lambda\leftrightarrow\mod\Lambda^{op}$. Let $\rin\Lambda:=D\pr\Lambda^{op}$ be the category of injective $\Lambda$-modules.

(2)[CR][Re] Let $R$ be a complete discrete valuation ring with the quotient field of $K$. An $R$-algebra $\Lambda$ is called an {\it $R$-order} if $\Lambda\in\pr R$. Then $\Lambda$ forms a subring of a finite dimensional $K$-algebra $\widetilde{\Lambda}:=\Lambda\otimes_RK$. For an $R$-order $\Lambda$, a left $\Lambda$-module $X$ is called a {\it $\Lambda$-lattice} if $X\in\pr R$. We denote by $\lat\Lambda$ the category of $\Lambda$-lattices. We have a functor $\widetilde{(\ )}:=(\ )\otimes_RK:\lat\Lambda\to\mod\widetilde{\Lambda}$, and we have the duality $D:=\hom_R(\ ,R):\lat\Lambda\leftrightarrow\lat\Lambda^{op}$. We call $\rin\Lambda:=D\pr\Lambda^{op}$ the category of {\it injective} $\Lambda$-lattices.

Let $\Lambda$ and $\Gamma$ be $R$-orders and $\phi:\Lambda\to\Gamma$ a morphism of $R$-algebras. We call $\Gamma$ an {\it overring} of $\Lambda$ if $\Cok\phi$ is an $R$-torsion module. An overring $\Gamma$ of $\Lambda$ with $\Ker\phi=0$ is called an {\it overorder}. An overorder of $\Lambda$ can be regarded as a subring of $\widetilde{\Lambda}$ containing $\Lambda$, and an overring of $\Lambda$ can be regarded as a subring of $\widetilde{\Lambda}/I$ containing $(\Lambda+I)/I$ for an ideal $I$ of $\widetilde{\Lambda}$. An order is called {\it maximal} if it has no proper overorder. It is well-known that maximal order is hereditary.

\vskip.5em{\bf\XZB\ Notations }
Let $\cc$ be an additive category, $\cc(X,Y):=\hom_{\cc}(X,Y)$, and $fg$ the composition of $f\in\cc(X,Y)$ and $g\in\cc(Y,Z)$. Throughout this paper, {\it any subcategory is assumed to be full and closed under isomorphisms, direct products, direct sums and direct summands.} We denote by $J_{\cc}$ the Jacobson radical of $\cc$. For a collection $S$ of objects in $\cc$, we denote by $\add S$ the smallest subcategory of $\cc$ containing $S$. We call $X\in\cc$ an {\it additive generator} of $\cc$ if $\add X=\cc$. We denote by $[S]$ the ideal of $\cc$ consisting of morphisms which factor through some object in $S$. For an ideal $I$ of $\cc$, a factor category $\cc/I$ of $\cc$ is defined by $\ob(\cc/I)=\ob\cc$ and $\cc/I(X,Y):=\cc(X,Y)/I(X,Y)$ for any $X,Y\in\cc$. We call $\cc$ {\it Krull-Schmidt} if any object is isomorphic to a finite direct sum of objects whose endomorphism rings are local.  We denote by $\ind\cc$ the set of isoclasses of indecomposable objects in $\cc$.

A {\it $\cc$-module} is a contravariant additive functors from $\cc$ to the category of abelian groups. We denote by $\Mod\cc$ the category of $\cc$-modules, where $(\Mod\cc)(M,M^\prime)$ consists of the natural transformations from $M$ to $M^\prime$. Then $\Mod\cc$ forms an abelian category. By Yoneda's Lemma, $\cc(\ ,X)$ is a projective object in $\Mod\cc$. We call $M\in\Mod\cc$ {\it finitely presented} if there exists an exact sequence $\cc(\ ,Y)\to\cc(\ ,X)\to M\to 0$. We denote by $\mod\cc$ the category of finitely presented $\cc$-modules.

\vskip.5em{\bf\XZC\ Definition }Let $\Lambda$ be an artin algebra or order (\XZA). For simplicity, we often put $\dn{M}_\Lambda:=\mod\Lambda$ if $\Lambda$ is an artin algebra, and $\dn{M}_\Lambda:=\lat\Lambda$ if $\Lambda$ is an order. We call $\Lambda$ {\it representation-finite} if $\ind\dn{M}_\Lambda$ is a finite set. Now let $\cc$ be a subcategory of $\dn{M}_\Lambda$. We call $\cc$ {\it closed under submodules} (resp. {\it factor modules}) if $X\in\cc$ (resp. $Z\in\cc$) holds for any exact sequence $0\to X\to Y\to Z\to0$ in $\mod\Lambda$ with $Y\in\cc$ and $X,Z\in\dn{M}_\Lambda$. We call $\cc$ {\it closed under extensions} if $Y\in\cc$ holds for any exact sequence $0\to X\to Y\to Z\to0$ in $\mod\Lambda$ with $X,Z\in\cc$. We call $\cc$ {\it closed under images} (resp. {\it kernels, cokernels}) if $\Im f\in\cc$ (resp. $\Ker f\in\cc$, $\Cok f\in\cc$) holds for any $f\in\cc(X,Y)$.

\vskip.5em{\bf\XA\ Approximation and Rejective subcategories}

\vskip.5em{\bf\XAA\ Definition }
(1)[AS1] Let $\cc$ be an additive category and $\cc^\prime$ a subcategory of $\cc$. We call $f\in\cc(Y,X)$ a {\it right $\cc^\prime$-approximation} of $X$ if $Y\in\cc^\prime$ and $\cc(\ ,Y)\stackrel{\cdot f}{\to}\cc(\ ,X)\to0$ is exact on $\cc^\prime$, or equivalently, $\cc(\ ,Y)\stackrel{\cdot f}{\to}[\cc^\prime](\ ,X)\to0$ is exact on $\cc$ (\XZB). We call $\cc^\prime$ {\it contravariantly finite} if any $X\in\cc$ has a right $\cc^\prime$-approximation. Dually, a {\it left $\cc^\prime$-approximation} and a {\it covariantly finite} subcategory are defined. We call $\cc^\prime$ {\it functorially finite} if it is contravariantly and covariantly finite.

We call a subcategory $\cc^\prime$ of $\cc$ {\it finite} (resp. {\it cofinite}) if $\cc^\prime$ (resp. $\cc/[\cc^\prime]$) has an additive generator (\XZB). If $\Lambda$ is an artin algebra or order in a semisimple algebra, then any finite (resp. cofinite) subcategory of $\dn{M}_\Lambda$ is functorially finite [AS1].

(2) Assume that an additive category $\cc$ has kernels and cokernels, and $\cc^\prime$ is a contravariantly finite subcategory of $\cc$. Then any $X\in\cc$ has a {\it right $\cc^\prime$-resolution}, which is a complex $\cdots\to Y_2\stackrel{f_2}{\to}Y_1\stackrel{f_1}{\to}Y_0\stackrel{f_0}{\to}X$ in $\cc$ such that $Y_i\in\cc^\prime$ and $\cdots\to\cc(\ ,Y_2)\stackrel{\cdot f_2}{\to}\cc(\ ,Y_1)\stackrel{\cdot f_1}{\to}\cc(\ ,Y_0)\stackrel{\cdot f_0}{\to}\cc(\ ,X)\to0$ is exact on $\cc^\prime$. We denote by $\Omega^n_{\cc^\prime}X$ the kernel of $f_{n-1}$. We write $\resdim{\cc^\prime}{X}\le n$ if $X$ has a right $\cc^\prime$-resolution with $Y_{n+1}=0$. We call $\resdim{\cc^\prime}{\cc}:=\sup\{\resdim{\cc^\prime}{X}\ |\ X\in\cc\}$ the {\it right resolution dimension} of $\cc^\prime$.\footnote{When $\cc=\mod\Lambda$, Sikko [Si] denote $\resdim{\cc^\prime}{\cc}$ by $\gl(\cc^\prime,\Lambda)$. In this paper, we shall use the notation $\resdim{\cc^\prime}{X}$ in [ABu] to consider arbitrary additive category $\cc$.} Dually, a {\it left $\cc^\prime$-resolution}, $\Omega^n_{\cc^\prime{}^{op}}X$, $\resdim{\cc^\prime{}^{op}}{X}$ and $\resdim{\cc^\prime{}^{op}}{\cc^{op}}$ are defined.

When $\cc$ is Krull-Schmidt, we call a right (left) $\cc^\prime$-resolution {\it minimal} if any $f_i$ is in $J_{\cc}$ (\XZB). One can easily show that any $X\in\cc$ has a minimal right (left) $\cc^\prime$-resolution, which is unique up to isomorphisms of complexes.

\vskip.5em{\bf\XAAA\ Proposition }{\it
Let $\cc$ and $\cc^\prime$ be in \XAA(2).

(1) $\resdim{\cc^\prime}{X}=\pd_{\cc^\prime}\cc(\ ,X)$ and $\resdim{\cc^\prime{}^{op}}{X}=\pd_{\cc^\prime{}^{op}}\cc(X,\ )$ hold for any $X\in\cc$. Thus $0\le \gl(\mod\cc^\prime)-\resdim{\cc^\prime}{\cc}\le2$ and $0\le\gl(\mod\cc^\prime{}^{op})-\resdim{\cc^\prime{}^{op}}{\cc^{op}}\le2$ hold.

(2) If any $X\in\cc$ is a kernel (resp. cokernel) of some $f\in\cc^\prime(Y,Z)$, then $\resdim{\cc^\prime}{\cc}=\max\{\gl(\mod\cc^\prime)-2,0\}$ (resp. $\resdim{\cc^\prime{}^{op}}{\cc^{op}}=\max\{\gl(\mod\cc^\prime{}^{op})-2,0\}$).

(3) If $\cc=\dn{M}_\Lambda$ for an artin algebra or order $\Lambda$ and $\cc^\prime$ is a functorially finite subcategory of $\cc$, then $\gl(\mod\cc^\prime)=\gl(\mod\cc^\prime{}^{op})$.}

\vskip.5em{\sc Proof }
(1) Since any right $\cc^\prime$-resolution of $X$ corresponds to a projective resolution of $\cc(\ ,X)$ in $\mod\cc^\prime$, the former assertion follows. For any $M\in\mod\cc^\prime$, take a projective resoution $\cc^\prime(\ ,Y)\stackrel{\cdot f}{\to}\cc^\prime(\ ,X)\to M\to 0$. Then $\pd_{\cc^\prime}M\le\resdim{\cc^\prime}{\Ker f}+2$ holds.

(2) Since $0\to\cc^\prime(\ ,X)\to\cc^\prime(\ ,Y)\stackrel{\cdot f}{\to}\cc^\prime(\ ,Z)$ is exact, we have $\resdim{\cc^\prime}{X}=\pd_{\cc^\prime}\cc(\ ,X)\le\max\{\gl(\mod\cc^\prime)-2,0\}$. Thus the assertion follows from (1).

(3) We only show the artin algebra case since the order case is similar. We have a natural duality $\Mod\cc\leftrightarrow\Mod\cc^{op}$ defined by $(DM)(X):=D(M(X))$. Since $\cc$ forms a dualizing variety by [AS1;2,3], we have an induced duality $D:\mod\cc\leftrightarrow\mod\cc^{op}$. Thus the assertion follows.\rule{5pt}{10pt}

\vskip.5em{\bf\XAAB\ }The following is an immediate consequence [A1]. We shall use it in \XDAB\ again.

\vskip.5em{\bf Corollary }{\it
Let $\Lambda$ be an artin algebra or order, and $\cc$ a functorially finite subcategory of $\dn{M}_\Lambda$. Then $\resdim{\cc^{op}}{\dn{M}_\Lambda^{op}}\ge\max\{\gl(\mod\cc)-2,0\}\le\resdim{\cc}{\dn{M}_\Lambda}$ holds, where the left (resp. right) equality holds if $\Lambda\in\cc$ (resp. $D\Lambda\in\cc$).}

\vskip.5em{\bf\XAAC\ Theorem }{\it
Let $\Lambda$ be an artin algebra or order, $\cc$ a contravariantly finite subcategory of $\dn{M}_\Lambda$ such that $\Lambda\in\cc$ and $\underline{\underline{\dn{M}}}_\Lambda:=\dn{M}_\Lambda/[\cc]$. For any $M\in\mod\underline{\underline{\dn{M}}}_\Lambda$, take a projective resolution $0\to\dn{M}_\Lambda(\ ,Z)\stackrel{\cdot b}{\to}\dn{M}_\Lambda(\ ,Y)\stackrel{\cdot a}{\to}\dn{M}_\Lambda(\ ,X)\to M\to0$ of $M$ in $\mod\dn{M}_\Lambda$. Then a projective resolution of $M$ in $\mod\underline{\underline{\dn{M}}}_\Lambda$ is given by the sequence below.
\begin{eqnarray*}
\cdots\to\underline{\underline{\dn{M}}}_\Lambda(\ ,\Omega_{\cc}^2X)\to\underline{\underline{\dn{M}}}_\Lambda(\ ,\Omega_{\cc} Z)\to\underline{\underline{\dn{M}}}_\Lambda(\ ,\Omega_{\cc} Y)\to\underline{\underline{\dn{M}}}_\Lambda(\ ,\Omega_{\cc} X)\to\ \ \ \ \ \ \ \ \ \ &&\\
\underline{\underline{\dn{M}}}_\Lambda(\ ,Z)\stackrel{\cdot b}{\to}\underline{\underline{\dn{M}}}_\Lambda(\ ,Y)\stackrel{\cdot a}{\to}\underline{\underline{\dn{M}}}_\Lambda(\ ,X)\to M\to0&&\end{eqnarray*}}

\vskip-.5em{\sc Proof }
One can easily show that $\underline{\underline{\dn{M}}}_\Lambda(\ ,Z)\stackrel{\cdot b}{\to}\underline{\underline{\dn{M}}}_\Lambda(\ ,Y)\stackrel{\cdot a}{\to}\underline{\underline{\dn{M}}}_\Lambda(\ ,X)\to M\to0$ is exact (e.g. [I2;II.1.3(4)]). By $M(\cc)=0$ and $\Lambda\in\cc$, we can take the following commutative diagram, where $f$ is a right $\cc$-approximation of $X$. 
\[\begin{diag}
0&\RA{}&Z&\RA{b}&Y&\RA{a}&X&\RA{}&0\\
&&\uparrow&&\uparrow&&\parallel\\
0&\RA{}&\Omega_{\cc}X&\RA{}&W&\RA{f}&X&\RA{}&0
\end{diag}\]

Taking the mapping cone, we have an exact sequence $0\to\Omega_{\cc}X\to Z\oplus W\stackrel{{b\choose*}}{\to}Y\to0$. Thus $0\to\dn{M}_\Lambda(\ ,\Omega_{\cc}X)\to\dn{M}_\Lambda(\ ,Z\oplus W)\stackrel{\cdot{b\choose*}}{\to}\dn{M}_\Lambda(\ ,Y)\to M^\prime\to0$ gives a projective resolution of $M^\prime\in\mod\underline{\underline{\dn{M}}}_\Lambda$. The assertion follows inductively.\rule{5pt}{10pt}

\vskip.5em{\bf\XAAD\ }By \XAAB, $\resdim{\cc}{\dn{M}_\Lambda}$ gave the value of $\gl(\mod\cc)$ for the subcategory $\cc$. Now \XAAC\ implies that $\resdim{\cc}{\dn{M}_\Lambda}$ also gives an upper bound of the value of $\gl(\mod\underline{\underline{\dn{M}}}_\Lambda)$ for the factor category $\underline{\underline{\dn{M}}}_\Lambda$.

\vskip.5em{\bf Corollary }{\it
Let $\Lambda$ be an artin algebra or order, and $\cc$ a contravariantly finite subcategory of $\dn{M}_\Lambda$ such that $\Lambda\in\cc$. Then $\gl(\mod\underline{\underline{\dn{M}}}_\Lambda)\le 3(\resdim{\cc}{\dn{M}_\Lambda})-1$ holds.}

\vskip.5em{\bf\XAB\ Definition }Let $\cc$ be an additive category. Recall that a subcategory $\cc^\prime$ of $\cc$ is called {\it coreflective} (resp. {\it reflective}) if the inclusion functor $\cc^\prime\rightarrow\cc$ has a right (resp. left) adjoint [St][HS]. Then $\cc^\prime$ is a coreflective (resp. reflective) subcategory of $\cc$ if and only if $\resdim{\cc^\prime}{\cc}=0$ (resp. $\resdim{\cc^\prime{}^{op}}{\cc^{op}}=0$) holds. We often denote by $(\ )^-:\cc\to\cc^\prime$ (resp. $(\ )^+:\cc\to\cc^\prime$) the right (resp. left) adjoint functor of the inclusion functor $\cc^\prime\to\cc$, and by $\epsilon^-$ (resp. $\epsilon^+$) the counit (resp. unit). Then $0\to X^-\stackrel{\epsilon^-_X}{\to}X$ (resp. $X\stackrel{\epsilon^+_X}{\to}X^+\to0$) gives a right (resp. left) $\cc^\prime$-resolution of $X\in\cc$. We call $\cc^\prime$ {\it bireflective} if it is reflective and coreflective. 

\vskip.5em{\sc Proof }
Since `only if' part follows from \XABA(1) below, we shall show `if' part. Fix $a\in\cc(X_1,X_2)$. Let $0\to Y_i\stackrel{f_i}\to X_i$ be a right $\cc^\prime$-resolution for $i=1,2$. It is easily checked that the functor $(\ )^-:\cc\to\cc^\prime$ is defined by $X_i^-:=Y_i$ and $a^-\in\cc^\prime(Y_1,Y_2)$ is the unique morphism such that $a^-f_2=f_1a$. Since $\cc^\prime(\ ,Y_i)\stackrel{\cdot f_i}{\to}\cc(\ ,X_i)$ is an isomorphism on $\cc^\prime$, we obtained a right adjoint functor $(\ )^-$ of the inclusion functor $\cc^\prime\rightarrow\cc$.\rule{5pt}{10pt}

\vskip.5em{\bf\XABA\ Proposition }{\it
Let $\cc$ and $\dd$ be additive categories. Assume that $\fff:\cc\to\dd$ is a left adjoint of $\ggg:\dd\to\cc$ with the unit $\epsilon^+$ and the counit $\epsilon^-$.

(1)[AR2;1.2] $\dd^\prime:=\add\fff\cc$ is a contravariantly finite subcategory of $\dd$ and $\epsilon^-_Y:\fff\circ\ggg Y\to Y$ gives a right $\dd^\prime$-approximation of $Y\in\dd$. Thus $Y$ is contained in $\dd^\prime$ if and only if $\epsilon^-_Y$ is a split epimorphism. If $\fff$ is full, then $\dd^\prime$ is a coreflective subcategory of $\dd$. Dually, $\cc^\prime:=\add\ggg\dd$ is a cavariantly finite subcategory of $\cc$ and $\epsilon^+_X:X\to\ggg\circ\fff X$ gives a left $\cc^\prime$-approximation of $X\in\cc$. Thus $X$ is contained in $\cc^\prime$ if and only if $\epsilon^+_X$ is a split monomorphism. If $\ggg$ is full, then $\cc^\prime$ is a reflective subcategory of $\cc$.

(2) Let $\cc^\prime$ a subcategory of $\cc$ and $\dd^\prime$ a subcategory of $\dd$ such that $\fff\cc^\prime\subseteq\dd^\prime$ and $\ggg\dd^\prime\subseteq\cc^\prime$. Assume that $\epsilon^+$ is an isomorphism on $\cc^\prime$ and $\epsilon^-$ is an isomorphism on $\dd^\prime$.

\strut\kern1em(i) If $\cc^\prime$ is a contravariantly finite subcategory of $\cc$, then $\dd^\prime$ is that of $\dd$. For $Y\in\dd$, if $\cdots\stackrel{f_2}{\to}X_2\stackrel{f_1}{\to}X_1\stackrel{f_0}{\to}\ggg Y$ is a right $\cc^\prime$-resolution of $\ggg Y$, then $\cdots\stackrel{\fff f_2}{\longrightarrow}\fff X_1\stackrel{\fff f_1}{\longrightarrow}\fff X_0\stackrel{\fff(f_0)\epsilon^-_Y}{\longrightarrow}Y$ is a right $\dd^\prime$-resolution of $Y$. Thus $\resdim{\dd^\prime}{\dd}\le\resdim{\cc^\prime}{\cc}$ holds.

\strut\kern1em(ii) If $\dd^\prime$ is a covariantly finite subcategory of $\dd$, then $\cc^\prime$ is that of $\cc$. For $X\in\cc$, if $\fff X\stackrel{g_0}{\to}Y_0\stackrel{g_1}{\to}Y_1\stackrel{g_2}{\to}\cdots$ is a left $\dd^\prime$-resolution of $\fff X$, then $X\stackrel{\epsilon^+_X\ggg g_0}{\longrightarrow}\ggg Y_0\stackrel{\ggg g_1}{\longrightarrow}\ggg Y_1\stackrel{\ggg g_2}{\longrightarrow}\cdots$ is a left $\cc^\prime$-resolution of $X$. Thus $\resdim{\cc^\prime{}^{op}}{\cc^{op}}\le\resdim{\dd^\prime{}^{op}}{\dd^{op}}$ holds.}

\vskip.5em{\sc Proof }
Let $\eta_{X,Y}:\cc(X,\ggg Y)\to\dd(\fff X,Y)$ be a functorial isomorphism for $X\in\cc$ and $Y\in\dd$. Then the following diagram is commutative.
\[\begin{diag}
\cc(X,\ggg Y)&\times&\cc(\ggg Y,\ggg Y)&\RA{}&\cc(X,\ggg Y)\\
\downarrow^{\fff}&&\downarrow^{\eta_{\ggg Y,Y}}&&\downarrow^{\eta_{X,Y}}\\
\dd(\fff X,\fff\circ\ggg Y)&\times&\dd(\fff\circ\ggg Y,Y)&\RA{}&\dd(\fff X,Y)
\end{diag}\]

\strut\kern-.5em
(1)(ii) Since $\epsilon^-_Y=\eta_{\ggg Y,Y}(1_{\ggg Y})$, the map $\dd(\fff X,\fff\circ\ggg Y)\stackrel{\cdot\epsilon^-_Y}{\longrightarrow}\dd(\fff X,Y)$ is surjective. It is bijective if $\fff$ is full.

(2)(i) Take $a\in\dd(Z,Y)$ with $Z\in\dd^\prime$. Then there exists $b$ such that $\ggg a=bf_0$. Then $\fff\circ\ggg(a)=\fff(b)\fff(f_0)$ holds. Since the left diagram below commutative, $a$ factors through $\fff(f_0)\epsilon^-_Y$.
\[\begin{diag}
\fff X_0&\RA{\fff f_0}&\fff\circ\ggg Y&\RA{\epsilon^-_Y}&Y\\
&\ULA{\fff b}&\UA{\fff\circ\ggg a}&&\UA{a}\\
&&\fff\circ\ggg Z&\RA{\epsilon^-_Z}&Z
\end{diag}\ \ \ \ \ \ \ \ \ \ \begin{diag}
X_0&\RA{f_0}&\ggg Y\\
&\ULA{b}&\UA{\ggg a}\\
&&\ggg Z
\end{diag}\]

We notice here that $\fff(\epsilon^+_X)\epsilon^-_{\fff X}=1_{\fff X}$ holds by the following commutative diagram.
\[\begin{diag}
(\epsilon^+_X,1_{\ggg\circ\fff X})\in&\cc(X,\ggg\circ\fff X)&\times&\cc(\ggg\circ\fff X,\ggg\circ\fff X)&\RA{}&\cc(X,\ggg\circ\fff X)&\ni\epsilon^+_X\\
\downarrow&\downarrow^{\fff}&&\downarrow^{\eta_{\ggg\circ\fff X,\fff X}}&&\downarrow^{\eta_{X,\fff X}}&\downarrow\\
(\fff\epsilon^+_X,\epsilon^-_{\fff X})\in&\dd(\fff X,\fff\circ\ggg\circ\fff X)&\times&\dd(\fff\circ\ggg\circ\fff X,\fff X)&\RA{}&\dd(\fff X,\fff X)&\ni 1_{\fff X}
\end{diag}\]

For $n>0$, we will show that $\dd(\ ,\fff X_{n+1})\stackrel{\cdot\fff f_{n+1}}{\longrightarrow}\dd(\ ,\fff X_n)\stackrel{\cdot\fff f_n}{\longrightarrow}\dd(\ ,\fff X_{n-1})$ is exact on $\dd^\prime$. Take $a\in\dd(Z,\fff X_n)$ with $Z\in\dd^\prime$ and $a\fff(f_n)=0$. Then $\ggg(a)\ggg\circ\fff(f_n)=0$ implies $\ggg(a)(\epsilon^+_{X_n})^{-1}f_n=0$ holds. Thus there exists $b$ which makes the right diagram below commutative. Since the left diagram below commutative, $a$ factors through $\fff f_{n+1}$.
\[\begin{diag}
&&\fff X_n\\
&\URA{\fff f_{n+1}}&\DA{\fff\epsilon^+_{X_n}}&\DRA{1_{\fff X_n}}&&&\\
\fff X_{n+1}&&\fff\circ\ggg\circ\fff X_n&\RA{\epsilon^-_{\fff X_n}}&\fff X_n&\RA{\fff f_n}&\fff X_{n-1}\\
&\ULA{\fff b}&\UA{\fff\circ\ggg a}&&\UA{a}&\URA{0}\\
&&\fff\circ\ggg Z&\RA{\epsilon^-_Z}&Z
\end{diag}\ \ \ \ \ \ \ \ \ \ \begin{diag}
&&X_{n}&\RA{f_n}&X_{n-1}\\
&\URA{f_{n+1}}&\DA{\epsilon^+_{X_{n}}}&&\DA{\epsilon^+_{X_{n-1}}}\\
X_{n+1}&&\ggg\circ\fff X_n&\RA{\ggg\circ\fff f_n}&\ggg\circ\fff X_{n-1}\\
&\ULA{b}&\UA{\ggg a}&\URA{0}\\
&&\ggg Z
\end{diag}\]

We can show that $\dd(\ ,\fff X_1)\stackrel{\cdot\fff f_1}{\longrightarrow}\dd(\ ,\fff X_0)\stackrel{\cdot\fff(f_0)\epsilon^-_Y}{\longrightarrow}\dd(\ ,Y)$ is exact on $\dd^\prime$ by a quite similar argument, where we may replace the right diagram above by the following one.
\[\begin{diag}
X_{0}&\RA{f_0}&\ggg Y\\
\DA{\epsilon^+_{X_{0}}}&&\DA{\epsilon^+_{Y}}&\DRA{1_{\ggg Y}}\\
\ggg\circ\fff X_0&\RA{\ggg\circ\fff f_0}&\ggg\circ\fff\circ\ggg Y&\RA{\ggg\epsilon^-_Y}&\ggg Y\\
&\ULA{\ggg a}&&\URA{0}\\
&&\ggg Z&\ \ \ \rule{5pt}{10pt}
\end{diag}\]

\vskip.5em{\bf\XAC\ }
Let $\phi:\Lambda\to\Gamma$ be a morphism of rings. We denote by $\phi^*:\Mod\Gamma\to\Mod\Lambda$ the natural induced functor. Then $\phi^*$ has a left adjoint $\Gamma\otimes_\Lambda\ :\Mod\Lambda\to\Mod\Gamma$ with an unit $\epsilon^+$ and a counit $\delta^-$, which are defined by $\epsilon^+_X:=\phi\otimes 1:X\to\Gamma\otimes_\Lambda X$ for $X\in\Mod\Lambda$ and $\delta^-_Y:\Gamma\otimes_\Lambda Y\to Y$, $\delta^-_Y(\gamma\otimes y):=\gamma y$ for $Y\in\Mod\Gamma$. Dually, $\phi^*$ has a right adjoint $\hom_\Lambda(\Gamma,\ ):\Mod\Lambda\to\Mod\Gamma$ with a counit $\epsilon^-$ and an unit $\delta^+$, which are defined by $\epsilon^-_X:=(\phi\cdot):\hom_\Lambda(\Gamma,X)\to X$ for $X\in\Mod\Lambda$ and $\delta^+_Y:Y\to\hom_\Lambda(\Gamma,Y)$, $\delta^+_Y(y)(\gamma):=\gamma y$ for $Y\in\Mod\Gamma$. 

(1) Define a $(\Lambda,\Lambda)$-module $C_\phi$ by an exact seqeuence $\Lambda\stackrel{\phi}{\to}\Gamma\stackrel{a}{\to}C_\phi\to0$. We have the associated complexes below given by $d^+_i(x_1\otimes x_2\otimes\cdots\otimes x_i):=1\otimes a(x_1)\otimes x_2\otimes\cdots\otimes x_i$ and $d^-_i(x_1\otimes\cdots\otimes x_{i-1}\otimes x_i):=x_1\otimes\cdots\otimes x_{i-1}\otimes a(x_i)\otimes1$.
{\small\begin{eqnarray*}
{\bf X}_\phi^+:&&0\to\Lambda\stackrel{\phi}{\longrightarrow}\Gamma\stackrel{d^+_1}{\longrightarrow}\Gamma\otimes_\Lambda C_\phi\stackrel{d^+_2}{\longrightarrow}\Gamma\otimes_\Lambda C_\phi\otimes_\Lambda C_\phi\stackrel{d^+_3}{\longrightarrow}\Gamma\otimes_\Lambda C_\phi\otimes_\Lambda C_\phi\otimes_\Lambda C_\phi\to\cdots\\
{\bf X}_\phi^-:&&0\to\Lambda\stackrel{\phi}{\longrightarrow}\Gamma\stackrel{d^-_1}{\longrightarrow}C_\phi\otimes_\Lambda\Gamma\stackrel{d^-_2}{\longrightarrow}C_\phi\otimes_\Lambda C_\phi\otimes_\Lambda\Gamma\stackrel{d^-_3}{\longrightarrow}C_\phi\otimes_\Lambda C_\phi\otimes_\Lambda C_\phi\otimes_\Lambda\Gamma\to\cdots
\end{eqnarray*}}

\vskip-1em
(2) Define a $(\Gamma,\Gamma)$-module $D_\phi$ by an exact seqeuence $0\to D_\phi\stackrel{b}{\to}\Gamma\otimes_\Lambda\Gamma\stackrel{c}{\to}\Gamma\to0$, where $c$ is the multiplication map. We have the associated complexes below given by $e^+_i(x_0\otimes x_1\otimes\cdots\otimes x_i):=b(x_0x_1)\otimes\cdots\otimes x_i\in\Gamma\otimes_\Lambda\Gamma\otimes_\Gamma D_\phi^{\otimes i-1}=\Gamma\otimes_\Lambda D_\phi^{\otimes i-1}$ and $e^-_i(x_0\otimes\cdots\otimes x_{i-1}\otimes x_{i}):=x_0\otimes\cdots\otimes b(x_{i-1}x_{i})\in D_\phi^{\otimes i-1}\otimes_\Gamma\Gamma\otimes_\Lambda\Gamma= D_\phi^{\otimes i-1}\otimes_\Lambda\Gamma$.
{\small\begin{eqnarray*}
{\bf Y}_\phi^+:&&\cdots\to\Gamma\otimes_\Lambda D_\phi\otimes_\Gamma D_\phi\otimes_\Gamma D_\phi\stackrel{e^+_3}{\longrightarrow}\Gamma\otimes_\Lambda D_\phi\otimes_\Gamma D_\phi\stackrel{e^+_2}{\longrightarrow}\Gamma\otimes_\Lambda D_\phi\stackrel{e^+_1}{\longrightarrow}\Gamma\otimes_\Lambda\Gamma\stackrel{c}{\to}\Gamma\to0\\
{\bf Y}_\phi^-:&&\cdots\to D_\phi\otimes_\Gamma D_\phi\otimes_\Gamma D_\phi\otimes_\Lambda\Gamma\stackrel{e^-_3}{\longrightarrow} D_\phi\otimes_\Gamma D_\phi\otimes_\Lambda\Gamma\stackrel{e^-_2}{\longrightarrow} D_\phi\otimes_\Lambda\Gamma\stackrel{e^-_1}{\longrightarrow}\Gamma\otimes_\Lambda\Gamma\stackrel{c}{\to}\Gamma\to0
\end{eqnarray*}}

\vskip-1em{\bf\XACA\ Theorem }{\it
Let $\phi:\Lambda\to\Gamma$ be a morphism of rings. Then the assertions below hold, where we can replace $\Mod$ by $\mod$ if $\Lambda$ and $\Gamma$ are artin algebras.

(1) $\xx_\phi:=\add\phi^*(\Mod\Gamma)$ is a functorially finite subcategory of $\Mod\Lambda$. A left and right $\xx_\phi$-resolutions of $X\in\Mod\Lambda$ are given by ${\bf X}_\phi^+\otimes_\Lambda X$ and $\hom_\Lambda({\bf X}_\phi^-,X)$ below respectively, where $\Omega^n_{\xx_\phi^{op}}X=C_\phi^{\otimes n}\otimes_\Lambda X$ and $\Omega^n_{\xx_\phi}X=\hom_\Lambda(C_\phi^{\otimes n},X)$ hold.
{\small\begin{eqnarray*}
&X\stackrel{\epsilon^+_X}{\longrightarrow}\Gamma\otimes_\Lambda X\stackrel{d^+_1\otimes 1}{\longrightarrow}\Gamma\otimes_\Lambda C_\phi\otimes_\Lambda X\stackrel{d^+_2\otimes 1}{\longrightarrow}\Gamma\otimes_\Lambda C_\phi\otimes_\Lambda C_\phi\otimes_\Lambda X\to\cdots&\\
&\cdots\to\hom_\Lambda(C_\phi\otimes_\Lambda C_\phi\otimes_\Lambda\Gamma,X)\stackrel{d^-_2\cdot}{\longrightarrow}\hom_\Lambda(C_\phi\otimes_\Lambda\Gamma,X)\stackrel{d^-_1\cdot}{\longrightarrow}\hom_\Lambda(\Gamma,X)\stackrel{\epsilon^-_X}{\longrightarrow}X&
\end{eqnarray*}}

\vskip-1em
(2) $\yy_\phi^+:=\add\hom_\Lambda(\Gamma,\Mod\Lambda)$ is a covariantly finite subcategory of $\Mod\Gamma$. A left $\yy_\phi^+$-resolutions of $Y\in\Mod\Gamma$ is given by an exact sequence $\hom_\Gamma({\bf Y}_\phi^+,Y)$ below, where $\Omega^n_{\yy_\phi^+{}^{op}}Y=\hom_\Gamma( D_\phi^{\otimes n},Y)$ holds.
{\small\[0\to Y\stackrel{\delta^+_Y}{\longrightarrow}\hom_\Lambda(\Gamma,Y)\stackrel{e^+_1\cdot}{\longrightarrow}\hom_\Lambda( D_\phi,Y)\stackrel{e^+_2\cdot}{\longrightarrow}\hom_\Lambda( D_\phi\otimes_\Gamma D_\phi,Y)\stackrel{e^+_3\cdot}{\longrightarrow}\hom_\Lambda( D_\phi\otimes_\Gamma D_\phi\otimes_\Gamma D_\phi,Y)\to\cdots\]}

\vskip-1em
(3) $\yy_\phi^-:=\add\Gamma\otimes_\Lambda(\Mod\Lambda)$ is a contravariantly finite subcategory of $\Mod\Gamma$. A right $\yy_\phi^-$-resolutions of $Y\in\Mod\Gamma$ is given by an exact sequence ${\bf Y}_\phi^-\otimes_\Gamma Y$ below, where $\Omega^n_{\yy_\phi^-}Y= D_\phi^{\otimes n}\otimes_\Gamma Y$ holds.
{\small\[\cdots\to D_\phi\otimes_\Gamma D_\phi\otimes_\Gamma D_\phi\otimes_\Lambda Y\stackrel{e^-_3\otimes1}{\longrightarrow} D_\phi\otimes_\Gamma D_\phi\otimes_\Lambda Y\stackrel{e^-_2\otimes1}{\longrightarrow} D_\phi\otimes_\Lambda Y\stackrel{e^-_1\otimes1}{\longrightarrow}\Gamma\otimes_\Lambda Y\stackrel{\delta^-_Y}{\longrightarrow}Y\to0\]}

\vskip-1em
(4) For any $Y\in\Mod\Gamma$, $\phi^*(D_\phi^{\otimes n}\otimes_\Gamma Y)=C_\phi^{\otimes n}\otimes_\Lambda\phi^*Y$ and $\phi^*\hom_\Gamma(D_\phi^{\otimes n},Y)=\hom_\Lambda(C_\phi^{\otimes n},\phi^*Y)$ hold.}

\vskip.5em{\sc Proof }
(1)(2)(3) All assertions follow immediately from \XABA(1), where the exactness of the sequence in (2) (resp. (3)) follows from the fact that $c:\Gamma\otimes_\Lambda\Gamma\to\Gamma$ is a split epimorphism of left (resp. right) $\Gamma$-modules.

(4) Since $\epsilon^+_{\phi^*Y}\delta^-_Y=1$ holds, we obtain $\phi^*(D_\phi\otimes_\Gamma Y)=\phi^*\Ker\delta^-_Y=\Cok\epsilon^+_{\phi^*Y}=C_\phi\otimes_\Lambda\phi^*Y$. Thus the assertion follows inductively.\rule{5pt}{10pt}

\vskip.5em{\bf\XAD\ }Another well-known example of functorially finite subcategories is given by Auslander-Buchweitz theory [ABu] and its application to the cotilting theory [AR1]. Let $\Lambda$ be an artin algebra and $T$ a cotilting $\Lambda$-module with $\id_\Lambda T\le n$. Then the Auslander-Buchweitz theory implies that $\xx_T:=\{ X\in\mod\Lambda\ |\ \ext^i_\Lambda(X,T)=0$ for any $i>0\}$ is a contravariantly finite subcategory of $\mod\Lambda$. More precisely, we shall show the following theorem.

\vskip.5em{\bf\XADA\ Theorem }{\it
Let $\Lambda$, $T$, $\xx_T$ and $n$ be those in (1)--(3) below. Then $\resdim{\xx_T}{(\mod\Lambda)}=n$, $\resdim{\xx_T^{op}}{(\mod\Lambda^{op})}=\max\{n-2,0\}$ and $n\le \gl(\mod\xx_T)=\gl(\mod\xx_T^{op})\le\max\{n,2\}$ hold.

(1) $\Lambda$ is an artin algebra, $T$ is a cotilting $\Lambda$-module with $\id_\Lambda T=n$ and $\xx_T:=\{ X\in\mod\Lambda\ |\ \ext^i_\Lambda(X,T)=0$ for any $i>0\}$.

(2) $\Lambda$ is an order, $T:=D\Lambda$, $\xx_T:=\lat\Lambda$ and $n:=1$.

(3) $\Lambda$ is an $n$-dimensional commutative local noetherian ring with a dualizing module $T$ and $\xx_T$ is the category of maximal Cohen-Macaulay $\Lambda$-modules.}

\vskip.5em{\sc Proof }
We only prove (1) since essentially (2) and (3) are special cases of (1).

(i) If there is an exact sequence $0\to C_n\to X_{n-1}\stackrel{f_{n-1}}{\to}\cdots\to X_1\stackrel{f_1}{\to}X_0\to C_0\to0$ with $X_i\in\xx_T$, then $C_n\in\xx_T$.

Putting $C_i:=\Ker f_{i-1}$, we obtain an exact sequence $0\to C_{i}\to X_{i-1}\stackrel{}{\to}C_{i-1}\to0$. Since $\id_\Lambda T=n$ and $X_i\in\xx_T$, we obtain $\ext^i_\Lambda(C_n,T)=\ext^{i+1}_\Lambda(C_{n-1},T)=\cdots=\ext^{i+n}_\Lambda(C_0,T)=0$ for any $i>0$. Thus $C_n\in\xx_T$ holds.

(ii) For any $C\in\mod\Lambda$, take a $\xx_T$-resolution $0\to\Omega_{\xx_T}^nC\to X_{n-1}\to\cdots\to X_0\to C\to0$, which is exact by $\Lambda\in\xx_T$. Then $\Omega_{\xx_T}^nC\in\xx_T$ holds by (i). Thus $\resdim{\xx_T}{(\mod\Lambda)}\le n$ holds. If $\Omega_{\xx_T}^nC=0$ holds, then the exact sequence $0\to X_{n-1}\to\cdots\to X_0\to C\to0$ with $X_i\in\xx_T$ implies $\ext^n_\Lambda(C,T)=0$. Thus $\id_\Lambda T=n$ implies $\resdim{\xx_T}{(\mod\Lambda)}=n$.

(iii) We will show $\resdim{\xx_T^{op}}{(\mod\Lambda^{op})}\le m:=\max\{n-2,0\}$. Then our proof completes since $n\le \gl(\mod\xx_T)=\gl(\mod\xx_T^{op})\le\resdim{\xx_T^{op}}{(\mod\Lambda^{op})}+2$ holds by \XAAA(1)(3).

 Put $\Gamma:=\endm_\Lambda(T)$. Then $T$ is a cotilting $\Gamma^{op}$-module with $\id_{\Gamma}T=n$ [M]. Put $\xx_T^\prime:=\{ X\in\mod\Gamma^{op}\ |\ \ext^i_\Gamma(X,T)=0$ for any $i>0\}$. Then $\resdim{\xx_T^\prime{}}{(\mod\Gamma^{op})}=n$ holds by (ii). We have functors $\fff:=\hom_\Lambda(\ ,T):\mod\Lambda\to\mod\Gamma^{op}$ and $\ggg:=\hom_{\Gamma^{op}}(\ ,T):\mod\Gamma^{op}\to\mod\Lambda$, which induce equivalences between $\xx_T$ and $\xx_T^\prime$. Since we have a functorial isomorphism $\hom_\Lambda(X,\ggg DY)=\hom_{(\Lambda,\Gamma)}(X\otimes_{\zzz}DY,T)=\hom_{\Gamma^{op}}(DY,\hom_\Lambda(X,T))=\hom_\Gamma(D\fff X,Y)$, $D\circ\fff$ is a left adjoint of $\ggg\circ D$.

For any $X\in\mod\Lambda$, take a projective resolution $\Lambda^l\to\Lambda^k\to X\to0$. Then we have an exact sequence $0\to\fff X\to T^k\to T^l$ of $\Gamma^{op}$-module. Then $\Omega_{\xx_T^\prime{}}^{m}\fff X\in\xx_T^\prime$ and $\resdim{\xx_T^\prime}{\fff X}\le m$ hold by (i)(ii). Thus \XABA(2) implies $\resdim{\xx_T^{op}}{X}\le m$.\rule{5pt}{10pt}

\vskip.5em{\bf\XADB\ }Let $\Lambda$ be an artin algebra and $\xx_n:=\add\Omega^n(\mod\Lambda)$. Auslander-Reiten [AR2,3] has shown that $\xx_n$ is a functorially finite subcategory of $\mod\Lambda$ for any $n$. If $\xx_n$ is closed under extensions, then there exists a cotilting $\Lambda$-module $T$ with $\id_\Lambda T\le n$ such that $\xx_T=\xx_n$ [AR1]. Thus we immediately obtain the corollary below, where the first inequality was obtained by Sikko [Si]. We notice that $n$-Gorenstein algebras satisfy that $\xx_n$ is closed under extensions.

\vskip.5em{\bf Corollary }{\it
Let $\Lambda$ be an artin algebra and $\xx_n:=\add\Omega^n(\mod\Lambda)$. If $\xx_n$ is closed under extensions, then $\resdim{\xx_n}{(\mod\Lambda)}\le n$, $\resdim{\xx_n^{op}}{(\mod\Lambda^{op})}\le\max\{n-2,0\}$ and $\gl(\mod\xx_n)=\gl(\mod\xx_n^{op})\le\max\{n,2\}$ hold.}

\vskip.5em{\bf\XAE\ Definition }(1) A subcategory $\cc^\prime$ of an additive category $\cc$ is called {\it right} (resp. {\it left}) {\it rejective} [I2,4] if the following equivalent conditions are satisfied (cf. \XAB).

\strut\kern1em(i) Any $X\in\cc$ has a monic right (resp. epic left) $\cc^\prime$-approximation.

\strut\kern1em(ii) The inclusion functor $\cc^\prime\rightarrow\cc$ has a right (resp. left) adjoint with a counit $\epsilon^-$ (resp. unit $\epsilon^+$) such that $\epsilon^-_X$ is monic (resp. $\epsilon^+_X$ is epic) for any $X\in\cc$.

In this case, $\epsilon^-_X$ (resp. $\epsilon^+_X$) gives a monic right (resp. epic left) $\cc^\prime$-approximation of $X\in\cc$ by \XAB. We call $\cc^\prime$ {\it rejective} if it is left and right rejective. Any right (resp. left) rejective subcategories are coreflective (resp. reflective), but the converse does not hold in general (see \XAFA).

(2) Let $\cc$ be a Krull-Schmidt category. We call $\cc$ {\it semisimple} if $J_{\cc}=0$ holds. We call a subcategory $\cc^\prime$ of $\cc$ {\it cosemisimple} if $\cc/[\cc^\prime]$ is semisimple, namely, any non-invertible morphism in $\cc$ between indecomposable objects factor through an object in $\cc^\prime$.

\vskip.5em{\bf\XAEA\ }Let $\cc^\prime$ be a subcategory of a Krull-Schmidt category $\cc$. Then $\cc^\prime$ is cosemisimple right rejective if and only if, for any $X\in\ind\cc\backslash\ind\cc^\prime$, there exists a morphism $f\in\cc(Y,X)$ such that $Y\in\cc^\prime$ and $\cc(\ ,Y)\stackrel{\cdot f}{\to}J_{\cc}(\ ,X)$ is an isomorphism on $\cc$.

\vskip.5em{\sc Proof }
Notice that $\cc^\prime$ is cosemisimple if and only if $[\cc^\prime](\ ,X)=J_{\cc}(\ ,X)$ holds for any $X\in\ind\cc\backslash\ind\cc^\prime$. Thus `only if' part follows. Similarly, `if' part follows from $J_{\cc}(\ ,X)=\cc(\ ,Y)f\subseteq[\cc^\prime](\ ,X)\subseteq J_{\cc}(\ ,X)$.\rule{5pt}{10pt}

\vskip.5em{\bf\XAEB\ Proposition }{\it
Let $\Lambda$ be an artin algebra or order, and $\cc$ a subcategory of $\dn{M}_\Lambda$. Then $\cc$ is a right (resp. left) rejective subcategory of $\dn{M}_\Lambda$ if and only if $\cc$ is closed under factor modules (resp. submodules) in the sense of \XZC. In this case, if $\cc$ is covariantly (resp. contravariantly) finite, then $\resdim{\cc^{op}}{(\mod\Lambda^{op})}\le 1$ (resp. $\resdim{\cc}{(\mod\Lambda)}\le 1$).}

\vskip.5em{\sc Proof }
We shall show `if' part. For $X\in\dn{M}_\Lambda$, put $X^-:=\sum_{Y\in\cc,\ f\in\hom_\Lambda(Y,X)}f(Y)$. Since $X^-$ is a factor module of some module in $\cc$, we obtain $X^-\in\cc$. Thus the natural inclusion $X^-\to X$ is a monic right $\cc$-approximation of $X$. We shall show `only if' part. Let $f:X\to Y$ be a surjection with $X\in\cc$ and $Y\in\dn{M}_\Lambda$. Then $f$ factors through the monic right $\cc$-approximation $\epsilon^-_Y$ of $Y$. Thus $\epsilon^-_Y$ is bijective and $Y\in\cc$ holds. The latter assertion is immediate.\rule{5pt}{10pt}

\vskip.5em{\bf\XAEC\ Example }Let $\Lambda$ be an artin algebra and $(\tt,\ff)$ a torsion theory on $\mod\Lambda$ [Ha]. Then $\tt$ is a right rejective and $\ff$ is a left rejective subcategory of $\mod\Lambda$ by \XAEB. For example, for a classical cotilting $\Lambda$-module $T$, $\tt:=\{ X\in\mod\Lambda\ |\ \hom_\Lambda(X,T)=0\}$ is right rejective and $\resdim{\tt^{op}}{(\mod\Lambda^{op})}\le 1$, and $\ff:=\{ X\in\mod\Lambda\ |\ \ext^1_\Lambda(X,T)=0\}$ is left rejective and $\resdim{\ff}{(\mod\Lambda)}\le 1$ (cf. \XAD) [As].

\vskip.5em{\bf\XAF\ Definition }
(1) Let $\phi:\Lambda\to\Gamma$ be a morphism of rings. Recall that the following conditions are equivalent [St].

\strut\kern1em(i) $\phi$ is an epimorphism in the category of rings.

\strut\kern1em(ii) $\phi^*:\Mod\Gamma\to\Mod\Lambda$ is full.

\strut\kern1em(iii) $D_\phi=0$, i.e. the multiplication map $\Gamma\otimes_\Lambda\Gamma\to\Gamma$ is bijective.

In this case, we can identify $\Mod\Gamma$ with the subcategory $\xx_\phi$ of $\Mod\Lambda$. Moreover, $\yy_\phi^-$ is a right rejective subcategory of $\Mod\Gamma$ and $\yy_\phi^+$ is a left rejective subcategory of $\Mod\Gamma$ by \XACA(2)(3). Any surjective ring morphism is epic, but the converse does not hold in general. For example, the natural inclusion $\def\arraystretch{.5}\left(\begin{array}{cc}{\scriptstyle k}&{\scriptstyle k}\\ {\scriptstyle 0}&{\scriptstyle k}\end{array}\right)\to\ma_2(k)$ is an epimorphism.

(2) We call a morphism $\phi:\Lambda\to\Gamma$ of artin algebras {\it quasi-split} if $\mod\Lambda=\add\phi^*(\mod\Gamma)$.

\strut\kern1em(i) If $\phi$ is quasi-split, then it splits as a morphism of right (resp. left) $\Lambda$-modules.

\strut\kern1em(ii) If $\phi$ splits as a morphism of two-sided $\Lambda$-modules, then it is quasi-split.

(3) Let $\phi_i:\Lambda\to\Gamma_i$ ($i=1,2$) be ring morphisms. We identify $\phi_1$ and $\phi_2$ if there exists a ring isomorphism $\psi:\Gamma_1\to\Gamma_2$ such that $\phi_1=\psi\circ\psi_2$. Notice that there exists a ring morphism $\psi:\Gamma_1\to\Gamma_2$ such that $\phi_1=\psi\circ\psi_2$ if and only if there exists a functor $\fff:\Mod\Gamma_2\to\Mod\Gamma_1$ such that $\phi_1^*\circ\fff$ is isomorphic to $\phi_2^*$.

\vskip.5em{\sc Proof }
(1)(iii)$\Rightarrow$(ii) $\hom_\Lambda(\phi^*X,\phi^*Y)=\hom_\Gamma(\Gamma\otimes_\Lambda\phi^*X,Y)=\hom_\Gamma(\Gamma\otimes_\Lambda\Gamma\otimes_\Gamma X,Y)=\hom_\Gamma(X,Y)$.

(ii)$\Rightarrow$(i) Let $\psi_i:\Gamma\to\Delta$ ($i=1,2$) be ring morphisms such that $\psi_1\circ\phi=\psi_2\circ\phi$. Then the identity map on $\Delta$ is contained in $\hom_\Lambda(\phi^*(\psi_1^*\Delta),\phi^*(\psi_2^*\Delta))$. Since this equals to $\hom_\Gamma(\psi_1^*\Delta,\psi_2^*\Delta)$, we obtain $\psi_1=\psi_2$.

(i)$\Rightarrow$(iii) Put $\Delta:=\def\arraystretch{.5}\left(\begin{array}{cc}{\scriptstyle\Gamma}&{\scriptstyle\Gamma\otimes_\Lambda\Gamma}\\ {\scriptstyle0}&{\scriptstyle\Gamma}\end{array}\right)$. Let $\psi_i:\Gamma\to\Delta$ ($i=1,2$) be the morphism of rings defined by $\psi_1(x)=\def\arraystretch{.5}\left(\begin{array}{cc}{\scriptstyle x}&{\scriptstyle 0}\\ {\scriptstyle0}&{\scriptstyle x}\end{array}\right)$ and $\psi_2(x)=\def\arraystretch{.5}\left(\begin{array}{cc}{\scriptstyle x}&{\scriptstyle x\otimes1-1\otimes x}\\ {\scriptstyle0}&{\scriptstyle x}\end{array}\right)$. Since $\psi_1\circ\phi=\psi_2\circ\phi$ holds, we obtain $\psi_1=\psi_2$. Thus $x\otimes1=1\otimes x$ holds for any $x\in\Gamma$, so we obtain $\Gamma\otimes_\Lambda\Gamma=\Gamma$.

(2)(i) We only have to put $X=\Lambda$ (resp. $X=D\Lambda$) in \XABA(1).

(ii) Assume that $a\in\hom_{(\Lambda,\Lambda)}(\Gamma,\Lambda)$ satisfies $\phi a=1$. Then $\epsilon^+_X(a\otimes 1_X)=1_X$ holds for any $X\in\mod\Lambda$. Thus $\epsilon^+_X$ is a split monomorphism, and $X\in\add\phi^*(\mod\Gamma)$ holds.

(3) We only have to construct $\psi$ from $\fff$.

(i) Assume that $\Gamma_i=\Gamma$ ($i=1,2$) and $\phi_1^*$ is isomorphic to $\phi_2^*$. Then there exists an inner automorphism $\psi$ of $\Gamma$ such that $\phi_1=\psi\circ\phi_2$.

Let $\alpha:\phi_1^*\to\phi_2^*$ be an isomorphism. Then $f:=\alpha_\Gamma\in\hom_\Lambda(\phi_1^*\Gamma,\phi_2^*\Gamma)$ makes the following diagram commutative for any $\gamma\in\Gamma$.
\[\begin{diag}
\phi_1^*\Gamma&\RA{f}&\phi_2^*\Gamma\\
\downarrow^{\cdot\gamma}&&\downarrow^{\cdot\gamma}\\
\phi_1^*\Gamma&\RA{f}&\phi_2^*\Gamma
\end{diag}\]

Since $f(\gamma^\prime\gamma)=f(\gamma^\prime)\gamma$ holds for any $\gamma^\prime\in\Gamma$, $f\in\endm_{\Gamma^{op}}(\Gamma)$. Thus there exists $a\in\Gamma$ such that $f=(a\cdot)$. Then $a$ is a unit of $\Gamma$ such that $a\phi_1(\lambda)=\phi_2(\lambda)a$ holds for any $\lambda\in\Lambda$.

(ii) We shall show the assertion. Since $\fff$ is an exact functor, we can take a $(\Gamma_2,\Gamma_1)$-bimodule $P$ such that $P\in\pr\Gamma_2$ and $\fff=\hom_{\Gamma_2}(P,\ )$. Since the identity functor is isomorphic to $\hom_{\Gamma_2}(P,\ )$ as a functor from $\Mod\Gamma_2$ to the category of abelian groups, $P$ is isomorphic to $\Gamma_2$ as a $\Gamma_2$-module by Yoneda's lemma. Thus the right action to $P$ induces a ring morphism $\psi:\Gamma_1\to\endm_{\Gamma_2}(P)=\Gamma_2$ such that $\fff$ is isomorphic to $\psi$. Then $\phi_1^*$ is isomorphic to $(\psi\circ\phi_2)^*$. Thus the assertion follows from (i).\rule{5pt}{10pt}

\vskip.5em{\bf\XAFA\ Theorem }{\it
Let $\Lambda$ be an artin algebra and $\cc$ a subcategory of $\mod\Lambda$. 

(1) The conditions below are equivalent, and there exists a bijection between rejective subcategories of $\mod\Lambda$ and factor algebras of $\Lambda$.

\strut\kern1em(i) $\cc$ is a rejective subcategory of $\mod\Lambda$.

\strut\kern1em(ii) $\cc=\xx_\phi:=\phi^*(\mod\Gamma)$ for a surjective ring morphism $\phi:\Lambda\to\Gamma$.

\strut\kern1em(iii) $\cc$ is closed under submodules and factor modules.

(2) The conditions below are equivalent, and there exists a bijection between bireflective subcategories of $\mod\Lambda$ and ring epimorphisms from $\Lambda$ to artin algebras.

\strut\kern1em(i) $\cc$ is a bireflective subcategory of $\mod\Lambda$.

\strut\kern1em(ii) $\cc=\xx_\phi:=\phi^*(\mod\Gamma)$ for a ring epimorphism $\phi:\Lambda\to\Gamma$ between artin algebras.

\strut\kern1em(iii) $\cc$ is functorially finite, and closed under kernels and cokernels.}

\vskip.5em{\sc Proof }
By \XAF(3), we only have to show the equivalences of conditions.

(1)(ii)$\Rightarrow$(iii) Put $\Gamma=\Lambda/I$ for an ideal $I$ of $\Lambda$. Then $X\in\mod\Lambda$ is contained in $\xx_\phi$ if and only if $IX=0$. Thus the assertion follows.

(iii)$\Rightarrow$(i) Immedaite from \XAEB.

(i)$\Rightarrow$(ii) Take $a\in \Lambda^+$ such that $(\cdot a)=\epsilon^+_\Lambda:\Lambda\to\Lambda^+$. Put $\Gamma:=\endm_{\Lambda}(\Lambda^+)$. Taking $\hom_\Lambda(\ ,\Lambda^+)$, we obtain a bijection $(a\cdot):\endm_\Lambda(\Lambda^+)=\Gamma\rightarrow\hom_{\Lambda}(\Lambda,\Lambda^+)=\Lambda^+$. Thus a map $\phi:\Lambda\rightarrow\Gamma$ is well-defined by $xa=a\phi(x)$ for any $x\in\Lambda$. Obviously $\phi$ is a ring morphism. Since $(a\cdot)$ is a bijection such that $(a\cdot)\circ \phi=(\cdot a)=\epsilon^+_\Lambda$, we can replace the left approximation $\Lambda\stackrel{\epsilon^+_\Lambda}{\to}\Lambda^+$ of $\Lambda$ by $\phi:\Lambda\to\Gamma$.

Since $\cc$ is left rejective, $\phi=\epsilon^+_\Lambda$ is surjective. For any $X\in\xx_\phi$, take a surjection $p:\Gamma^n\to X$. Then $p$ factors through the monic right $\cc$-approximation $\epsilon^-_X$. Thus $\epsilon^-_X$ is bijective and $X\in\cc$ holds. For any $X\in\cc$, take a surjection $p:\Lambda^n\to X$. Then $p^+:\Gamma^n\to X$ is also a surjection, and $X\in\xx_\phi$ holds. Thus $\cc=\xx_\phi$.

(2)(ii)$\Rightarrow$(i) Immedite from \XAC.

(i)$\Rightarrow$(iii) Let $0\to X\stackrel{a}{\to}Y\stackrel{b}{\to}Z$ be an exact sequence in $\mod\Lambda$ with $Y,Z\in\cc$. Then there exists $a^\prime$ such that $a=\epsilon^+_Xa^\prime$. Since $\epsilon^+_Xa^\prime b=ab=0$ holds and $Z\in\cc$, we obtain $a^\prime b=0$. Thus there exists $a^{\prime\prime}$ such that $a^\prime=a^{\prime\prime}a$. Now $a=\epsilon^+_Xa^{\prime\prime}a$ implies that $\epsilon^+_X$ is a split monomorphism, so $X\in\cc$.

(iii)$\Rightarrow$(i) Take a minimal right $\cc$-approximation $f:Y\to X$. Let $0\to Z\stackrel{g}{\to}Y\stackrel{f}{\to}X$ be an exact sequence in $\mod\Lambda$. We only have to show that $\hom_\Lambda(\ ,Y)=0$ holds on $\cc$. Take $a\in\hom_\Lambda(M,Z)$ with $M\in\cc$ and consider the following commutative diagram of exact sequences.
\[\begin{diag}
0&\RA{}&Z&\RA{g}&Y&\RA{f}&X\\
&&\uparrow^{a}&&\parallel&&\uparrow^{c}\\
&&M&\RA{ag}&Y&\RA{b}&L&\RA{}&0
\end{diag}\]

Since $L\in\cc$ holds by our assumption, there exists $c^\prime$ such that $c=c^\prime f$. Thus $f=bc^\prime f$ holds. Since $f$ is minimal, $bc^\prime$ is an automorphism of $Y$. Thus $b$ is an isomorphism, and $a=0$ holds.

(i)+(iii)$\Rightarrow$(ii) By the argument in the proof of (1)(i)$\Rightarrow$(ii), $\epsilon^+_\Lambda:\Lambda\to\Lambda^+$ is given by a ring morphism $\phi:\Lambda\to\Gamma$. For any $X\in\mod\Gamma$, take a projective resolution $\Gamma^m\to\Gamma^n\to X\to0$. Then $\phi^*X\in\cc$ and $\xx_\phi\subseteq\cc$ holds by (iii). For another $Y\in\mod\Gamma$, $\hom_\Lambda(\Lambda,\phi^*Y)=\phi^*Y\stackrel{\phi\cdot}{\to}\hom_\Lambda(\Gamma,\phi^*Y)$ is an isomorphism by $\phi^*Y\in\cc$. Thus we obtain the following commutative diagram of exact seqeunces.
\[\begin{diag}
0&\RA{}&\hom_\Gamma(X,Y)&\RA{}&Y^n&\RA{}&Y^m\\
&&\downarrow^{\phi^*}&&\downarrow&&\downarrow\\
0&\RA{}&\hom_\Lambda(\phi^*X,\phi^*Y)&\RA{}&\hom_\Lambda(\Gamma,\phi^*Y)^n&\RA{}&\hom_\Lambda(\Gamma,\phi^*Y)^m
\end{diag}\]

Since the middle and right maps are isomorphisms, so is the left one. Thus $\phi^*:\mod\Gamma\to\mod\Lambda$ is full, and $\phi$ is a ring epimorphism.

For any $X\in\cc$, take a surjection $p:\Lambda^n\to X$. Then $p^+:\Gamma^n\to X$ is also surjection. Applying same argument to $\Ker p^+\in\cc$, we obtain an exact sequence $\Gamma^m\stackrel{a}{\to}\Gamma^n\to X\to0$ for $a\in\hom_\Lambda(\Gamma^m,\Gamma^n)=\hom_\Gamma(\Gamma^m,\Gamma^n)$. Thus $\phi^*Y=X$ holds for the cokernel $Y$ of $a$ in $\mod\Gamma$. Hence $\cc=\xx_\phi$ holds.\rule{5pt}{10pt}

\vskip.5em{\bf\XAFB\ Corollary }{\it
Let $\phi:\Lambda\to\Gamma$ be a morphism of artin algebras and $\xx_\phi:=\add\phi^*(\mod\Gamma)$. Then $\xx_\phi$ is a bireflective subcategory of $\mod\Lambda$ if and only if $\phi=\phi_2\circ\phi_1$ holds for a ring epimorphism $\phi_1:\Lambda\to\Delta$ and a quasi-split morphism $\phi_2:\Delta\to\Gamma$ between artin algebras.}

\vskip.5em{\sc Proof }
We only have to show `if' part. By \XAFA(2), there exists a ring epimorphism $\phi_1:\Lambda\to\Delta$ with an artin algebra $\Delta$ such that $\xx_{\phi_1}=\xx_\phi$. Since $\phi_1^*:\mod\Delta\to\mod\Lambda$ is full faithful, there exists a functor $\fff:\mod\Gamma\to\mod\Delta$ such that $\phi^*$ is isomorphic to $\phi_1^*\circ\fff$. Then the assertion follows from \XAF(3).\rule{5pt}{10pt}

\vskip.5em{\bf\XAFC\ }We notice that \XAFA(1) holds even if we replace $\mod$ by $\Mod$. While the corresponding result of \XAFA(2) for $\Mod\Lambda$ is the theorem below, where $\Gamma$ in (ii) below is not necessarily an artin algebra. For example, the natural inclusion $\def\arraystretch{.5}\left(\begin{array}{cc}{\scriptstyle k}&{\scriptstyle k\oplus kx}\\ {\scriptstyle 0}&{\scriptstyle k}\end{array}\right)\to\ma_2(k[x])$ is an epimorphism. We shall omit the proof since it is shown by a parallel argument. The equivalence of (ii) and (iii) below were given in [GD].

\vskip.5em{\bf Theorem }{\it
Let $\Lambda$ be an artin algebra and $\cc$ a subcategory of $\Mod\Lambda$. The conditions below are equivalent, and there exists a bijection between bireflective subcategories of $\Mod\Lambda$ and ring epimorphisms from $\Lambda$.

\strut\kern1em(i) $\cc$ is a bireflective subcategory of $\Mod\Lambda$.

\strut\kern1em(ii) $\cc=\phi^*(\Mod\Gamma)$ for a ring epimorphism $\phi:\Lambda\to\Gamma$.

\strut\kern1em(iii) $\cc$ is closed under kernels and cokernels.}

\vskip.5em{\bf\XAFD\ }Let us consider an analogy of \XAFA\ for orders. For an overring $\Gamma$ of $\Lambda$ (\S\XZA), it is easily checked that $\phi^*:\mod\Gamma\to\mod\Lambda$ induces a full faithful functor $\phi^*:\lat\Gamma\rightarrow\lat\Lambda$. Then $\phi^*$ has a right adjoint dunctor $(\ )^-:=\hom_\Lambda(\Gamma,\ ):\lat\Lambda\to\lat\Gamma$ and a left adjoint functor $(\ )^+:=\Gamma\overline{\otimes}_\Lambda\ :\lat\Lambda\to\lat\Gamma$, where $\Gamma\overline{\otimes}_\Lambda X$ is a factor of $\Gamma\otimes_\Lambda X$ by its torsion submodule. We can prove the following theorem similarly (cf. [I2;II.5.3]).

\vskip.5em{\bf Theorem }{\it Let $\Lambda$ be an order and $\cc$ a subcategory of $\lat\Lambda$. Then the conditions below are equivalent, and there exists a bijection between rejective subcategories of $\lat\Lambda$ and overrings of $\Lambda$.

(i) $\cc$ is a rejective subcategory of $\lat\Lambda$.

(ii) $\cc=\phi^*(\lat\Gamma)$ for an overring $\phi:\Lambda\to\Gamma$ of $\Lambda$.

(iii) $\cc$ is closed under submodules and factor modules (\S\XZB).}

\vskip.5em{\bf\XAFE\ }
Let $\cc$ be a Krull-Schmidt category. A subset of $\ind\cc$ is called {\it rejectable} if it has the form $S=\ind\cc\backslash\ind\cc^\prime$ for some rejective subcategory $\cc^\prime$ of $\cc$. It is a quite interesting problem to study rejectable subsets of $\dn{M}_\Lambda$ for an artin algebra or order $\Lambda$.  There is a close relationship with Auslander-Reiten theory, and one can characterize finite rejectable subsets of $\dn{M}_\Lambda$ in terms of Auslander-Reiten quivers [I1,2]. The simplest case is the lemma below of Drozd-Kirichenko [DK][HN;2.2.1,2.2.2], which characterizes one-point rejectable subsets of $\dn{M}_\Lambda$.

\vskip.5em{\bf Proposition }{\it
Let $\Lambda$ be an artin algebra or order, and $X\in\ind\dn{M}_\Lambda$. Then $\{ X\}$ is rejectable if and only if $X\in\pr\Lambda\cap\rin\Lambda$. Then the corresponding rejective subcategory of $\dn{M}_\Lambda$ is cosemisimple except the case $\Lambda$ is an order such that $\Lambda=\Gamma\times\Delta$ for some maximal order $\Delta$ with $\ind(\lat\Delta)=\{ X\}$.}

\vskip.5em{\sc Proof }
Define a subcategory $\cc$ of $\dn{M}_\Lambda$ by $\ind\cc=\ind\dn{M}_\Lambda\backslash\{ X\}$. It is easily shown that $\cc$ is closed under submodules (resp. factor modules) if and only if $X\in\rin\Lambda$ (resp. $X\in\pr\Lambda$). Thus the equivalence follows from \XAFA(1) and \XAFD. If $J_\Lambda X\in\cc$ holds, then $\cc$ is cosemisimple since any non-invertible endomorphism of $X$ factors through $J_\Lambda X$. If $J_\Lambda X\notin\cc$, then $\Lambda$ is an order and $J_\Lambda X$ is isomorphic to $X$. This implies that $\Lambda=\Gamma\times\Delta$ for some maximal order $\Delta$ with $\ind(\lat\Delta)=\{ X\}$. See the proof (ii) of \XBCA(2).\rule{5pt}{10pt}

\vskip.5em{\bf\XAG\ }An extremely interesting example of functorially finite subcategories is given by radical embeddings [EHIS], which also appeared in [N].

\vskip.5em{\bf Theorem }{\it
Let $\Lambda\stackrel{\phi}{\subset}\Gamma$ be artin algebras with $J_\Lambda=J_\Gamma$.

(1) $\xx_\phi:=\add\phi^*(\mod\Gamma)$ satisfies $\resdim{\xx_\phi}{(\mod\Lambda)}\le 1$ and $\resdim{\xx_\phi^{op}}{(\mod\Lambda^{op})}\le 1$.

(2) $\yy^+_\phi:=\add\hom_\Lambda(\Gamma,\mod\Lambda)$ and $\yy^-_\phi:=\add\Gamma\otimes_\Lambda(\mod\Lambda)$ coincide with $\mod\Gamma$.

(3) $\phi^*$ induces a full faithful functor $\cc^\prime:=\mod\Gamma/[\mod\Gamma/J_\Gamma]\to\cc:=\mod\Lambda/[\mod\Lambda/J_\Lambda]$, and $\cc^\prime$ forms a rejective subcategory of $\cc$.}

\vskip.5em{\sc Proof }
We shall use the notations in \XAC. We only show the right-hand side assertion.

(1) Since $C_\phi J_\Lambda=0$ holds, $J_\Lambda\hom_\Lambda(C_\phi,X)=0$ holds for any $X\in\mod\Lambda$. Hence $\hom_\Lambda(C_\phi,X)\in\add_\Lambda(\Lambda/J_\Lambda)\subseteq\add_\Lambda(\Gamma/J_\Gamma)\subseteq\xx_\phi$ implies $\resdim{\xx_\phi}{(\mod\Lambda)}\le 1$ by \XACA(1).

(2)(3) A $\Gamma$-module is semisimple if and only if it is semisimple as a $\Lambda$-module. For any $X\in\mod\Lambda$, we have the following commutative diagram of exact sequences, where the upper sequence is a $\xx_\phi$-resolution and we put $X^-:=\hom_\Lambda(\Gamma,X)$.
\[\begin{diag}
&0&\RA{}&\hom_\Lambda(C_\phi,X)&\RA{}&X^-&\RA{\epsilon^-_X}&X\\
&&&\parallel&&\cup^b&&\cup^a\\
(*)\ \ \ &0&\RA{}&\hom_\Lambda(C_\phi,X)&\RA{}&\soc X^-&\RA{}&\soc X
\end{diag}\]

Since $a$ factors through $\epsilon^-_X$ by $\soc X\in\xx_\phi$, the sequence $(*)$ is split exact, and the diagram is pull-back. We shall show (2). By \XACA, we have an exact sequence $0\to Y\stackrel{\delta^+_Y}{\to}Y^-\to\hom_\Gamma(D_\phi,Y)\to0$ for any $Y\in\mod\Gamma$. The induced sequence $0\to\soc Y\stackrel{\delta^+_Y}{\to}\soc Y^-\to\hom_\Gamma(D_\phi,Y)$ should be split exact since $\length_\Lambda\phi^*\hom_\Gamma(D_\phi,Y)=\length_\Lambda\hom_\Lambda(C_\phi,\phi^*Y)=\length_\Lambda\soc Y^--\length_\Lambda\soc Y$ holds by $(*)$ and \XACA(4). Thus $\delta^+_Y:Y\to Y^-$ is also a split monomorphism.

We shall show (3). Fix any $X\in\mod\Lambda$ and $Y\in\mod\Gamma$. Assume that $g\in\hom_\Gamma(Y,X^-)$ corresponds to $f\in\hom_\Lambda(\phi^*Y,X)$. Then $f$ factors through $\mod\Lambda/J_\Lambda$ if and only if $f$ factors through $a$ if and only if $g$ factors through $b$ if and only if $g$ factors through $\mod\Gamma/J_\Gamma$. This means that $\epsilon^-_X$ induces an isomorphism $\cc^\prime(Y,X^-)\to\cc(Y,X)$. Thus the functor $\phi^*:\cc^\prime\to\cc$ and its right adjoint $(\ )^-:\cc\to\cc^\prime$ are well-defined. If $X\in\mod\Gamma$, then we have isomorphisms $\cc^\prime(Y,X)\stackrel{\delta^+_X}{\to}\cc^\prime(Y,X^-)\stackrel{\epsilon^-_X}{\to}\cc(Y,X)$ by (2), so the functor $\phi^*:\cc^\prime\to\cc$ is full faithful. To show that $\epsilon^-_X$ is monic in $\cc$, assume $h\epsilon^-_X=0$ in $\cc$. Then there exists $h^\prime$ such that $h\epsilon^-_X=h^\prime a$. Since above diagram is pull-back, $h$ factors through $b$. Thus $\cc^\prime$ is a right rejective subcategory of $\cc$.\rule{5pt}{10pt}

\vskip.5em{\bf\XB\ Rejective chain }

In this section, we will study a chain of (left, right) rejective subcategories. We will collect general facts on (left, right) rejective subcategories. 

\vskip.5em{\bf\XBA\ }Let $\cc$ be an additive category and $\cc^{\prime\prime}\subseteq\cc^\prime$ subcategories of $\cc$.

(1) If $\cc^{\prime\prime}$ is a (left, right) rejective subcategory of $\cc$, then it is a (left, right) rejective subcategory of $\cc^\prime$.

(2) If $\cc^\prime$ is a (left, right) rejective subcategory of $\cc$, then $\cc^\prime/[\cc^{\prime\prime}]$ is a (left, right) rejective subcategory of $\cc/[\cc^{\prime\prime}]$.

(3) If $\cc^{\prime\prime}$ is a rejective subcategory of $\cc^\prime$ and $\cc^\prime$ is a rejective subcategory of $\cc$, then $\cc^{\prime\prime}$ is a rejective subcategory of $\cc$.

(4) Assume that $\cc^{\prime\prime}$ is a right (resp. left) rejective subcategory of $\cc^\prime$ with a counit $\epsilon^\prime{}^-$ (resp. unit $\epsilon^\prime{}^+$) such that $\epsilon^\prime{}^-_X$ (resp. $\epsilon^\prime{}^+_X$) is monic (resp. epic) in $\cc$ for any $X\in\cc^\prime$. If $\cc^\prime$ is a right (resp. left) rejective subcategory of $\cc$, then $\cc^{\prime\prime}$ is a right (resp. left) rejective subcategory of $\cc$.

(5) Assume that $\cc^\prime$ is a semisimple functorially finite subcategory of $\cc$. Then $\cc^\prime$ is a right rejective subcategory of $\cc$ if and only if $\cc^\prime$ is a left rejective subcategory of $\cc$.

(6) Assume that both of $\cc^\prime$ and $\cc^{\prime\prime}$ are (left, right) rejective subcategory of $\cc$, and $\cc^{\prime\prime}$ is a cosemisimple subcategory of $\cc^\prime$. Then any subcategory $\dd$ with $\cc^{\prime\prime}\subseteq\dd\subseteq\cc^\prime$ is a (left, right) rejective subcategory of $\cc$.

\vskip.5em{\sc Proof }(1) Immediate. 

(2) Since $\cc(\ ,X^-)\stackrel{\cdot\epsilon^-_X}{\longrightarrow}[\cc^\prime](\ ,X)$ is an isomorphism for any $X\in\cc$, so is $[\cc^{\prime\prime}](\ ,X^-)\stackrel{\cdot\epsilon^-_X}{\longrightarrow}[\cc^{\prime\prime}](\ ,X)$. Thus so is $\cc/[\cc^{\prime\prime}](\ ,X^-)\stackrel{\cdot\epsilon^-_X}{\longrightarrow}\cc/[\cc^{\prime\prime}](\ ,X)$.

(4) Immediate since $\epsilon^\prime{}^-_{X^-}\epsilon^-_X$ gives a monic right $\cc^{\prime\prime}$-approximation of $X$.

(3) By (4), we only have to show that any monomorphism $a\in\cc^\prime(X,Y)$ in $\cc^\prime$ is still monic in $\cc$. Assume that $b\in\cc(Z,X)$ satisfies $ba=0$. Take $b^\prime$ such that $b=\epsilon^+_Zb^\prime$. Then $\epsilon^+_Zb^\prime a=0$ implies $b^\prime a=0$. Since $b^\prime$ is a morphism in $\cc^\prime$, we obtain $b^\prime=0$ and $b=0$.

(5) We only show the `if' part. Take a minimal right $\cc^\prime$-approximation $Y\stackrel{f}{\to}X$ of $X\in\cc$. Assume that $a\in\cc(Z,Y)$ satisfies $af=0$. Take $a^\prime$ such that $a=\epsilon^+_Za^\prime$. Then $\epsilon^+_Za^\prime f=0$ implies $a^\prime f=0$. Since $f$ is minimal, $a^\prime$ is in $J_{\cc^\prime}=0$. Thus $a=0$.

(6) For any $X\in\cc$, decompose a monic right $\cc^\prime$-approximation $\epsilon^-_X={f_1\choose f_2}:X^-=Y_1\oplus Y_2\to X$ such that $Y_1\in\dd$ and $Y_2$ has no direct summand in $\ind\dd$. Take a monic right $\cc^\prime$-approximation $g\in\cc(Z,Y_2)$. Since $\cc^{\prime\prime}$ is a cosemisimple rejective subcategory of $\cc^\prime$ by (1), $\cc(\ ,Z)\stackrel{\cdot g}{\to}J_{\cc}(\ ,Y_2)$ is an isomorphism on $\cc^\prime$ by \XAEA. Thus ${f_1\choose gf_2}:Y_1\oplus Z\to X$ gives a monic right $\dd$-approximation.\rule{5pt}{10pt}

\vskip.5em{\bf\XBB\ Definition }
Let $\cc$ be an additive category and $\cc^\prime=\cc_m\subseteq\cc_{m-1}\subseteq\cdots\subseteq\cc_0=\cc$ a chain of subcategories. We call it a ({\it half, left, right}) {\it rejective chain} of {\it length} $m$ if $\cc_{n+1}$ is a cosemisimple (left or right, left, right) rejective subcategory of $\cc_n$ for any $n$ ($0\le n<m$) (\S\XAE). A left (resp. right) rejective chain is called {\it total} if each $\cc_{n}$ is a left (resp. right) rejective subcategory of $\cc$.\footnote{Right rejective chains in [I5;2.6] mean our total right rejective chains, and right rejective chains in [I3;2.2] mean our $\Lambda$-total right rejective chains (see \XBCC). We note here that, for the validity of the assertion [I4;2.1], one should assume that the chain there is total (see \XCEA(1)). The main result [I4;2.2] and its proof are valid since the chain constructed there was $\Lambda$-total (see \XBEA(1)).}
Any rejective chain is total left and total right by \XBA(3). By \XBBA\ below, any chain defined above can be refined to a {\it saturated chain} which satisfies $\#(\ind\cc_n\backslash\ind\cc_{n+1})=1$ for any $n$.

Now assume that $\Lambda$ is an artin algebra or order, and $\cc$ is a subcategory of $\dn{M}_\Lambda$. We call $\cc$ {\it maximal} if $\cc=0$ holds (artin algebra case) or $\cc=\lat\Gamma$ holds for some hereditary overring $\Gamma$ of $\Lambda$ (order case). We call the chain above {\it complete} if $\cc^\prime$ is maximal. On the other hand, assume that the chain above is left (resp. right) rejective with a unit $\epsilon_n^+$ (resp. $\epsilon_n^-$) of the natural inclusion $\cc_{n+1}\to\cc_n$. We call the chain {\it $\Lambda$-total} if $\epsilon_{n,X}^+$ (resp. $\epsilon_{n,X}^-$) is epic (resp. monic) in $\dn{M}_\Lambda$ for any $X\in\cc_n$.

Any $\Lambda$-total left (resp. right) rejective chain is total by \XBA(4), and the converse holds if $D\Lambda\in\cc$ (resp. $\Lambda\in\cc$). Although totality depends only on the categorical structure of $\cc$, $\Lambda$-totality depends on the exact structure in $\dn{M}_\Lambda$. If $\Lambda$ is an order in a semisimple algebra, then any left (resp. right) rejective chain is $\Lambda$-total, since a morphism $f\in\hom_\Lambda(X,Y)$ is epic (resp. monic) in $\dn{M}_\Lambda$ if and only if $0\to\endm_\Lambda(Y)\stackrel{f\cdot}{\to}\hom_\Lambda(X,Y)$ (resp. $0\to\endm_\Lambda(X)\stackrel{\cdot f}{\to}\hom_\Lambda(X,Y)$) is exact.

\vskip.5em{\bf\XBBA\ }
Let $\cc^\prime=\cc_m\subseteq\cc_{m-1}\subseteq\cdots\subseteq\cc_0=\cc$ be a chain of subcategories of $\cc$. Take any chain $\cc_{n+1}=\cc_{n,l_n}\subseteq\cdots\subseteq\cc_{n,1}\subseteq\cc_{n,0}=\cc_n$ of subcategories of $\cc_n$ for each $n$ ($0\le n<m$), and consider the refined chain $\cc^\prime=\cc_{m-1,l_{m-1}}\subseteq\cc_{m-1,l_{m-1}-1}\subseteq\cdots\subseteq\cc_{0,1}\subseteq\cc_{0,0}=\cc$. If the original chain is a (half, left, right, total right, total left) rejective chain, then so is the refined chain by \XBA(6).

\vskip.5em{\bf\XBBB\ }The theorem below plays an important role in this paper. We shall give in \XCC\ corresponding results for (not necessarily total) left (resp. right) rejective chains.

\vskip.5em{\bf Theorem }{\it
Let $\Lambda$ be an artin algebra or order, and $\cc$ a finite subcategory of $\dn{M}_\Lambda$. Assume that $\cc$ has a complete $\Lambda$-total left (resp. right) rejective chain of length $m>0$. Then $\resdim{\cc}{\dn{M}_\Lambda}<m(+1)$ (resp. $\resdim{\cc^{op}}{\dn{M}_\Lambda^{op}}<m(+1)$) and $\gl(\mod\cc)\le m(+1)$ hold, where $+1$ are added if $\Lambda$ is an order.}

\vskip.5em{\sc Proof }
We only show the left-hand side assertion.

(i) Let $0\to Z\stackrel{f}{\to}Y\to X$ be an exact sequence with $X\in\dn{M}_\Lambda$, $Y\in\cc_n$ and $f\in J_{\dn{M}_\Lambda}$. We shall show that $Z$ has a right $\cc$-approximation of the form $W\to Z$ with $W\in\cc_{n+1}$.

We only have to show that any morphism $a:W\to Z$ with $W\in\cc_l$ ($l\ge n$) factors through $\cc_{l+1}$. Since $af\in J_{\dn{M}_\Lambda}(W,Y)$, there exists $b$ such that $af=\epsilon^+_{l,W}b$ for the left $\cc_{l+1}$-approximation $\epsilon^+_{l,W}$ of $W\in\cc_l$. Since $\epsilon^+_{l,W}$ is epic in $\dn{M}_\Lambda$, there exists $c$ such that $b=cf$. Thus $a=\epsilon^+_{l,W}c$ holds.

(ii) $\cc_m$ is a rejective subcategory of $\dn{M}_\Lambda$ by \XAFD. By (i), any $X\in\dn{M}_\Lambda$ has a right $\cc$-resolution $0\to Y_m\to Y_{m-1}\to\cdots\to Y_0\to X$ with $Y_n\in\cc_n$ for any $n$. Thus $\resdim{\cc}{\dn{M}_\Lambda}<m(+1)$ holds. For the latter assertion, we can assume either $m>1$ holds or $\Lambda$ is an order. The proof of \XAAA(1) implies $\gl(\mod\cc)\le 2+\sup\{\resdim{\cc}{Z}\}$, where $Z\in\dn{M}_\Lambda$ has an exact sequence $0\to Z\stackrel{a}{\to}Y\to X$ with $X,Y\in\cc$ and $a\in J_{\dn{M}_\Lambda}$. Again by (i), $Z$ has a right $\cc$-resolution $0\to Y_m\to\cdots\to Y_1\to Z$ with $Y_n\in\cc_n$ for any $n$. Thus $\resdim{\cc}{Z}<m-1(+1)$ and $\gl(\mod\cc)\le m(+1)$.\rule{5pt}{10pt}

\vskip.5em{\bf\XBBC\ Example }
(1) Let $\Lambda$ be an artin algebra and $0=I_m\subseteq I_{m-1}\subseteq\cdots\subseteq I_0=\Lambda$ a chain of two-sided ideals of $\Lambda$ such that $I_nJ_\Lambda\subseteq I_{n+1}$ for any $n$. Put $\cc_n:=\add_\Lambda(\bigoplus_{i=0}^{m-n}\Lambda/I_i)$. Then $0=\cc_m\subseteq\cc_{m-1}\subseteq\cdots\subseteq\cc_0$ gives a $\Lambda$-total left rejective chain. Thus $\gl\endm_\Lambda(\bigoplus_{i=0}^{m}\Lambda/I_i)\le m$ holds by \XBBB. The case when $I_n=J_\Lambda^n$ is a classical example of Auslander [A1], and Dlab-Ringel [DR3] proved that $\endm_\Lambda(\bigoplus_{i=0}^{m}\Lambda/J_\Lambda^i)$ is a quasi-hereditary algebra (see \XCEA).

(2) Let $\Lambda$ be an order in a semisimple algebra and $\Lambda=\Lambda_0\subseteq\Lambda_1\subseteq\cdots\subseteq\Lambda_m$ a chain of overorders of $\Lambda$ such that $\Lambda_m$ is hereditary and $J_{\Lambda_n}\in\lat\Lambda_{n+1}$ for any $n$. Put $\cc_n:=\add_\Lambda(\bigoplus_{i=0}^{m-n}\Lambda_i)$. Then $\cc_m\subseteq\cc_{m-1}\subseteq\cdots\subseteq\cc_0$ gives a $\Lambda$-total left rejective chain. Thus $\gl\endm_\Lambda(\bigoplus_{i=0}^{m}\Lambda_i)\le m+1$ holds by \XBBB. Any order has such a chain by putting $\Lambda_{n+1}:=O_l(J_{\Lambda_n})=\{x\in\widetilde{\Lambda}\ |\ xJ_{\Lambda_n}\subseteq J_{\Lambda_n}\}$ [K].

\vskip.5em{\sc Proof }Any $X\in\ind\cc_0\backslash\ind\cc_1$ can be written $X=\Lambda e$ for a primitive idempotent $e$ of $\Lambda$. We shall show that $\cc_1$ is a cosemisimple left rejective subcategory of $\cc_0$ by \XAEA.

(1) The natural surjection $f:\Lambda e\to\Lambda e/I_{m-1}e$ induces a surjection $\cc_0(\Lambda e/I_{m-1}e,\ )\stackrel{f\cdot}{\to}J_{\cc_1}(\Lambda e,\ )$ since any $a\in J_{\cc_0}(\Lambda e,\Lambda/I_i)$ satisfies $a(I_{m-1}e)=0$.

(2) Any $a\in J_{\cc_0}(\Lambda e,\Lambda_i)$ satisfies $a(\Lambda_1 e)\subseteq \Lambda_i$ by $a(\Lambda e)\subseteq J_\Lambda$ ($i=0$) and $\Lambda_1\subseteq\Lambda_i$ ($i>0$). Thus the natural injection $f:\Lambda e\to\Lambda_1 e$ indues a surjection $\cc_0(\Lambda_1 e,\ )\stackrel{f\cdot}{\to}J_{\cc_1}(\Lambda e,\ )$.\rule{5pt}{10pt}

\vskip.5em{\bf\XBC\ Definition }Let $\Lambda$ be an artin algebra or order, and $\cc$ a subcategory of $\dn{M}_\Lambda$. Although in general (left, right) rejective subcategories of $\cc$ cannot be characterized simply as in \XAEB, we will give a simple criterion for a subcategory $\cc^\prime$ of $\cc$ to be cosemisimple right (resp. left) rejective.

(1) Define a functor $\fff_{\cc}:\dn{M}_\Lambda\to\dn{M}_\Lambda$ with a natural transformation $\epsilon^-:\fff_{\cc}\to 1$ by $\fff_{\cc}X:=\sum_{Y\in\cc,\ f\in J_{\dn{M}_\Lambda}(Y,X)}f(Y)=\sum_{Y\in\cc}YJ_{\dn{M}_\Lambda}(Y,X)$ and the natural inclusion $\epsilon^-_X:\fff_{\cc}X\to X$. Moreover, using the duality $D:\dn{M}_\Lambda\leftrightarrow\dn{M}_{\Lambda^{op}}$, we define a functor $\ggg_{\cc}:\dn{M}_\Lambda\to\dn{M}_\Lambda$ with a natural transformation $\epsilon^+:1\to\ggg_{\cc}$ by $\ggg_{\cc}X:=D\fff_{D\cc}(DX)$ and $\epsilon^+_X:=D(\epsilon^-_{DX})$.

(2) For $X\in\ind\cc$, $\epsilon^-_X$ (resp. $\epsilon^+_X$) is an isomorphism if and only if any exact sequence $0\to Z\to Y\to X\to0$ (resp. $0\to X\to Y\to Z\to0$) in $\mod\Lambda$ with $Y\in\cc$ and $Z\in\dn{M}_\Lambda$ splits. We call such $X$ {\it splitting projective} (resp. {\it injective}) in $\cc$ [AS2]. Moreover, if $X$ is not isomorphic to $\fff_{\cc}X$ (resp. $\ggg_{\cc}X$), then we define a subcategory $\cc^\prime$ of $\cc$ by $\ind\cc^\prime=\ind\cc\backslash\{ X\}$. We call this process to get $\cc^\prime$ from $\cc$ a {\it right} (resp. {\it left}) {\it cancellation}, and we call a left or right cancellation a {\it half} cancellation. Our construction generalize preprojective partition [AS2,3] of Auslander-Smalo to arbitrary $\cc$.

\vskip.5em{\bf\XBCA\ Proposition }{\it
Let $\Lambda$ be an artin algebra or order in a semisimple algebra, and $\cc$ a subcategory of $\dn{M}_\Lambda$.

(1) Let $\cc^\prime$ be a subcategory of $\cc$. If $\fff_{\cc}X\in\cc^\prime$ (resp. $\ggg_{\cc}X\in\cc^\prime$) holds for any $X\in\ind\cc\backslash\ind\cc^\prime$, then $\cc^\prime$ is a cosemisimple right (resp. left) rejective subcategory of $\cc$.

(2) If $\cc$ is a non-maximal (\XBB) covariantly (resp. contravariantly) finite subcategory of $\dn{M}_\Lambda$, then $\cc$ has a splitting projective (resp. injective) object $X\in\ind\cc$ which is not isomorphic to $\fff_{\cc}X$ (resp. $\ggg_{\cc}X$). Thus a right (resp. left) cancellation is applicable to $\cc$.

(3) Assume that $\fff_{\cc}$ (resp. $\ggg_{\cc}$) gives an endofunctor of $\cc$. Then any subcategory of $\cc$ obtained by a right (resp. left) cancellation is cosemisimple right (resp. left) rejective.}

\vskip.5em{\sc Proof }
(1) Fix $X\in\ind\cc\backslash\ind\cc^\prime$. Then $\fff_{\cc}X\in\cc^\prime$ implies $\cc(\ ,\fff_{\cc}X)\epsilon^-_X\subseteq J_{\cc}(\ ,X)$. By the construction of $\fff_{\cc}$, $\cc(\ ,\fff_{\cc}X)\stackrel{\cdot\epsilon^-_X}{\longrightarrow}J_{\cc}(\ ,X)$ is an isomorphism on $\cc$. Thus $\cc^\prime$ is a cosemisimple right rejective subcategory of $\cc$ by \XAEA.

(2)(i) We will show the existence of splitting projective object [AS2].

Let $f:\Lambda\to X$ be a minimal left $\cc$-approximation of $\Lambda$. Then $X\neq0$ holds by $\cc\neq0$. We will show that $\fff_{\cc}X\neq X$ holds. Otherwise, there exists a surjection $a\in J_{\cc}(Y,X)$. We can take $b$ such that $f=ba$. Since $f$ is a left $\cc$-approximation, there exists $c$ such that $b=fc$. then $f(1-ca)=0$ implies that $ca$ is an automorphism, a contradiction.

(ii) Assume that $X\in\ind\cc$ is splitting projective and isomorphic to $\fff_{\cc}X$. This is impossible for the artin algebra case, so we consider the order case. Put $X_0:=X$ and $X_{n+1}:=\fff_{\cc}X_n$ for $n\ge0$. Let $f\in J_{\endm_\Lambda(X)}$ be the composition of the isomorphism $X\to\fff_{\cc}X$ and $\epsilon^-_X$. Since $X_n=f^n(X)$ holds, we obtain $\bigcap_{n\ge0}X_n=0$. Now let $\Delta$ be a maximal overorder of $\endm_\Lambda(X)$ and assume $Y:=X\Delta\neq X$. Since there exists a surjection $X^l\to Y$ for some $l>0$, any morphism $Y\to X_n$ factors through $X_{n+1}$. This implies $\hom_\Lambda(Y,X)=0$, a contradiction. Thus $\endm_\Lambda(X)$ is a maximal order, and $\add X=\lat\Delta$ holds for some maximal overring $\Delta$ of $\Lambda$. Since any morphism $Z\to X_n$ with $Z\in\ind\cc\backslash\{X\}$ factors through $X_{n+1}$, we obtain $\hom_\Lambda(Z,X)=0$. This implies $\cc=\cc^\prime\times\lat\Delta$ for $\cc^\prime:=\add(\ind\cc\backslash\{X\})$. Applying the same argument to $\cc^\prime$, we shall obtain a desired splitting projective object by the non-maximality of $\cc$.

(3) Since $\fff_{\cc}X\in\cc^\prime$ holds by the construction, the assertion follows from (1).\rule{5pt}{10pt}

\vskip.5em{\bf\XBCB\ }
Let $\Lambda$ be an artin algebra or order, and $\cc$ a subcategory of $\dn{M}_\Lambda$.

(1) Assume that $\cc$ is closed under factor modules (resp. submodules). Then $\fff_{\cc}$ (resp. $\ggg_{\cc}$) gives an endofunctor of $\cc$. Any subcategory of $\cc$ obtained by a right (resp. left) cancellation is closed under factor modules (resp. submodules).

(2) Assume that $\Lambda$ is an artin algebra and $\cc$ is closed under images. Then $\fff_{\cc}$ (resp. $\ggg_{\cc}$) gives an endofunctor of $\cc$. Any subcategory of $\cc$ obtained by a half cancellation is closed under images.

\vskip.5em{\sc Proof }
(1) Since there exists a surjection $Y\to \fff_{\cc}X$ with $Y\in\cc$, the former assertion follows. Let $P\in\ind\cc$ be splitting projective and $\cc^\prime$ a subcategory of $\cc$ such that $\ind\cc^\prime=\ind\cc\backslash\{ P\}$. Take any surjection $Y\to X$ with $Y\in\cc^\prime$ and $X\in\dn{M}_\Lambda$. Then $X\in\cc$ holds, and $P$ is not a direct summand of $X$ by $\fff_{\cc}P\neq P$. Thus $X\in\cc^\prime$.

(2) Similar argument as in the proof of (1) works.\rule{5pt}{10pt}

\vskip.5em{\bf\XBCC\ }By (1) below, the definition of right rejective chains given in [I3;2.2] is equivalent to our $\Lambda$-total right rejective chain in \XBB.

\vskip.5em{\bf Theorem }{\it
Let $\Lambda$ be an artin algebra or order, and $\cc$ a finite subcategory of $\dn{M}_\Lambda$. 

(1) Assume that $\cc$ is closed under factor modules (resp. submodules). Then any successive right (resp. left) cancellation gives a complete $\Lambda$-total right (resp. left) rejective chain consisting of subcategories which are closed under factor modules (resp. submodules). Conversely, any saturated complete $\Lambda$-total right (resp. left) rejective chain of $\cc$ is obtained by a successive right (resp. left) cancellation.

(2) If $\Lambda$ is an artin algebra and $\cc$ is closed under images, then any successive half (resp. left, right) cancellation gives a complete half (resp. $\Lambda$-total right, $\Lambda$-total left) rejective chain consisting of subcategories which are closed under images.}

\vskip.5em{\sc Proof }
We only have to show the latter part of (1) since other assertions are immediate from \XBCA(2)(3) and \XBCB. Assume that $\cc^\prime\supset\cc^{\prime\prime}$ is a consecutive two terms in a saturated $\Lambda$-total right rejective chain of $\cc$ and $\{X\}:=\ind\cc^\prime\backslash\ind\cc^{\prime\prime}$. Let $f:X^-\to X$ be a right $\cc^{\prime\prime}$-approximation of $X$ which is monic in $\dn{M}_\Lambda$. Since $X^-\subsetneq X$ holds, we obtain $\fff_{\cc^\prime}X=X^-\neq X$ and $X$ is splitting projective.\rule{5pt}{10pt}

\vskip.5em{\bf\XBCD\ Example }
In this subsection, we shall see that a few well-known examples of (left, right) rejective chains are given by \XBCC. We notice that Rejection Lemma of Drozd-Kirichenco \XAFE\ is a left and right cancellation simultaneously.

(1) Let $\Lambda$ be an artin algebra of finite representation type. In [DR4], Dlab-Ringel constructed chains of subcategories of $\mod\Lambda$ by using (i) Roiter measure, (ii) dual Roiter measure [Roi][G][R3], (iii) preinjective partition, or (iv) preprojective partition [AS2,3]. In the chain constructed by (i) or (iii) (resp. (ii) or (iv)), each subcategories were closed under submodules (resp. factor modules). Thus the chain is a $\Lambda$-total left (resp. right) rejective by \XBCC. On the other hand, Dlab-Ringel proved that the chain gives a heredity chain of the Auslander algebra of $\Lambda$. This follows from their more general theorem on splitting filtrations [DR4], which will be explained in \XCG\ from our categorical viewpoint.

(2) Let $\Lambda_0$ be a cyclic Nakayama artin algebra. Then any indecomposable $\Lambda_0$-module is local-colocal. Take $X\in\ind(\mod\Lambda_0)$ with a maximal length. Then one can easily show $X\in\pr\Lambda_0\cap\rin\Lambda_0$. Thus there exists a factor algebra $\Lambda_1$ of $\Lambda_0$ such that $\{ X\}=\ind(\mod\Lambda_0)-\ind(\mod\Lambda_1)$ by \XAFE, and $\Lambda_1$ is a cyclic Nakayama again. Repeating this process, we obtain a complete rejective chain $0=\mod\Lambda_m\subset\cdots\subset\mod\Lambda_1\subset\mod\Lambda_0$.

(3) We call an order $\Lambda$ {\it Gorenstein} if $\Lambda\in\rin\Lambda$, and {\it Bass} if any overorder of $\Lambda$ is Gorenstein [DKR]. By definition, any overorder of a Bass order is Bass.

Let $\Lambda_0$ be a Bass order in a semisimple algebra. Then  $\pr\Lambda_0=\rin\Lambda_0$ holds since $\Lambda_0$ is Gorenstein. Thus, if $\Lambda_0$ is not maximal, then there exist $X\in\ind(\pr\Lambda_0)$ and an overorder $\Lambda_1$ of $\Lambda_0$ such that $\{ X\}=\ind(\lat\Lambda_0)-\ind(\lat\Lambda_1)$ by \XAFE. Repeating this process, we obtain a complete rejective chain $\lat\Lambda_m\subset\cdots\subset\lat\Lambda_1\subset\lat\Lambda_0$ since any increasing chain of overorders of $\Lambda_0$ stops. Such a chain was called a {\it Bass chain}, and plays a crucial role in the theory of Bass orders [DK][HN][Ro].

\vskip.5em{\bf\XBD\ }In this section, we shall apply \XBCA\ to construct rejective chains for more general classes of subcategories of $\dn{M}_\Lambda$, which are not necessarily closed under factor modules (resp. submodules). In fact, the class of subcategories satisfying the condition in \XBDA\ below is much larger than that in \XBCC. Our results \XBDA\ and \XBDB\ are immediate coonsequence of the following kew theorem.

\vskip.5em{\bf Theorem }{\it
Let $\Lambda$ be an artin algebra or order, and $\cc$ a finite subcategory of $\dn{M}_\Lambda$. Assume that $\fff_{\cc^\prime}(\cc^\prime)\subseteq\cc$ (resp. $\ggg_{\cc^\prime}(\cc^\prime)\subseteq\cc$) holds for any subcategory $\cc^\prime$ of $\cc$. Then any successive right (resp. left) cancellation gives a complete $\Lambda$-total right (resp. left) rejective chain of $\cc$.}

\vskip.5em{\sc Proof }
Let $P$ be splitting projective in $\cc$ and $\cc^\prime$ a subcategory of $\cc$ such that $\ind\cc^\prime=\ind\cc\backslash\{P\}$. By \XBCA(2)(3), we only have to show that $\fff_{\cc^{\prime\prime}}(\cc^{\prime\prime})\subseteq\cc^\prime$ holds for any subcategory $\cc^{\prime\prime}$ of $\cc^\prime$. For any $X\in\cc^{\prime\prime}$, $P$ is not a direct summand of $\fff_{\cc^{\prime\prime}}X$ by $\fff_{\cc}P\neq P$. Since $\fff_{\cc^{\prime\prime}}X\in\cc$ holds by our assumption, we obtain $\fff_{\cc^{\prime\prime}}X\in\cc^\prime$.\rule{5pt}{10pt}

\vskip.5em{\bf\XBDA\ Corollary }{\it
(1) Let $\Lambda$ be an artin algebra and $\cc$ a finite subcategory of $\mod\Lambda$. Assume that any submodule (resp. factor module) of any $X\in\ind\cc$ is contained in $\cc$. Then $\cc$ has a complete $\Lambda$-total right (resp. left) rejective chain.

(2) Let $\Lambda$ be an order and $\cc$ a finite subcategory of $\lat\Lambda$. Assume that $Y\in\cc$ holds for any $X\in\ind\cc$ and $Y\in\lat\Lambda$ such that $\widetilde{Y}$ is a submodule (resp. factor module) of $\widetilde{X}$. Then $\cc$ has a complete $\Lambda$-total right (resp. left) rejective chain.}

\vskip.5em{\bf\XBDB\ }The result (1) below played a crucial role in the proof of Solomon's conjecture on zeta functions of orders [I3,5].

\vskip.5em{\bf Corollary }{\it
Let $\Lambda$ be an artin algebra (resp. order).

(1) For $n\in\nnn$, put $\cc^{(n)}:=\add\{ X\in\ind\dn{M}_\Lambda\ |\ \length_\Lambda X<n$ (resp. $\length_{\widetilde{\Lambda}}\widetilde{X}<n$)$\}$. If $\#\ind\cc^{(n)}<\infty$, then $\cc^{(n)}$ has a complete $\Lambda$-total left rejective chain and a complete $\Lambda$-total right rejective chain.

(2) For $M\in\dn{M}_\Lambda$, put $\cc_M:=\add\{ X\in\ind\dn{M}_\Lambda\ |\ X$ (resp. $\widetilde{X}$) is a submodule of $M$ (resp. $\widetilde{M}$)$\}$ and $\cc^M:=\add\{ X\in\ind\dn{M}_\Lambda\ |\ X$ (resp. $\widetilde{X}$) is a factor module of $M$ (resp. $\widetilde{M}$)$\}$. If $\#\ind\cc_M<\infty$, then $\cc_M$ has a complete $\Lambda$-total right rejective chain. Dually, if $\#\ind\cc^M<\infty$, then $\cc^M$ has a complete $\Lambda$-total left rejective chain.}

\vskip.5em{\bf\XBE\ }
Let $\Lambda$ be an artin algebra or order. We shall construct a rejective chain for an arbitrary given $M_0:=M\in\dn{M}_\Lambda$. Define $M_n$ inductively by (i) (resp. (ii)) below. Fix $m>0$ and put $\cc_n:=\add_\Lambda(\bigoplus_{i=n}^{m}M_i)$. Then $\cc_m\subseteq\cc_{m-1}\subseteq\cdots\subseteq\cc_0=\cc$ forms a $\Lambda$-total right (resp. left) rejective chain.

(i) Decompose $M_n=X_n\oplus Y_n$ arbitrary, and take arbitrary factor module $Z_n\in\dn{M}_\Lambda$ of a direct sum of copies of $M_n$. Put $M_{n+1}:=M_nJ_{\dn{M}_\Lambda}(M_n,X_n)\oplus Y_n\oplus Z_n$.

(ii) Decompose $M_n=X_n\oplus Y_n$ arbitrary, and take arbitrary submodule $Z_n\in\dn{M}_\Lambda$ of a direct sum of copies of $M_n$. Put $M_{n+1}:=D((DM_n)J_{\dn{M}_\Lambda}(DM_n,DX_n))\oplus Y_n\oplus Z_n$.

\vskip.5em{\sc Proof }
$M_{n+1}$ is a factor module of a direct sum of copies of $M_n$. Thus so is $M_i$ for any $i$ ($i>n$). In particular, $M_iJ_{\dn{M}_\Lambda}(M_i,X)\subseteq M_nJ_{\dn{M}_\Lambda}(M_n,X)$ holds for any $X\in\dn{M}_\Lambda$ and $i$ ($i>n$), and we obtain $\fff_{\cc_n}(X)=M_nJ_{\dn{M}_\Lambda}(M_n,X)$. Fix $X\in\ind\cc_n\backslash\ind\cc_{n+1}$. Then $X\in\add X_n$ holds. Thus $\fff_{\cc_n}(X)=M_nJ_{\dn{M}_\Lambda}(M_n,X)\in\add M_{n+1}\subseteq\cc_{n+1}$ implies that $\cc_{n+1}$ is a cosemisimple right rejective subcategory of $\cc_n$ by \XBCA(1).\rule{5pt}{10pt}

\vskip.5em{\bf\XBEA\ Theorem }(cf. [I4;2.2]){\it 
Let $\Lambda$ be an artin algebra or order in a semisimple algebra, and $M_0:=M\in\dn{M}_\Lambda$. Put $M_{n+1}:=M_nJ_{\endm_\Lambda(M_n)}$ (resp. $M_{n+1}:=D((DM_n)J_{\endm_\Lambda(DM_n)})$) inductively. Then there exists $m>0$ such that $\cc_n:=\add_\Lambda(\bigoplus_{i=n}^{m}M_i)$ gives a complete $\Lambda$-total right (resp. left) rejective chain $\cc_m\subseteq\cc_{m-1}\subseteq\cdots\subseteq\cc_0$.}

\vskip.5em{\sc Proof }
The chain is $\Lambda$-total right (resp. left) rejective since we can obtain it by putting $X_n=M_n$, $Y_n=Z_n=0$ in each step in \XBE. We only have to show the completeness. For the artin algebra case, $M_m=0$ holds for sufficiently large $m$ and such $m$ gives a complete chain. Let us consider the order case. Put $\Delta_n:=\endm_\Lambda(M_n)$. By definition, $\Delta_{n+1}=\endm_\Lambda(M_nJ_{\Delta_n})\supseteq\endm_{\Delta_n}(J_{\Delta_n})\supseteq\Delta_n$ holds. Since any increasing chain of orders in a semisimple algebra stops, there exists $m$ such that $\Delta_m=\Delta_n$ holds for any $n>m$. In this case, $\endm_{\Delta_m}(J_{\Delta_m})=\Delta_m$ holds, and this implies that $\Delta_m$ is a hereditary order [CR]. Thus $\Gamma:=\endm_{\Delta_m^{op}}(M_m)$ is a hereditary overring of $\Lambda$ such that $\lat\Gamma=\add M_m$, and we have shown the assertion.\rule{5pt}{10pt}

\vskip.5em{\bf\XBEB\ Corollary }{\it
Let $\Lambda$ be an artin algebra or order in a semisimple algebra. Then any $M\in\dn{M}_\Lambda$ is contained in a finite subcategories $\cc$ of $\dn{M}_\Lambda$ such that $\cc$ has a complete $\Lambda$-total right (resp. left) rejective chain. Moreover, $\gl(\mod\cc)<\infty$ holds, and $\gl(\mod\cc)\le \length_\Lambda M_{\endm_\Lambda(M)}$ holds if $\Lambda$ is an artin algebra.}

\vskip.5em{\sc Proof }The inequality follows from \XBBB\ and the construction in \XBEA\ since $M_m=0$ holds for $m:=\length_\Lambda M_{\endm_\Lambda(M)}$.\rule{5pt}{10pt}

\vskip.5em{\bf\XC\ Quasi-hereditary algebras and Global dimension }

\vskip.5em{\bf\XCA\ }
There exists a bijection between equivalence classes of Krull-Schmidt categories $\cc$ with additive generators $M$ and Morita-equivalence classes of semiperfect rings $\Gamma$, which is given by $\cc\mapsto\cc(M,M)$ and the converse is given by $\Gamma\mapsto\pr\Gamma$.

Assume that $\cc$ corresponds to $\Gamma$. For idempotents $e$ and $f$ of $\Gamma$, we say that $e$ is {\it equivalent} to $f$ if $\add_\Gamma(\Gamma e)=\add_\Gamma(\Gamma f)$ holds. Then the set of subcategories $\cc^\prime$ of $\cc$, the set of idempotent ideals $I$ of $\Gamma$, and the set of equivalence classes of idempotents $e$ correspond bijectively. This is given by $\cc^\prime\mapsto I:=[\cc^\prime](M,M)$, and $I\mapsto e$ such that $I=\Gamma e\Gamma$. In this case, $\cc^\prime$ corresponds to the ring $e\Gamma e$.

\vskip.5em{\bf\XCB\ Theorem }{\it
Let $\Gamma$ be a basic semiperfect ring, $e$ an idempotent $e$ of $\Gamma$, $\Gamma^\prime:=e\Gamma e$ and $\overline{\Gamma}:=\Gamma/\Gamma e\Gamma$. Assume that $\cc$ and $\cc^\prime$ correspond to $\Gamma$ and $\Gamma^\prime$ respectively by \XCA.

(1) $\cc^\prime$ is a coreflective (resp. reflective) subcategory of $\cc$ if and only if $e\Gamma\in\pr\Gamma^\prime$ (resp. $\Gamma e\in\pr\Gamma^{\prime}{}^{op}$).

(2) $\cc^\prime$ is a right (resp. left) rejective subcategory of $\cc$ if and only if $\Gamma e\Gamma\in\pr\Gamma$ (resp. $\Gamma e\Gamma\in\pr\Gamma^{op}$). Then $0\le \gl\Gamma-\gl\Gamma^\prime\le\gl\overline{\Gamma}+2$.

(3) $\cc^\prime$ is a cosemisimple right (resp. left) rejective subcategory of $\cc$ if and only if $J_\Gamma(1-e)\in\add_\Gamma(\Gamma e)$ (resp. $(1-e)J_\Gamma\in\add(e\Gamma)_\Gamma$). Then $0\le \gl\Gamma-\gl\Gamma^\prime\le2$ holds, where the right inequality is strict if $\hom_\Gamma(\overline{\Gamma}(1-e),\Gamma e)=0$.

(4) If $\cc^\prime$ is a rejective subcategory of $\cc$, then $\gl\Gamma\le\max\{\gl\overline{\Gamma}+2,\gl\Gamma^\prime\}$.}

\vskip.5em{\bf\XCBA\ }
Let ${\bf X}:0\to X_n\to X_{n-1}\to\cdots\to X_0\stackrel{f}{\to}X_{-1}\to 0$ be a complex of $\Gamma$-modules.

(1) If ${\bf X}$ is exact, then $\pd_\Gamma X_{-1}\le\max\{i+\pd_\Gamma X_i\ |\ 0\le i\le n\}$.

(2) If $X_i\in\pr\Gamma$ for any $i\ge0$, then $\pd_\Gamma X_{-1}\le\max\{i+1+\pd_\Gamma H^i({\bf X})\ |\ -1\le i\le n\}$.

\vskip.5em{\sc Proof }Since both are easily shown inductively, we only show (2). Considering the complex $0\to X_n\to X_{n-1}\to\cdots\to X_1\to\Ker f\to 0$, we obtain $\pd_\Gamma\Ker f\le\max\{i+\pd_\Gamma H^i({\bf X})\ |\ 0\le i\le n\}$. Then the exact sequence $0\to\Ker f\to X_0\stackrel{f}{\to}X_{-1}\to H^{-1}({\bf X})\to0$ shows the assertion.\rule{5pt}{10pt}

\vskip.5em{\bf\XCBB\ Lemma }{\it Let $\Gamma$ be a ring, $e$ an idempotent of $\Gamma$, $\Gamma^\prime:=e\Gamma e$ and $\overline{\Gamma}:=\Gamma/\Gamma e\Gamma$. Then $-\pd_{\Gamma^\prime}(e\Gamma)\le\gl\Gamma-\gl\Gamma^\prime\le\gl\overline{\Gamma}+\pd_\Gamma\overline{\Gamma}+1$ holds, where the right inequality is strict if $\hom_\Gamma(\overline{\Gamma}(1-e),\Gamma e)=0$ (cf. [APT]).}

\vskip.5em{\sc Proof }
(i) We will show the left inequality. We can assume $n:=\gl\Gamma<\infty$. Take a projective resolution $0\to P_n\to P_{n-1}\to\cdots \to P_0\to X\to0$ of $X\in\mod\Gamma$. Taking $e\Gamma\otimes_\Gamma\ $, we obtain an exact sequence $0\to eP_n\to eP_{n-1}\to\cdots \to eP_0\to eX\to0$ of $\Gamma^\prime$-modules. Thus $\pd{}_{\Gamma^\prime}(eX)\le n+\pd_{\Gamma^\prime}(e\Gamma)$ holds by \XCBA(1). This implies $\gl\Gamma^\prime\le\gl\Gamma+\pd_{\Gamma^\prime}(e\Gamma)$ since the functor $e\Gamma\otimes_\Gamma\ :\mod\Gamma\to\mod\Gamma^\prime$ is dense.

(ii) We will show the right inequality. We can assume $n:=\gl\Gamma^\prime<\infty$. For $X\in\mod\Gamma$, take a projective resolution ${\bf A}:0\to P_n\to P_{n-1}\to\cdots \to P_0\to eX\to0$ of $\Gamma^\prime$-modules. From ${\bf A}$, we obtain a complex ${\bf B}:0\to\Gamma e\otimes_{\Gamma^\prime} P_n\to\Gamma e\otimes_{\Gamma^\prime} P_{n-1}\to\cdots \to\Gamma e\otimes_{\Gamma^\prime} P_0\to X\to0$. Since $e\Gamma\otimes_\Gamma{\bf B}={\bf A}$ holds, any homology of ${\bf B}$ are $\overline{\Gamma}$-modules. Thus we obtain $\pd_\Gamma X\le\sup\{\pd_\Gamma Y\ |\ Y\in\mod\overline{\Gamma}\}+n+1$ by \XCBA(2). Taking a projective resolution of $Y\in\mod\overline{\Gamma}$ as a $\overline{\Gamma}$-module, we obtain $\pd_\Gamma Y\le\gl\overline{\Gamma}+\pd_\Gamma\overline{\Gamma}$ by \XCBA(1). Thus we obtain $\gl\Gamma\le\gl\Gamma^\prime+\gl\overline{\Gamma}+\pd_\Gamma\overline{\Gamma}+1$.\rule{5pt}{10pt}

\vskip.5em{\bf\XCBC\ Proof of \XCB\ }
(1) $\cc^\prime$ is a coreflective subcategory of $\cc$ if and only if there exists $a\in\cc(P,\Gamma)$ with $P\in\cc^\prime$ such that $\cc(\Gamma e,P)\stackrel{\cdot a}{\to}\cc(\Gamma e,\Gamma)$ is an isomorphism if and only if there exists an isomorphism $eP\to e\Gamma$ of $\Gamma^\prime$-modules with $P\in\cc^\prime$ if and only if $e\Gamma\in\pr\Gamma^\prime$.

(2) $\cc^\prime$ is a right rejective subcategory of $\cc$ if and only if there exists $a\in\cc(P,\Gamma)$ with $P\in\cc^\prime$ such that $P=\cc(\Gamma,P)\stackrel{\cdot a}{\to}[\cc^\prime](\Gamma,\Gamma)=\Gamma e\Gamma$ is an isomorphism if and only if $\Gamma e\Gamma\in\add_\Gamma(\Gamma e)$ if and only if $\Gamma e\Gamma\in\pr\Gamma$. The latter assertion follows from \XCBB\ since $\Gamma e\Gamma\in\pr\Gamma$ implies $\pd_\Gamma\overline{\Gamma}\le 1$ and $\pd_{\Gamma^\prime}(e\Gamma)=0$ by $e\Gamma e\Gamma=e\Gamma$.

(3) Put $f:=1-e$. By \XAEA, $\cc^\prime$ is a cosemisimple right rejective subcategory of $\cc$ if and only if there exists $a\in\cc(P,\Gamma f)$ with $P\in\cc^\prime$ such that $P=\cc(\Gamma,P)\stackrel{\cdot a}{\to}J_{\cc}(\Gamma,\Gamma f)=J_\Gamma f$ is an isomorphism if and only if $J_\Gamma f\in\add_\Gamma(\Gamma e)$.


(4) We can assume $n:=\gl\Gamma^\prime<\infty$. For any $X\in\mod\Gamma$, take a projective resolution ${\bf A}:0\to P_n\to P_{n-1}\to\cdots \to P_0\to eX\to0$ of a $\Gamma^\prime$-module $eX$. From ${\bf A}$, we obtain a complex ${\bf B}:0\to\Gamma e\otimes_{\Gamma^\prime}P_n\to\Gamma e\otimes_{\Gamma^\prime}P_{n-1}\to\cdots\to\Gamma e\otimes_{\Gamma^\prime}P_0\to X\to0$. Take an epic left $\cc^\prime$-approximation $f:\Gamma\to P$ with $P\in\cc^\prime=\add_\Gamma(\Gamma e)$. Then we obtain the following commutative diagram, whose vertical maps are isomorphisms.
\[\begin{diag}
{\bf B}:&0&\to&\Gamma e\otimes_{\Gamma^\prime}P_n&\to&\cdots&\stackrel{g}{\to}&\Gamma e\otimes_{\Gamma^\prime}P_0&\to&X&\to&0\\
&&&\uparrow^{f\cdot}&&&&\uparrow^{f\cdot}\\
&0&\to&\hom_\Gamma(P,\Gamma e\otimes_{\Gamma^\prime}P_n)&\to&\cdots&\to&\hom_\Gamma(P,\Gamma e\otimes_{\Gamma^\prime}P_0)
\end{diag}\]

Since $\hom_\Gamma(\Gamma e,{\bf B})={\bf A}$ holds, the lower sequence is exact by $P\in\add_\Gamma(\Gamma e)$. Thus we obtain $H^i({\bf B})=0$ except $i=0,1$. We have an exact sequence $0\to H^0({\bf B})\to\Cok g\to X\to H^{-1}({\bf B})\to0$ with $\pd_\Gamma\Cok g\le n$. Since $H^i({\bf B})$ is a $\overline{\Gamma}$-module, $\pd_\Gamma H^i({\bf B})\le\gl\overline{\Gamma}+1$ holds. Thus we have $\pd_\Gamma X\le\max\{\gl\overline{\Gamma}+2,n\}$.\rule{5pt}{10pt}

\vskip.5em{\bf\XCC\ Theorem }{\it
Assume that $\cc$ corresponds to $\Gamma$ by \XCA.

(1) Assume that $\Gamma$ is an artin algebra. If $\cc$ has a complete total right (resp. left) rejective chain of length $m$, then $\gl\Gamma\le m$. If $\cc$ has a complete half rejective chain of length $m$, then $\gl\Gamma\le 2m-2$.

(2) Assume that $\Gamma$ is an order in a semisimple algebra. If $\cc$ has a complete half rejective chain of length $m$, then $\gl\Gamma\le m+1$.}

\vskip.5em{\sc Proof }All assertions follow from \XCB(3) except the former assertion of (1). It follows from \XBBB\ since the equivalence $\cc=\pr\Gamma$ (resp. $\cc=\rin\Gamma$) makes the chain $\Gamma$-total.\rule{5pt}{10pt}

\vskip.5em{\bf\XCCA\ Example }
Let $\left(\def\arraystretch{.5}\begin{array}{cccccc}{\scriptstyle R}&{\scriptstyle R}&{\scriptstyle R}&{\scriptstyle\cdots}&{\scriptstyle R}&{\scriptstyle R}\\ {\scriptstyle  J}&{\scriptstyle R}&{\scriptstyle R}&{\scriptstyle\cdots}&{\scriptstyle R}&{\scriptstyle R}\\ {\scriptstyle  J^2}&{\scriptstyle  J}&{\scriptstyle R}&{\scriptstyle\cdots}&{\scriptstyle R}&{\scriptstyle R}\\ {\scriptstyle\cdots}&{\scriptstyle\cdots}&{\scriptstyle\cdots}&{\scriptstyle\cdots}&{\scriptstyle\cdots}\\ {\scriptstyle  J^{n-2}}&{\scriptstyle  J^{n-3}}&{\scriptstyle  J^{n-4}}&{\scriptstyle\cdots}&{\scriptstyle R}&{\scriptstyle R}\\ {\scriptstyle  J^{n-1}}&{\scriptstyle  J^{n-2}}&{\scriptstyle  J^{n-3}}&{\scriptstyle\cdots}&{\scriptstyle J}&{\scriptstyle R}\end{array}\right)\subseteq\Gamma\subseteq\ma_n(R)$ be $R$-orders with $J=J_R$. Then $\Gamma$ has a complete total right (resp. left) rejective chain of length $n-2$. Thus $\gl\Lambda\le n-1$ holds, a classical result of Jategaonkar [J].

\vskip.5em{\sc Proof }
Without loss of generality, we can assume that $\Gamma$ is basic. Let $e_i$ be an idempotent of $\Gamma$ with $(i,i)$-th entry $1$ and other entries $0$. Then $J_\Gamma e_n=\Gamma e_{n-1}$ holds. Thus $\add_\Gamma(\Gamma(1-e_n))$ is a cosemisimple right rejective subcategory of $\pr\Gamma$ by \XCB(3). Since $(e_1+e_2)\Gamma(e_1+e_2)$ is hereditary, the assertion follows inductively.\rule{5pt}{10pt}

\vskip.5em{\bf\XCD\ }We recall two classes of algebras closely related to rejective chains. One is quasi-hereditary algebras of Cline-Parshall-Scott [CPS1,2], and another is neat algebras of Agoston-Dlab-Wakamatsu [ADW].

(1) A two-sided ideal $I$ of an artin algebra $\Gamma$ is called {\it heredity} if $I^2=I$, $IJ_\Gamma I=0$ and $I\in\pr\Gamma$. This condition is left-right symmetric since the last condition is equivalent to $I\in\pr\Gamma^{op}$ by \XBA(5). An artin algebra $\Gamma$ is called {\it quasi-hereditary} if it has a {\it heredity chain}, which is a chain $0=I_m\subseteq I_{m-1}\subseteq\cdots\subseteq I_0=\Gamma$ of ideals of $\Gamma$ such that $I_{n-1}/I_n$ is a heredity ideal of $\Gamma/I_n$ for any $n$ ($0<n\le m$).

An order $\Gamma$ is called {\it quasi-hereditary} if there exists an idempotent $e$ of $\Gamma$ such that $e\Gamma e$ is a maximal order and $\Gamma/\Gamma e\Gamma$ is a quasi-hereditary artin algebra [KW].

(2) An idempotent $f$ of an artin algebra $\Gamma$ is called {\it neat} if $\ext^i_\Gamma((\Gamma/J_\Gamma)f,(\Gamma/J_\Gamma)f)=0$ holds for any $i>0$. This is equivalent to that $\ext^i_{\Gamma^{op}}(f(\Gamma/J_\Gamma),f(\Gamma/J_\Gamma))=0$ holds for any $i>0$ since $D((\Gamma/J_\Gamma)f)=f(\Gamma/J_\Gamma)$. An artin algebra $\Gamma$ is called {\it neat} if it has a {\it neat sequence} $(f_1,f_2,\cdots,f_m)$, which is a complete set of orthogonal idempotents of $\Gamma$ such that $f_t$ is a neat idempotent of $\epsilon_t\Gamma\epsilon_t$ for $\epsilon_t:=f_t+f_{t+1}+\cdots+f_m$ for any $t$ ($1\le t\le m$). It is known that any quasi-hereditary algebra is neat [ADW].

\vskip.5em{\bf\XCDA\ }Let us recall theorems in [CPS1,2] and [ADW;prop.2] concerning on the global dimension. For completeness, we shall give a quick proof of them.

\vskip.5em{\bf Theorem }{\it Let $\Gamma$ be a basic artin algebra.

(1) If $I$ is a heredity ideal of $\Gamma$, then $0\le \gl\Gamma-\gl\Gamma/I\le 2$ holds. Consequently, if $\Gamma$ is a quasi-hereditary algebra with a chain $0=I_m\subseteq I_{m-1}\subseteq\cdots\subseteq I_0=\Gamma$, then $\gl\Gamma\le 2m-2$ holds.

(2) If $f$ be a neat idempotent of $\Gamma$, then $\frac{1}{2}(\gl\Gamma^\prime+1)\le\gl\Gamma\le2(\gl\Gamma^\prime+1)$ holds for $e:=1-f$ and $\Gamma^\prime:=e\Gamma e$. Consequently, if $\Gamma$ is a neat algebra with a neat sequence $(f_1,f_2,\cdots,f_m)$, then $\gl\Gamma\le 2^m-2$ holds.}

\vskip.5em{\sc Proof }
(1) Put $I=\Gamma e\Gamma$. Then the assertion follows from \XCBB\ since $\Gamma^\prime$ is semisimple.

(2) Let ${\bf P}$ be a minimal projective resolution of a $\Gamma$-module $J_\Gamma f$. Since $f$ is neat, all terms in ${\bf P}$ are in $\add_\Gamma(\Gamma e)$. Since the functor $e\Gamma\otimes_\Gamma\ :\add_\Gamma(\Gamma e)\to\pr\Gamma^\prime$ is an equivalence, $e\Gamma\otimes_\Gamma{\bf P}$ gives a minimal projective resolution of a $\Gamma^\prime$-module $eJ_\Gamma f$. Thus we obtain $\pd_\Gamma(J_\Gamma f)=\pd_{\Gamma^\prime}(eJ_\Gamma f)=\pd_{\Gamma^\prime}(e\Gamma)$.

By \XCBB, we obtain $\gl\Gamma^\prime-\gl\Gamma\le\pd_{\Gamma^\prime}(e\Gamma)=\pd_\Gamma(J_\Gamma f)\le\gl\Gamma-1$. On the other hand, since $\overline{\Gamma}$ is semisimple, we obtain $\gl\Gamma-\gl\Gamma^\prime\le\gl\overline{\Gamma}+\pd_\Gamma\overline{\Gamma}+1\le\pd_\Gamma(\Gamma e\Gamma)+2=\pd_\Gamma(\Gamma e\Gamma f)+2=\pd_\Gamma(J_\Gamma f)+2=\pd_{\Gamma^\prime}(e\Gamma)+2\le\gl\Gamma^\prime+2$, where we used $\Gamma e\Gamma f=J_\Gamma f$ since $\add_\Gamma(\Gamma e)$ is a cosemisimple subcategory of $\pr\Gamma$.\rule{5pt}{10pt}

\vskip.5em{\bf\XCE\ }
Assume that $\cc$ corresponds to an artin algebra or order $\Gamma$ by \XCA.

(1) We have a bijection between semisimple rejective subcategories $\cc^\prime$ of $\cc$ and heredity ideals $I$ of $\Gamma$. A heredity chain corresponds to a chain $\cc_m\subseteq\cc_{m-1}\subseteq\cdots\subseteq\cc_0=\cc$ of subcategories of $\cc$ such that $\cc_n/[\cc_{n+1}]$ is a semisimple rejective subcategory of $\cc/[\cc_{n+1}]$ for any $n$ ($0\le n<m$) and $\cc_m$ is maximal in the sense of \XBB.

(2) For a chain $0=\cc_m\subseteq\cc_{m-1}\subseteq\cdots\subseteq\cc_0=\cc$ of subcategories of $\cc$, take an idempotent $e_n$ of $\Gamma$ which corresponds to the subcategory $\cc_n$ of $\cc$ by \XCA. We call the chain {\it neat} if $(f_1,f_2,\cdots,f_m)$ is an neat sequence for $f_n:=e_{n-1}-e_n$. We have a bijection between neat sequences of $\Gamma$ and neat chains of $\cc$.

\vskip.5em{\bf\XCEA\ Theorem }{\it
(1) Any complete total right (resp. left) rejective chain is a heredity chain. As a conclusion, the categories appeared in \XBBC, \XBCC(1), \XBCD, \XBD, \XBDA, \XBDB, \XBEA\ and \XBEB\ are quasi-hereditary.

(2) Any complete half rejective chain is a neat chain. As a conclusion, the category appeared in \XBCC(2) is neat.}

\vskip.5em{\sc Proof }
(1) Let $0=\cc_m\subseteq\cc_{m-1}\subseteq\cdots\subseteq\cc_0=\cc$ be a total right rejective chain. By \XBA(1)(2), $\cc_{n_2}/[\cc_{n_3}]$ is a right rejective subcategory of $\cc_{n_1}/[\cc_{n_3}]$ for any $n_1\ge n_2\ge n_3$. In particular, the chain is heredity.

(2) Let $\cc^\prime$ be a cosemisimple right rejective subcategory of $\cc$ corresponding to an idempotent $e$ of $\Gamma$ by \XCA. We only have to show that $f:=1-e$ is an neat idempotent. By \XCB(3), $0\to J_\Gamma f\to\Gamma f\to(\Gamma/J_\Gamma)f\to0$ is a projective resolution of $(\Gamma/J_\Gamma)f$ with $J_\Gamma f\in\add_\Gamma(\Gamma e)$. This implies $\pd_\Gamma((\Gamma/J_\Gamma)f)\le 1$ and $\ext^1_\Gamma((\Gamma/J_\Gamma)f,(\Gamma/J_\Gamma)f)=0$.\rule{5pt}{10pt}

\vskip.5em{\bf\XCEB\ }Immediately we obtain the theorem below, where the artin algebra case was conjectured by Ringel [Xi1] and proved in [I4;1.1].

\vskip.5em{\bf Corollary }{\it
Let $\Lambda$ be an artin algebra or order in a semisimple algebra. For any $M\in\dn{M}_\Lambda$, there exists $N\in\dn{M}_\Lambda$ such that $\endm_\Lambda(M\oplus N)$ is quasi-hereditary.}

\vskip.5em{\bf\XCEC\ }
Let $\Gamma$ be an artin algebra or order. By \XCEA, the following implications hold, where (QH) menas that $\Gamma$ is quasi-hereditary, (NT) means that $\Gamma$ is neat, and (RC) (resp. (TRRC),(RRC),(HRC),...) means that $\Gamma$ has a rejective (resp. total right rejective, right rejective, half rejective,...) chain. We shall see in \XCED\ that all implications are strict.
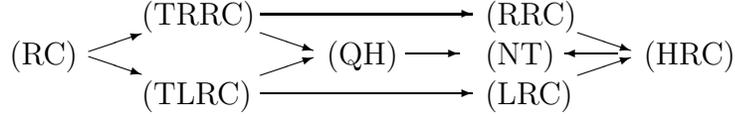
\[\begin{picture}(300,35)
\put(0,15){(RC)}
\put(30,20){\vector(3,1){20}}
\put(30,18){\vector(3,-1){20}}
\put(50,30){(TRRC)}
\put(50,0){(TLRC)}
\put(120,15){(QH)}
\put(180,15){(NT)}
\put(180,30){(RRC)}
\put(180,0){(LRC)}
\put(240,15){(HRC)}
\put(150,19){\vector(1,0){20}}
\put(230,19){\vector(-1,0){20}}
\put(95,34){\vector(1,0){80}}
\put(95,4){\vector(1,0){80}}
\put(95,27){\vector(3,-1){20}}
\put(95,11){\vector(3,1){20}}
\put(215,27){\vector(3,-1){20}}
\put(215,11){\vector(3,1){20}}
\end{picture}\]

\vskip0em{\bf\XCED\ Examples }
(1) Let $\Gamma$ be the algebra defined by the quiver below with relations $bc=ca=dab=0$. Then the indecomposable left and right projective modules have the form below. Then $\Gamma$ is (RRC) with a right rejective chain $\{3\}\subset\{2,3\}\subset\{1,2,3\}$, but not (QH)(LRC)(TRRC).
\[\begin{picture}(30,20)
\put(0,0){\vector(1,0){30}}
\put(0,4){\vector(1,0){30}}
\put(30,2){\vector(-1,1){15}}
\put(15,17){\vector(-1,-1){15}}
\put(-5,0){${\scriptstyle 3}$}
\put(30,0){${\scriptstyle 1}$}
\put(12,18){${\scriptstyle 2}$}
\put(20,15){${\scriptstyle a}$}
\put(5,15){${\scriptstyle b}$}
\put(12,5){${\scriptstyle c}$}
\put(12,-5){${\scriptstyle d}$}
\end{picture}
\ \ \ \ \ \ \ \ \ \ \ \ \ \ \ 
\left.\def\arraystretch{.5}\begin{array}{c}
{\scriptstyle 1}\\ {\scriptstyle 2}\\ {\scriptstyle 3}\\ {\scriptstyle\ \ \ 1}\\ {\scriptstyle\ \ \ 2}\end{array}\right.\ \ \ \ \ 
\left.\def\arraystretch{.5}\begin{array}{c}
{\scriptstyle 2}\\ {\scriptstyle 3}\\ {\scriptstyle\ \ \ 1}\\ {\scriptstyle\ \ \ 2}\end{array}\right.\ \ \ \ \ 
\left.\def\arraystretch{.5}\begin{array}{c}
{\scriptstyle 3}\\ {\scriptstyle 1\ \ 1}\\ {\scriptstyle\ \ \ 2}\end{array}\right.
\ \ \ \ \ \ \ \ \ \ \ \ \ \ \ 
\left.\def\arraystretch{.5}\begin{array}{c}
{\scriptstyle 1}\\ {\scriptstyle 3\ \ 3}\\ {\scriptstyle\ \ \ 2}\\ {\scriptstyle\ \ \ 1}\end{array}\right.\ \ \ \ \ 
\left.\def\arraystretch{.5}\begin{array}{c}
{\scriptstyle 2}\\ {\scriptstyle 1}\\ {\scriptstyle\ \ \ 3}\\ {\scriptstyle\ \ \ 2}\\ {\scriptstyle\ \ \ 1}\end{array}\right.\ \ \ \ \ 
\left.\def\arraystretch{.5}\begin{array}{c}
{\scriptstyle 3}\\ {\scriptstyle 2}\\ {\scriptstyle 1}\end{array}\right.\]

(2) Let $\Gamma$ be the algebra defined by the quiver below with relations $dc=bca=bcd=ca=0$. Then the indecomposable left and right projective modules have the form below. Thus $\Gamma$ is (QH) with a heredity chain $\{2\}\subset\{1,2\}\subset\{1,2,3\}$, but not (HRC) since it has no cosemisimple left (resp. right) rejective subcategory.
\[\begin{picture}(30,20)
\put(30,0){\vector(-1,0){30}}
\put(0,4){\vector(1,0){30}}
\put(30,2){\vector(-1,1){15}}
\put(15,17){\vector(-1,-1){15}}
\put(-5,0){${\scriptstyle 3}$}
\put(30,0){${\scriptstyle 1}$}
\put(12,18){${\scriptstyle 2}$}
\put(20,15){${\scriptstyle a}$}
\put(5,15){${\scriptstyle b}$}
\put(12,5){${\scriptstyle c}$}
\put(12,-5){${\scriptstyle d}$}
\end{picture}
\ \ \ \ \ \ \ \ \ \ \ \ \ \ \ 
\left.\def\arraystretch{.5}\begin{array}{c}
{\scriptstyle 1}\\ {\scriptstyle 2\ \ 3}\\ {\scriptstyle 3\ \ \ }\\ {\scriptstyle 1\ \ \ }\end{array}\right.\ \ \ \ 
\left.\def\arraystretch{.5}\begin{array}{c}
{\scriptstyle 2}\\ {\scriptstyle 3}\\ {\scriptstyle 1}\end{array}\right.\ \ \ \ 
\left.\def\arraystretch{.5}\begin{array}{c}
{\scriptstyle 3}\\ {\scriptstyle 1}\\ {\scriptstyle 3}\end{array}\right.
\ \ \ \ \ \ \ \ \ \ \ 
\left.\def\arraystretch{.5}\begin{array}{c}
{\scriptstyle 1}\\ {\scriptstyle 3}\\ {\scriptstyle 2}\\ {\scriptstyle 1}\end{array}\right.\ \ \ \ 
\left.\def\arraystretch{.5}\begin{array}{c}
{\scriptstyle 2}\\ {\scriptstyle 1}\end{array}\right.\ \ \ \ 
\left.\def\arraystretch{.5}\begin{array}{c}
{\scriptstyle 3}\\ {\scriptstyle 2\ \ 1}\\ {\scriptstyle 1\ \ 3}\end{array}\right.\]

(3) Let $k\subsetneq K$ be a finite extension of fields and $J:=tK[[t]]\subset\Delta:=k+J\subset\oo:=K[[t]]$. Let $\Gamma:=\def\arraystretch{.5}\left(\begin{array}{cccc}{\scriptstyle\oo}&{\scriptstyle\oo}&{\scriptstyle\oo}&{\scriptstyle\oo}\\ {\scriptstyle J}&{\scriptstyle\Delta}&{\scriptstyle\oo}&{\scriptstyle\oo}\\ {\scriptstyle J^2}&{\scriptstyle t\Delta}&{\scriptstyle\Delta}&{\scriptstyle\Delta}\\ {\scriptstyle J^2}&{\scriptstyle J^2}&{\scriptstyle J}&{\scriptstyle\Delta}\end{array}\right)$ be a $k[[t]]$-order. Then $\Gamma$ is (HRC) with a half rejective chain $\{1\}\subset\{1,2\}\subset\{1,2,3\}\subset\{1,2,3,4\}$, where $i$ corresponds to an idempotent with $(i,i)$-th entry $1$ and other entries $0$. But $\Gamma$ is neither (RRC) nor (LRC).

\vskip.5em{\bf\XCF\ }Any artin algebra $\Gamma$ with $\gl\Gamma\le2$ is quasi-hereditary by Dlab-Ringel [DR2], and any order $\Gamma$ with $\gl\Gamma\le2$ is also quasi-hereditary by K\"onig-Wiedemann [KW]. More strongly, we will that such ring has a complete total right rejective chain (\XCEA).

\vskip.5em{\bf Theorem }{\it 
Let $\Gamma$ be an artin algebra or order, and $\cc:=\pr\Gamma$. If $\gl\Gamma\le 2$ holds, then $\cc$ has a complete total right (resp. left) rejective chain.}

\vskip.5em{\bf\XCFA\ Lemma }{\it Let $\Gamma$ be an artin algebra with $\gl\Gamma=n$ ($2\le n<\infty$). Then there exists a simple $\Gamma$-module $S$ with $\pd_\Gamma S=n-1$.}

\vskip.5em{\sc Proof }
We can take a non-zero $\Gamma$-module $X$ with $\pd_\Gamma X=n-1$. Assume that $X$ is not simple. Thus $X$ has a proper simple submodule $S$, and we obtain an exact sequence $0\to S\to X\to X/S\to0$. By an exact sequence $\ext^{n}_\Gamma(X,\ )\to\ext^n_\Gamma(S,\ )\to\ext^{n+1}_\Gamma(X/S,\ )$, we obtain $\pd_\Gamma S<n$. If $\pd_\Gamma S<n-1$, then we obtain $\pd_\Gamma(X/S)=n-1$ by an exact sequence $\ext^{n-1}_\Gamma(X/S,\ )\to\ext^{n-1}_\Gamma(X,\ )\to\ext^{n-1}_\Gamma(S,\ )\to\ext^n_\Gamma(X/S,\ )\to\ext^{n}_\Gamma(X,\ )$. Replacing $X$ by $X/S$ and repeating this argument, we obtain the assertion.\rule{5pt}{10pt}

\vskip.5em{\bf\XCFB\ Proof of \XCF\ }
We can assume $\Gamma$ is non-semisimple basic. By \XCFA, we can take a primitive idempotent $f$ of $\Gamma$ such that $\pd_\Gamma S=1$ for $S:=(\Gamma/J_\Gamma)f$. Put $e:=1-f$, $\cc^\prime:=\add_\Gamma(\Gamma e)\subset\cc=\pr\Gamma$ and $\Gamma^\prime:=e\Gamma e$. Then $\cc^\prime$ is a cosemisimple right rejective subcategory of $\cc$ and $\gl\Gamma^\prime\le 2$ holds by \XCB(3). By \XBA(4), we only have to show that any monomorphism $a:P_1\to P_0$ in $\cc^\prime$ is still monic in $\cc$. Since $\gl\Gamma\le2$, we have an exact sequence $0\to P_2\to P_1\stackrel{a}{\to}P_0$ in $\mod\Gamma$ with $P_2\in\pr\Gamma$. Since $a$ is monic in $\cc^\prime$, $eP_2=\hom_\Gamma(\Gamma e,P_2)=0$ holds. Thus $P_2$ is a module over a semisimple ring $\Gamma/\Gamma e\Gamma$. This implies that $P_2$ is isomorphic to $S^n$ for some $n\ge0$. If $n>0$, then $S\in\pr\Gamma$, a contradiction to $\pd_\Gamma S=1$. Thus $n=0$ and $P_2=0$ holds.\rule{5pt}{10pt}

\vskip.5em{\bf\XCG\ Example }Let $\Lambda$ be an artin algebra of finite representation type. We call a chain $0=\cc_m\subseteq\cc_{m-1}\subseteq\cdots\subseteq\cc_0=\mod\Lambda$ a {\it splitting chain} if $\cc_{n+1}$ is closed under either factor modules or submodules, and is a cosemisimple subcategory of $\cc_n$ for any $n$ ($0\le n<m$). Any splitting filtration of Dlab-Ringel [DR4] gives a splitting chain. Now, \XAEB\ implies that any $\cc_n$ is a left or right rejective subcategory of $\mod\Lambda$. Thus \XBA(1) implies that any splitting chain is a half rejective chain. On the other hand, \XBA(2)(5) implies that $\cc_n/[\cc_{n+1}]$ is a semisimple rejective subcategory of $\mod\Lambda/[\cc_{n+1}]$ for any $n$. By \XCE(1), we conclude that any splitting chain is a heredity chain, a theorem of Dlab-Ringel [DR4].

\vskip.5em{\bf\XD\ The function $r_\Lambda$ and Representation dimension}

\vskip.5em{\bf\XDA\ Definition }
Again assume that $\Lambda$ is an artin algebra and $\dn{M}_\Lambda:=\mod\Lambda$, or $\Lambda$ is an order and $\dn{M}_\Lambda:=\lat\Lambda$. Put $g_\Lambda(M):=\gl\endm_\Lambda(M)$ and $r_\Lambda(M):=\inf\{g_\Lambda(M\oplus N)\ |\ N\in\dn{M}_\Lambda\}$. By \XDAA\ below (cf. [I4;1.2]), $r_\Lambda$ gives a function $\dn{M}_\Lambda\to\nnn_{\ge0}$. We often regard $g_\Lambda$ and $r_\Lambda$ as functions on the set of finite subcategories of $\dn{M}_\Lambda$. Then $r_\Lambda(\cc^\prime)\le r_\Lambda(\cc)$ holds for finite subcategories $\cc^\prime\subseteq\cc$. Put $|r_\Lambda|:=\sup\{ r_\Lambda(X)\ |\ X\in\dn{M}_\Lambda\}$. We call $\rdim\Lambda:=r_\Lambda(\Lambda\oplus D\Lambda)$ the {\it representation dimension} of $\Lambda$, which was introduced by Auslander [A1] for artin algebras.

\vskip.5em{\bf\XDAA\ Theorem }{\it
Let $\Lambda$ be an artin algebra or order in a semisimple algebra. Then $r_\Lambda(X)<\infty$ holds for any $X\in\dn{M}_\Lambda$.}

\vskip1em{\sc Proof }Immediate from \XBEB.\rule{5pt}{10pt}

\vskip.5em{\bf\XDAB\ }
For a subcategory $\cc$ of $\dn{M}_\Lambda$, recall that we write $\resdim{\cc}{\dn{M}_\Lambda}\le m$ if, for any $X\in\dn{M}_\Lambda$, there exists a complex $M_{m}\to\cdots\to M_1\to M_0\to X$ such that $M_i\in\cc$ and $0\to\dn{M}_\Lambda(\ ,M_{m})\to\cdots\to\dn{M}_\Lambda(\ ,M_1)\to\dn{M}_\Lambda(\ ,M_0)\to\dn{M}_\Lambda(\ ,X)\to0$ is exact on $\cc$ (\XAA). Let us recall \XAAB\ below, which gives us a method to calculate $r_\Lambda$ [A1][EHIS].

\vskip.5em{\bf Theorem }{\it
Let $\Lambda$ be an artin algebra or order in a semisimple algebra, $M\in\dn{M}_\Lambda$ and $\cc:=\add M$. Then $\resdim{\cc^{op}}{\dn{M}_\Lambda^{op}}\ge\max\{g_\Lambda(M)-2,0\}\le\resdim{\cc}{\dn{M}_\Lambda}$, where the left (resp. right) equality holds if $\Lambda\in\cc$ (resp. $D\Lambda\in\cc$).}

\vskip.5em{\bf\XDAC\ }Immediately, we obtain Auslander's theorem [A1] below. It is easily checked that $\rdim\Lambda\le1$ and $|r_\Lambda|\le1$ occurs only when $\Lambda$ is semisimple (algebra case) or $\Lambda$ is hereditary (order case).

\vskip.5em{\bf Theorem }{\it
Let $\Lambda$ be an artin algebra or order in a semisimple algebra. Then $\rdim\Lambda\le 2$ if and only if $|r_\Lambda|\le 2$ if and only if $\Lambda$ is representation-finite (\S\XZC).}

\vskip.5em{\sc Proof }
Assume that a finite subcategory $\cc$ of $\dn{M}_\Lambda$ satisfies $\Lambda\oplus D\Lambda\in\cc$. By \XDAB, $g_\Lambda(\cc)\le 2$ is equivalent to $\resdim{\cc}{\dn{M}_\Lambda}=0$ which means that, for any $X\in\dn{M}_\Lambda$, there exists a morphism $Y\stackrel{f}{\to}X$ such that $Y\in\cc$ and $\dn{M}_\Lambda(\ ,Y)\stackrel{\cdot f}{\to}\dn{M}_\Lambda(\ ,X)$ is an isomorphism on $\cc$. This is equivalent to that $f$ is an isomorphism by $\Lambda\in\cc$. Thus $g_\Lambda(\cc)\le 2$ is equivalent to $\cc=\dn{M}_\Lambda$. Thus we obtain the assertion.\rule{5pt}{10pt}

\vskip.5em{\bf\XDAD\ Example }
Let $\Lambda$ be an artin algebra with the Loewy length $\LL(\Lambda)$.

(1) $r_\Lambda(\Lambda)\le g_\Lambda(\bigoplus_{i=0}^{\LL(\Lambda)}\Lambda/J_\Lambda^i)\le\LL(\Lambda)$ holds by \XBBC(1). In particular, $\rdim\Lambda\le\LL(\Lambda)$ holds if $\Lambda$ is selfinjective [A1].

(2) If $\Lambda$ is hereditary, then $\rdim\Lambda\le g_\Lambda(\Lambda\oplus D\Lambda)\le 3$ holds [A1]. On the other hand, if $J_\Lambda^2=0$, then $\Lambda$ is stably equivalent to a hereditary algebra [ARS], and we can obtain Auslander's result $\rdim\Lambda\le3$ [A1] by \XDB(2) below. Several classes of algebras are known to have the representation dimension at most three [Xi1][CP][H] (see also \XDD).

\vskip.5em{\sc Proof }
(2) By \XDAB, we only have to show that $\cc:=\add(\Lambda\oplus D\Lambda)$ satisfies $\resdim{\cc}{\dn{M}_\Lambda}\le 1$. Take $X\in\ind\dn{M}_\Lambda\backslash\ind\cc$. Then a projective resolution $0\to P_1\to P_0\to X\to0$ gives a right $\cc$-resolution since $\hom_\Lambda(D\Lambda,X)=0$ holds by $\gl\Lambda\le1$.\rule{5pt}{10pt}

\vskip.5em{\bf\XDB\ }Let $\Lambda$ be an artin algebra or order and $\Gamma$ an artin algebra or order. We say that $\Lambda$ is {\it finitely equivalent} to $\Gamma$ if there exists finite subcategories (\XAA) $\xx$ and $\xx^\prime$ of $\dn{M}_\Lambda$ and $\dn{M}_\Gamma$ respectively such that $\dn{M}_\Lambda/[\xx]$ is equivalent to $\dn{M}_\Gamma/[\xx^\prime]$. Especially, when $\xx=[\pr\Lambda]$ and $\xx^\prime=\pr\Gamma$, we say that $\Lambda$ is {\it stably equivalent} to $\Gamma$ [ARS]. We shall prove the theorem below in \XEF, where (2) is a result of Xiangqian [X]. It is also valid even if $\Lambda$ is an artin algebra and $\Gamma$ is an order.

\vskip.5em{\bf Theorem }{\it
Let $\Lambda$ be an artin algebra or order and $\Gamma$ an artin algebra or order. Assume that they are not representation-finite.

(1) If $\Lambda$ is finitely equivalent to $\Gamma$, then $|r_\Lambda|=|r_\Gamma|$.

(2) If $\Lambda$ is stably equivalent to $\Gamma$, then $\rdim\Lambda=\rdim\Gamma$.}

\vskip.5em{\bf\XDC\ }
Let $\Lambda$ be a hereditary artin algebra with the underlying valued quiver $Q$. Then $\Lambda$ is representation-finite if and only if $Q$ is Dynkin. We call $\Lambda$ {\it tame} if $Q$ is extended Dynkin, and {\it wild} otherwise [DR1][R1]. When $\Lambda$ is tame, we call $\endm_\Lambda(T)$ a {\it tame concealed} algebra for any preprojective tilting $\Lambda$-module $T$ [R2;4.3].

\vskip.5em{\bf Theorem }{\it
If $\Lambda$ is a tame hereditary algebra, then $|r_\Lambda|=\rdim\Lambda=3$. More generally, if $\Gamma$ is a tame concealed algebra, then $|r_\Gamma|=\rdim\Gamma=3$.}

\vskip.5em{\sc Proof }
(1) We will show $|r_\Lambda|\le 3$. For $X\in\mod\Lambda$, we denote by $\sub X$ the subcategory of $\mod\Lambda$ consisting submodules of a direct sum of copies of $X$. We denote by $\pp$ (resp. $\rr$, $\ii$) the category of preprojective (resp. regular, preinjective) $\Lambda$-modules [DR1]. Put $d(X,Y):=\dim\hom_\Lambda(X,Y)$ for $X,Y\in\mod\Lambda$. We denote by $\tau^+_\Lambda$ the Auslander translate with the inverse $\tau^-_\Lambda$ [ARS].

(i) We will show $\#\ind(\rr\cap\sub X)<\infty$ and $\#\ind(\pp\cap\sub X)<\infty$ for any $X\in\rr$.

The former assertion is easy since $\rr$ is a direct sum of subcategories corresponding to tubes [DR1][R1]. We shall show the latter one. We can take $p>0$ such that $X=\tau_\Lambda^{+p}X$. Put $c:=\max\{ d(Y,X)\ |\ Y\in\ind(\add(\bigoplus_{i=0}^{p-1}\tau_\Lambda^{-i}\Lambda))\}$. Then $d(Y,X)\le c$ holds for any $Y\in\ind(\pp\cap\sub X)$ since $d(Y,X)=d(\tau^+_\Lambda Y,\tau^+_\Lambda X)$ holds for any $Y\in\ind(\mod\Lambda)\backslash\pr\Lambda$. Fix $Y\in\ind(\pp\cap\sub X)$ and take a basis $\{ f_i\}_{1\le i\le d(Y,X)}$ of $\hom_\Lambda(Y,X)$. Then $(f_1,\cdots,f_{d(Y,X)}):Y\to X^{d(Y,X)}$ should be injective. Thus $\dim Y\le d(Y,X)\dim X\le c\dim X$ holds. Since there are only finitely many $Y\in\ind\pp$ with $\dim Y\le c\dim X$, the assertion follows.

(ii) Fix any $P\oplus N\oplus I\in\mod\Lambda$ with $P\in\pp$, $N\in\rr$ and $I\in\ii$. We can put $\rr\cap\sub N=\add N^\prime$ and $\pp\cap\sub N^\prime=\add Q$ by (i). Let $P^\prime$ (resp. $I^\prime$) be the direct sum of $X\in\ind\pp$ (resp. $X\in\ind\ii$) such that $\hom_\Lambda(X,Q)\neq0$ (resp. $\hom_\Lambda(I,X)\neq0$). By \XDAB, we only have to show $\resdim{\cc}{(\mod\Lambda)}\le 1$ for $\cc:=\add(P^\prime\oplus N^\prime\oplus I^\prime)$.

Fix any $X\in\ind(\mod\Lambda)$. We can assume $X\notin\add I^\prime$, so $\hom_\Lambda(I^\prime,X)=0$ holds. Consider the following commutative diagram of exact sequences, where $f$ is a right $(\add N^\prime)$-approximation of $X$ and $g$ is a right $(\add P^\prime)$-approximation of $X$.
\[\begin{diag}
0&\RA{}&Y&\RA{}&N_0&\RA{f}&X\\
&&\parallel&&\uparrow&&\uparrow^g\\
0&\RA{}&Y&\RA{}&Z&\RA{}&P_0
\end{diag}\]

Then $Y\in\sub N^\prime\subseteq\add(P^\prime\oplus N^\prime)$ holds by our choice. Since any indecomposable direct summand $W$ of $Z$ satifsies $W\in\add Y$ or $\hom_\Lambda(W,P_0)\neq0$. In both cases, $W\in\add(P^\prime\oplus N^\prime)$ holds by our choice. Thus $Z\in\add(P^\prime\oplus N^\prime)$ holds, and $0\to Z\to P_0\oplus N_0\stackrel{{g\choose f}}{\to}X$ gives a right $\cc$-resolution of $X$.

(2) We will show $|r_\Gamma|\le 3$. Let $\Gamma=\endm_\Lambda(T)$ for a preprojective tilting $\Lambda$-module $T$, which defines a torsion theory $(\tt,\ff)$ on $\mod\Lambda$ and a splitting torsion theory $(\xx,\yy)$ on $\mod\Gamma$. Let $\cc$ be a subcategory of $\mod\Lambda$ such that $\ind\cc=\ind(\mod\Lambda)\backslash\ind\tt$. Then $\cc$ is finite [R2]. By Tilting theorem, $\xx$ is equivalent to $\ff\subseteq\cc$ and $\yy=\mod\Gamma/[\xx]$ is equivalent to $\tt=\mod\Lambda/[\cc]$. Thus $\Lambda$ and $\Gamma$ are finitely equivalent, and $|r_\Gamma|=|r_\Lambda|=3$ holds by (1) and \XDB(1).\rule{5pt}{10pt}

\vskip.5em{\bf\XDD\ }The theorem of Erdmann-Holm-Schr\"oer and the author [EHIS] is generalized as follows. For the convenience of readers, we shall give a proof. For simplicity, put $\cc+\dd:=\add(\cc\cup\dd)$ for subcategories $\cc$ and $\dd$ of $\dn{M}_\Lambda$.

\vskip.5em{\bf Theorem }{\it
Let $\Lambda\stackrel{\phi}{\subset}\Gamma$ be artin algebras or orders. Assume that $\Gamma$ is representation-finite and $J_\Lambda$ is an ideal (resp. left ideal, right ideal) of $\Gamma$. Then $\rdim\Lambda\le 3$ (resp. $r_\Lambda(\Lambda)\le 3$, $r_\Lambda(D\Lambda)\le 3$) holds.}

\vskip.5em{\sc Proof }
We shall use the notations in \XAC. Put $\cc:=\pr\Lambda+\xx_\phi+\rin\Lambda$. (For the case when $J_\Lambda$ is a left (resp. right) ideal of $\Gamma$, put $\cc:=\pr\Lambda+\xx_\phi$ (resp. $\cc:=\xx_\phi+\rin\Lambda$) and apply a similar argument). By \XDAB, we only have to show $\resdim{\cc}{\dn{M}_\Lambda}\le 1$. Fix $X\in\ind\dn{M}_\Lambda\backslash\ind\cc$. By \XAG(1), a right $\xx_\phi$-approximation of $X$ gives an exact sequence $0\to\hom_\Lambda(C_\phi,X)\to X^-\stackrel{\epsilon^-_X}{\to}X\to\ext^1_\Lambda(C_\phi,X)$. Since $C_\phi J_\Lambda=0$ holds, $\ext^i_\Lambda(C_\phi,X)$ is a semisimple $\Lambda$-module. Taking the projective cover of $\Cok\epsilon^-_X$, we obtain the following pull-back diagram.
\[\begin{diag}
0&\RA{}&\hom_\Lambda(C_\phi,X)&\RA{}&X^-&\RA{\epsilon^-_X}&X&\RA{}&\Cok\epsilon^-_X&\RA{}&0\\
&&\parallel&&\uparrow&&\uparrow^{f}&&\parallel\\
0&\RA{}&\hom_\Lambda(C_\phi,X)&\RA{}&Y&\RA{}&P&\RA{}&P/J_\Lambda P&\RA{}&0
\end{diag}\]

Then $f$ induces a morphism $g:J_\Lambda P\to\Im\epsilon^-_X$. Since $J_\Lambda P\in\xx_\phi$ holds by our assumption, $g$ factors through $X^-$. Thus $Y$ is isomorphic to $J_\Lambda P\oplus \hom_\Lambda(C_\phi,X)\in\xx_\phi$. Taking the mapping cone, we obtain a right $(\pr\Lambda+\xx_\phi)$-resolution $0\to Y\to P\oplus X^-\stackrel{{f\choose\epsilon^-_X}}{\to}X\to0$ $(*)$. Now take any $a\in\hom_\Lambda(D\Lambda,X)$. Then $Da\in\hom_{\Lambda^{op}}(DX,\Lambda)$ factors through $J_\Lambda$ since $DX$ is not projective. Thus $a$ factors through the natural map $D\Lambda\to DJ_\Lambda$. Since $DJ_\Lambda\in\xx_\phi$ holds by our assumption, $a$ factors through $\epsilon^-_X$. Thus $(*)$ is a right $\cc$-resolution.\rule{5pt}{10pt}

\vskip.5em{\bf\XDDA\ }Let us give a proof of [EHIS;1.4] from our categorical viewpoint.

\vskip.5em{\bf Corollary }{\it
Let $\Lambda\stackrel{\phi}{\subset}\Gamma$ be artin algebras with $J_\Lambda=J_\Gamma$. Assume that $\Gamma$ is representation-finite and $N$ is an additive generator of $\mod\Gamma$. Then $\endm_\Lambda(\Lambda\oplus\phi^*N\oplus D\Lambda)$ is quasi-hereditary.}

\vskip.5em{\sc Proof }
We have a chain $\cc^\prime:=\mod\Lambda/J_\Lambda\subseteq\cc_2:=\xx_\phi\subseteq\cc_1:=\xx_\phi+\rin\Lambda\subseteq\cc:=\pr\Lambda+\xx_\phi+\rin\Lambda$. By \XAFA(1), $\cc^\prime$ is a semisimple rejective subcategory of $\cc$. By \XAG(3), $\cc_2/[\cc^\prime]$ is equivalent to $\mod\Gamma/[\mod\Gamma/J_\Gamma]$. By \XCB(2), $\gl(\mod\Gamma)\le 2$ implies $\gl(\mod\cc_2/[\cc^\prime])\le 2$. By \XCF, $\cc_2$ has a heredity chain $0\subseteq\cc^\prime=\cc_m\subseteq\cdots\subseteq\cc_3\subseteq\cc_2$. Since $\cc_2/[\cc^\prime]$ is a rejective subcategory of $\cc/[\cc^\prime]$ by \XAG(3), the induced chain $0\subseteq\cc_{m-1}/[\cc^\prime]\subseteq\cdots\subseteq\cc_3/[\cc^\prime]\subseteq\cc_2/[\cc^\prime]$ consists of rejective subcategories of $\cc/[\cc^\prime]$ by \XBA(3). Now $\cc_2$ is a cosemisimple (left rejective) subcategory of $\cc_1$ by $DJ_\Lambda\in\cc_2$, and $\cc_1$ is a cosemisimple right rejective subcategory of $\cc$ by $J_\Lambda\in\cc_2$. Thus $\cc_1/[\cc_2]$ is a semisimple rejective subcategory of $\cc/[\cc_2]$ by \XBA(2)(5). Consequently, $\cc$ has a heredity chain $0\subseteq\cc^\prime=\cc_m\subseteq\cdots\subseteq\cc_3\subseteq\cc_2\subseteq\cc_1\subseteq\cc$.\rule{5pt}{10pt}

\vskip.5em{\bf\XDDB\ Example }
(1) An order $\Lambda$ is called {\it B\"ackstr\"om} if there exists a hereditary overorder $\Gamma$ such that $J_\Lambda=J_\Gamma$. By \XDD, $\rdim\Lambda\le3$ folds for any B\"ackstr\"om orders. The representation theory of B\"ackstr\"om orders were studied by Ringel-Roggenkamp [RR], and it is shown that the representation type of $\Lambda$ can be determined by a certain associated valued quiver of $\Lambda$. This result can be explained by the result \XDDC\ below in [I8] since the associated valued quiver of $\Lambda$ is defined as the associated valued quiver of the hereditary algebra $\def\arraystretch{.5}\left(\begin{array}{cc}{\scriptstyle \Gamma/J_\Gamma}&{\scriptstyle \Gamma/J_\Gamma}\\ {\scriptstyle 0}&{\scriptstyle \Lambda/J_\Lambda}\end{array}\right)$.

(2) In [EHIS], the representation dimension of a special biserial algebras is shown to be at most $3$ by applying \XDD. Now let us consider the representation dimension of a clannish algebra $A$ over a field $k$ [CB]. Then there exists a B\"ackstr\"om $k[[t]]$-order $\Lambda$ with a hereditary overorder $\Gamma$ and an ideal $I$ of $\Gamma$ such that with $J_\Lambda=J_\Gamma$ and $A=\Lambda/I$ by [I7;1.3(3)]. Since $B:=\Gamma/I$ is a cyclic Nakayama algebra with $J_A=J_B\subset A\subset B$, we obtain that $B$ is representation-finite and $\rdim A\le3$ by \XDD.

\vskip.5em{\bf\XDDC\ Theorem }{\it
Let $\Lambda$ be a B\"ackstr\"om order and $\Gamma$ the hereditary overorder of $\Lambda$ such that $J_\Lambda=J_\Gamma$. Then $\Lambda$ is stably equivalent to the finite dimensional hereditary algebra $\def\arraystretch{.5}
\left(\begin{array}{cc}{\scriptstyle \Gamma/J_\Gamma}&{\scriptstyle \Gamma/J_\Gamma}\\ {\scriptstyle 0}&{\scriptstyle \Lambda/J_\Lambda}\end{array}\right)$.}

\vskip.5em{\bf\XDE\ Definition }
We shall introduce two homological invariant of $\Lambda$ which is closely related to the function $r_\Lambda$. Let $\Lambda$ and $\Gamma$ be artin algebras or orders.

(1) We write $\Lambda\preceq\Gamma$ if there exists $P\in\pr\Gamma$ such that $\endm_\Gamma(P)$ is Morita-equivalent to $\Lambda$. Obviously, $\preceq$ gives a partial order on the set of Morita-equivalence classes of artin algebras and those of orders. Define the {\it expanded dimension} of $\Lambda$ by $\expdim\Lambda:=\inf\{\gl\Gamma\ |\ \Lambda\preceq\Gamma\}$. This concept first appeared in Auslander's observation in [A1] such that $\expdim\Lambda$ is finite for any artin algebra $\Lambda$ by \XBBC(1).

(2) Let $\cc$ be a subcategory of $\dn{M}_\Lambda$. We define the {\it weak resolution dimension} $\wresdim\cc$ as the minimal number $n\ge0$ which satisfies the following equivalent conditions (cf. \XAA).

\strut\kern1em(i) There exists $M\in\dn{M}_\Lambda$ such that, for any $X\in\cc$, there exists an exact sequence $0\to M_n\to\cdots\to M_0\to Y\to0$ with $M_i\in\add M$ and $X\in\add Y$.

\strut\kern1em(ii) There exists $M\in\dn{M}_\Lambda$ such that, for any $X\in\cc$, there exists an exact sequence $0\to Y\to M_0\to\cdots\to M_n\to0$ with $M_i\in\add M$ and $X\in\add Y$.

Here $\wresdim(\mod\Lambda)+2$ coincides with $\rwrdim\Lambda$ in [Rou]. 

\vskip.5em{\sc Proof }
(2)(i)$\Rightarrow$(ii) We shall show that $N:=\Omega^{-n}M\oplus D\Lambda$ satisfies the condition (ii). Put $M_{-1}:=Y$ for simplicity. Assume that we have an exact sequence $0\to M_{n-i}^\prime\to\cdots\to M_{-1}^\prime\to N_{-2}\to\cdots\to N_{-i-1}\to0$ with $M_j^\prime\in\add(M_j\oplus D\Lambda)$ and $N_j\in\add N$. Consider the following commutative diagram, where the lower sequence is an injective resolution.
\[\begin{diag}
0&\to&M_{n-i}^\prime&\to&M_{n-i-1}^\prime&\to&\cdots&\to&M_{-1}^\prime&\to&N_{-2}&\to&\cdots&\to&N_{-i}&\to&N_{-i-1}&\to&0\\
&&\parallel&&\downarrow&&&&\downarrow&&\downarrow&&&&\downarrow&&\downarrow&&\\
0&\to&M_{n-i}^\prime&\to& I_{n-i-1}&\to&\cdots&\to&I_{-1}&\to&I_{-2}&\to&\cdots&\to&I_{-i}&\to&\Omega^{-n}M_{n-i}^\prime&\to&0
\end{diag}\]

Taking the mapping cone, we obtain an exact sequence $0\to M_{n-i-1}^\prime\to\cdots\to M_{-1}^\prime\to N_{-2}\to\cdots\to N_{-i-2}\to0$ with $M_j^\prime\in\add(M_j\oplus D\Lambda)$ and $N_j\in\add N$. The sequence for $i=n+1$ shows (ii).\rule{5pt}{10pt}

\vskip.5em{\bf\XDEA\ }
Let $\Lambda$ be an artin algebra or order.

(1) For $n=0,1$, $\expdim\Lambda=n$ if and only if $r_\Lambda(\Lambda)=n$ if and only if $\gl\Lambda=n$.

(2) $\expdim\endm_\Lambda(X)\le r_\Lambda(X)$ holds for any $X\in\dn{M}_\Lambda$.

(3) $\wresdim(\mod\Lambda)\le\expdim\Lambda\le r_\Lambda(\Lambda)\le\min\{\gl\Lambda,\rdim\Lambda,\LL(\Lambda)\}$.

(4) $\wresdim\xx_2\le\max\{\expdim\Lambda-2,0\}$ holds for $\xx_2:=\add\Omega^2(\mod\Lambda)$.

(5) If $\phi:\Lambda\to\Gamma$ is a quasi-split morphism of artin algebras (\XAF), then $\wresdim\dn{M}_\Lambda\le\wresdim\dn{M}_\Gamma$.

\vskip.5em{\sc Proof }(1) If $\Gamma\preceq\Lambda$ and $\gl\Lambda\le 1$, then $\gl\Gamma\le 1$.

(2) For any $Y$, put $\Gamma:=\endm_\Lambda(X\oplus Y)$ and $P:=\hom_\Lambda(X\oplus Y,X)\in\pr\Gamma$. Since $\endm_{\Gamma}(P)=\endm_\Lambda(X)$ holds, we obtain $\endm_\Lambda(X)\preceq\Gamma$.

(4) Put $n:=\expdim\Lambda$ and $m:=\max\{n-2,0\}$. Take $\Gamma$ and $P\in\pr\Gamma$ with $\gl\Gamma=n$ and $\Lambda=\endm_\Lambda(P)$. We will show that $M:=\hom_\Gamma(P,\Gamma)$ satisfies \XDE(2)(i). For any $X\in\xx_2$, take a morphism $f\in\hom_\Gamma(P^\prime,P^{\prime\prime})$ in $\add P$ such that $0\to Y\to\hom_\Gamma(P,P^\prime)\stackrel{\cdot f}{\to}\hom_\Gamma(P,P^{\prime\prime})$ is exact and $X\in\add Y$. By $\pd_\Gamma\Ker f\le m$, we can take a projective resolution $0\to P_m\to\cdots\to P_0\to\Ker f\to0$. Then we obtain the desired sequence $0\to\hom_\Gamma(P,P_m)\to\cdots\to\hom_\Gamma(P,P_0)\to Y\to0$.

(3) The inequality $\wresdim\mod\Lambda\le\expdim\Lambda$ follows from the argument in the proof of (4). Other inequalities follows from (2) and \XDAD(1).

(5) If $M\in\dn{M}_\Gamma$ satisfies the condition in \XDE(2)(i), then so does $\phi^*M\in\dn{M}_\Lambda$.\rule{5pt}{10pt}

\vskip.5em{\bf\XDF\ }Rouquier's theorem below [Rou;6.10,6.9] gave the first example of algebras with representation dimension greater than $3$. Although his quite interesting proof uses Koszul duality in derived categories essentially, it would be a natural problem to give a direct proof. Consequently, $\expdim\Lambda$ and $r_\Lambda(\Lambda)$ is also $n+1$ by \XDEA(3)(4).

\vskip.5em{\bf Theorem }{\it
Let $k$ be a field and $\Lambda=\wedge(k^n)$ the exterior algebra with $n>0$. Then $\wresdim(\mod\Lambda)+2=\rdim\Lambda=n+1$.}

\vskip.5em{\bf\XDFA\ Example }We can obtain many classes of artin algebras with large representation dimension by using \XDF. Again let $\Lambda=\wedge(k^n)$ be the exterior algebra with $n>0$.

(1) An artin algebra $\Gamma$ with $\Lambda\preceq\Gamma$ satisfies $n+1=\expdim\Lambda\le\expdim\Gamma\le\rdim\Gamma$ by \XDEA(3).

(2) If $\phi:\Lambda\to\Gamma$ is a quasi-split morphism of artin algebras, then $n-1=\wresdim(\mod\Lambda)\le\wresdim(\mod\Gamma)\le\rdim\Gamma$ by \XDEA(5).

\vskip.5em{\bf\XDFB\ }The concept of controlled wild algebras was introduced in [Han]. Any controlled wild algebra is wild, and Ringel conjectured that the converse holds. It is known that any wild hereditary algebra is controlled wild.

\vskip.5em{\bf Theorem }{\it
If $\Lambda$ is a controlled wild algebra, then $|r_\Lambda|=\infty$.}

\vskip.5em{\sc Proof }
For any artin algebra $\Gamma$, there exists $M\in\mod\Lambda$ and an ideal $I$ of $\endm_\Lambda(M)$ such that $\endm_\Lambda(M)=\Gamma\oplus I$ [Han;2.3]. Considering the case $\Gamma$ is the exterior algebra over $n$-dimensional vector space, we obtain $n-1\le\wresdim(\mod\endm_\Lambda(M))\le\expdim\endm_\Lambda(M)\le r_\Lambda(M)\le |r_\Lambda|$ by \XDFA(2) and \XDEA(2)(3). Thus $|r_\Lambda|=\infty$ holds.\rule{5pt}{10pt}

\vskip.5em{\bf\XDFC\ Corollary }{\it
Let $\Lambda$ be a finite dimensional hereditary algebra or a B\"ackstr\"om order (\XDDB). Then the value of $|r_\Lambda|$ is given as follows.
{\scriptsize \[\begin{array}{|c||c|c|c|}\hline
\mbox{\rm associated valued quiver}&\mbox{\rm Dynkin}&\mbox{\rm extended Dynkin}&\mbox{\rm else}\\ \hline
|r_\Lambda|&\le 2&3&\infty\\ \hline
\end{array}\]}}

\vskip-1em{\sc Proof }For the hereditary case, the assertion follows from \XDAB, \XDC\ and \XDFB. Thus \XDB\ and \XDDC\ shows the B\"ackstr\"om case.\rule{5pt}{10pt}

\vskip.5em{\bf\XDFD\ Question }Does any wild algebra satify $|r_\Lambda|=\infty$? By \XDFB, this is true if Ringel's conjecture is true. On the other hand, does any tame algebra satify $|r_\Lambda|<\infty$? Also, it is an interesting question whether any tame algebra satifies $\rdim\Lambda\le3$ or not. These questions can be regarded as a part of the study of tame algebras in terms of endomorphism rings.

\vskip.5em{\bf\XDG\ }
Let $\Lambda$ be an artin algebra or order and $n\ge0$. We call $\Lambda$ {\it reflexive-finite} if $\Lambda$ has only finitely many indecomposable reflexive modules. Put $\xx_n:=\add\Omega^n(\mod\Lambda)$. Recall that $X\in\mod\Lambda$ is called {\it $n$-torsionfree} if $\ext^i_\Lambda(\tr X,\Lambda)=0$ holds for any $i$ ($0<i\le n$) [AB]. It is well-known that $X$ is $1$-torsionfree (resp. $2$-torsionfree) if and only if it is torsionless (resp. reflexive). One can easily check that any $n$-tosionfree module is contained in $\Omega^n(\mod\Lambda)$. Conversely, it is known that $\xx_n$ consists of $n$-torsionfree modules if $\Lambda$ is an $(n-1)$-Gorenstein ring [AR3,4][FGR], for example.

\vskip.5em{\bf\XDGA\ Theorem }{\it
Let $\Lambda$ be an artin algebra or order, and $n\ge2$. Assume that $\xx_n$ has an additive generator $M$ that is $n$-torsionfree. Then $g_\Lambda(M)\le n$ and $r_\Lambda(\Lambda)\le n$ hold.}

\vskip.5em{\sc Proof }Put $\Gamma:=\endm_\Lambda(M)$. For any $X\in\mod\Gamma$, take a projective resolution ${\bf P}:0\to \Omega^nX\to P_{n-1}\to\cdots\to P_0\to X\to0$. Then there exists a complex ${\bf A}:M_{n-1}\to\cdots\to M_0$ is $\mod\Lambda$ in $\add_\Lambda M$ and $\hom_\Lambda(M,{\bf A})$ gives $P$. By $\Lambda\in\add M$, ${\bf A}$ is exact. Since each $M_i$ is $n$-torionfree, there exists an exact sequence $0\to M_i\to P_{i,0}\to\cdots\to P_{i,n}$ which is a left $(\pr\Lambda)$-resolution. Thus we can consider the following commutative diagram.
{\scriptsize\[\begin{diag}
0&\RA{}&\Ker f&\RA{}&M_{n-1}&\RA{f}&M_{n-2}&\RA{}&\cdots&\RA{}&M_1&\RA{}&M_0\\
&&&&\downarrow&&\downarrow&&&&\downarrow&&\downarrow&&\\
&&&&P_{n-1,0}&\RA{}&P_{n-2,0}&\RA{}&\cdots&\RA{}&P_{1,0}&\RA{}&P_{0,0}&&\\
&&&&\downarrow&&\downarrow&&&&\downarrow&&&&\\
&&&&P_{n-1,1}&\RA{}&P_{n-2,1}&\RA{}&\cdots&\RA{}&P_{1,1}&&\\
&&&&\downarrow&&\downarrow&&&&&&&&\\
&&&&\cdots&&\cdots&&&&&&&&\\
&&&&\downarrow&&\downarrow&&&&&&&&\\
&&&&P_{n-1,n-2}&\RA{}&P_{n-2,n-2}\\
&&&&\downarrow\\
&&&&P_{n-1,n-1}
\end{diag}\]}

\vskip-1em
Taking a mapping cone, we obtain $\Ker f\in\xx_n=\add M$. Thus we obtain $\Omega^nX=\hom_\Lambda(M,\Ker f)\in\pr\Gamma$ and $\pd_\Gamma X\le n$.\rule{5pt}{10pt}

\vskip.5em{\bf\XDGB\ Theorem }{\it
(1) Let $\Lambda$ be an artin algebra. Then $r_\Lambda(\Lambda)\le2$ implies $\ind\xx_2<\infty$, and the converse holds if $\xx_2$ consists of reflexive $\Lambda$-modules. In particular, if $\Lambda$ is $1$-Gorenstein, then $r_\Lambda(\Lambda)\le2$ if and only if $\Lambda$ is reflexive-finite (\XDG).

(2) Let $\Lambda$ be an order in a semisimple algebra. Then $r_\Lambda(\Lambda)\le2$ if and only if $\expdim\Lambda\le 2$ if and only if $\Lambda$ is reflexive-finite.}

\vskip.5em{\sc Proof }
(1) By \XDGA, we only have to show the former assertion. Assume that $\Gamma:=\endm_\Lambda(\Lambda\oplus X)$ satisfies $\gl\Gamma\le 2$. Then $P:=\hom_\Lambda(\Lambda\oplus X,\Lambda)$ satisfies $\endm_\Gamma(P)=\Lambda$, and $\ppp:=\hom_\Gamma(P,\ ):\pr\Gamma\to\mod\Lambda$ is full faithful. We only have to show $\ppp(\pr\Gamma)\supseteq\Omega^2(\mod\Lambda)$. For any projective resolution $0\to \Omega^2X\to P_1\stackrel{f}{\to}P_0\to X\to 0$ in $\mod\Lambda$, there exists a morphism $g:Q_1\to Q_0$ in $\add_\Gamma P$ such that $f=\ppp g$. Since $\gl\Gamma\le 2$, $\Ker g\in\pr\Gamma$ holds. Thus we have $\Omega^2X=\ppp(\Ker g)\in\ppp(\pr\Gamma)$.

(2) Since $\Lambda$ is an order in semisimple algebra, $\Lambda$ is $1$-Gorenstein (e.g. \XDI) and $\xx_2$ is the category of reflexive $\Lambda$-modules. By \XDGA\ and \XDEA(3), we only have to show that $\expdim\Lambda\le 2$ implies that $\xx_2$ has an additive generator. Assume that an order $\Gamma$ and $P\in\pr\Gamma$ satisfy $\gl\Gamma\le2$ and $\endm_\Gamma(P)=\Lambda$. Assume that an idempotent $e$ of $\Gamma$ satisfies $P\in\add_\Gamma(\Gamma e)$ and $J_\Gamma(1-e)\in\pr\Gamma$. Then $\Gamma^\prime:=e\Gamma e$ and $P^\prime:= eP\in\pr\Gamma^\prime$ satisfy $\gl\Gamma^\prime\le 2$ and $\endm_{\Gamma^\prime}(P^\prime)=\Lambda$ by \XCB(3), so we can replace $(\Gamma,P)$ by $(\Gamma^\prime,P^\prime)$. Repeating this procedure, we can assume that any simple $\Gamma$-modules $S$ with $\pd_\Gamma S\le 1$ is a factor module of $P$. Then it is easily checked that any $X\in\mod\Gamma$ with $\hom_\Gamma(P,X)=0$ satisfies $\ext^i_\Gamma(X,\Gamma)=0$ for $i=0,1$ by [I2;III.2.2.1]. Thus $\ppp:=\hom_\Gamma(P,\ ):\pr\Gamma\to\mod\Lambda$ is full faithful by [I6;2.2.1]. Now $\ppp(\pr\Gamma)\supseteq\Omega^2(\mod\Lambda)$ follows from the argument in the proof of (1).\rule{5pt}{10pt}

\vskip.5em{\bf\XDH\ }Let $\fdim\Lambda:=\sup\{\pd X\ |\ X\in\mod\Lambda,\ \pd X<\infty\}$ be the {\it finitistic dimension} of $\Lambda$ [B]. The finitistic dimension conjecture (FDC) asserts that $\fdim\Lambda<\infty$ holds for any artin algebra $\Lambda$. We refer to [Z] for known results and the relationship to other homological conjecture. Recently, Igusa-Todorov [IT] introduced a function $\psi_\Lambda$ and applied it to prove (FDC) for artin algebras with $\rdim\Lambda\le 3$. Their result is valid for much larger class of rings including orders. We shall give its proof for the completeness. We refer [EHIS] and [Xi1,2] for approach to (FDC) using Igusa-Todorov's theorem.

\vskip.5em{\bf Theorem }{\it
Let $\Lambda$ be an artin algebra or order, and $\xx_n:=\add\Omega^n(\mod\Lambda)$. If $\wresdim\xx_n\le 1$ holds for some $n$, then $\fdim\Lambda<\infty$. Thus $\rdim\Lambda\le3$ (resp. $r_\Lambda(\Lambda)\le3$, $\expdim\Lambda\le3$) implies $\fdim\Lambda<\infty$.}

\vskip.5em{\bf\XDHA\ Lemma }{\it
Let $\Lambda$ be a noetherian ring such that $\mod\Lambda$ forms a Krull-Schmidt category and $\ext^i_\Lambda(X,Y)\in\mod\endm_\Lambda(X)$ holds for any $X,Y\in\mod\Lambda$ and $i\ge0$. Then there exists a function $\psi_\Lambda:\mod\Lambda\to\nnn_{\ge0}$ with the following properties.

(i) If $\pd X<\infty$, then $\psi_\Lambda(X)=\pd X$.

(ii) $\add X\subseteq\add Y$ implies $\psi_\Lambda(X)\le\psi_\Lambda(Y)$.

(iii) If $0\to X\to Y\to Z\to0$ is exact with $\pd Z<\infty$, then $\pd Z\le\psi_\Lambda(X\oplus Y)+1$.}

\vskip.5em{\sc Proof }
For $X\in\mod\Lambda$, put $a_X(i):=i+\sup\{\pd Y\ |\ Y\in\add\Omega^iX,\ \pd Y<\infty\}\in\nnn_{\ge0}$. Then $a_X(0)\le a_X(1)\le a_X(2)\le\cdots$ holds. We denote by $G_X$ the free abelian group with the basis $\ind(\add X)\backslash\ind(\pr\Lambda)$, which is a finite set since $\mod\Lambda$ is Krull-Schmidt. Then $\Omega$ gives an element of $\endm_{\zzz}(G_X)$. Since $G_X$ is a finitely generated $\zzz$-module, there exists $m\in\nnn_{\ge0}$ such that $\Omega:\Omega^nG_X\to\Omega^{n+1}G_X$ is an isomorphism for any $n$ ($n\ge m$). We denote $\phi_\Lambda(X)$ the minimal value of $m$ with this property. We shall show that $\psi_\Lambda(X):=a_X(\phi_\Lambda(X))$ satisfies the desired properties.

If $\pd X<\infty$, then $\phi_\Lambda(X)=\pd X$ and $\phi_\Lambda(X)=a_X(\pd X)=\pd X$ holds. If $\add X\subseteq\add Y$, then $G_X\subseteq G_Y$ implies $\phi_\Lambda(X)\le\phi_\Lambda(Y)$ and $\psi_\Lambda(X)=a_X(\phi_\Lambda(X))\le a_Y(\phi_\Lambda(X))\le a_Y(\phi_\Lambda(Y))=\psi_\Lambda(Y)$. We shall show (iii). Considering the $(\pd Z)$-th syzygies of the given exact sequences, we obtain $\Omega^{\pd Z}(X-Y)=0$ in $G_{X\oplus Y}$. By the definition of $e:=\phi_\Lambda(X\oplus Y)$, we have $\Omega^e(X-Y)=0$ in $G_{X\oplus Y}$. Thus the exact sequence of $e$-th syzygies has the form $0\to W\oplus P\stackrel{f}{\to} W\oplus P^\prime\to\Omega^eX\to0$ for some $W\in\mod\Lambda$ and $P,P^\prime\in\pr\Lambda$. Then $\pd\Omega^eX\le \pd W+1$ holds. Cancelling the trivial direct summand of complexes, we can assume $f\in J_{\mod\Lambda}$. Now we have a long exact sequence $\cdots\to\ext^{i}_\Lambda(\Omega^eX,\ )\to\ext^i_\Lambda(W,\ )\stackrel{f}{\to}\ext^i_\Lambda(W,\ )\to\cdots$. Using Nakayama's lemma on $\endm_\Lambda(W)$-modules, we obtain $\pd W\le \pd\Omega^eX$. Since $W$ is a direct summand of $\Omega^eY$, we have $\pd X-1\le e+\pd\Omega^eX-1\le e+\pd W\le a_{X\oplus Y}(e)=\psi_\Lambda(X\oplus Y)$ holds.\rule{5pt}{10pt}

\vskip.5em{\bf\XDHB\ Proof of \XDH\ }
By \XDEA(3)(4), we only have to show the former assertion. Let $M\in\dn{M}_\Lambda$ be in \XDE(2)(i). For any $X\in\mod\Lambda$ with $\pd X<\infty$, take an exact sequence $0\to M_1\to M_0\to Y\to0$ with $M_i\in\add M$ and $\Omega^nX\in\add Y$. Then $\pd X\le\pd Y+n\le \psi_\Lambda(M_1\oplus M_0)+n+1=\psi_\Lambda(M)+n+1$ holds by \XDHA(ii)(iii).\rule{5pt}{10pt}

\vskip.5em{\bf\XDI\ }
Let $\Gamma$ be an order in a semisimple algebra. Any $X\in\lat\Lambda$ has a minimal relative-injective resolution $0\to\Gamma\to I_0(X)\to I_1(X)\to\cdots$ with $I_i(X)\in\rin\Gamma$, which is a dual of a minimal projective resolution. In this case, Fujita [F] proved that a minimal injective resolution of $X$ is given by $0\to X\to X\otimes_RK\to I_0(X)\otimes_R(K/R)\to I_1(X)\otimes_R(K/R)\to\cdots$. Putting $X:=\Gamma$, we obtain that $\Gamma$ is always $1$-Gorenstein, and $\Gamma$ is $2$-Gorenstein if and only if $I_0(\Gamma)\in\pr\Gamma$ [FGR].

We denote by ${\bf A}(\Lambda)$ the set of Morita-equivalence classes of orders $\Gamma$ in semisimple algebras such that $\Gamma$ is $2$-Gorenstein and $\endm_\Gamma(I_0(\Gamma))^{op}$ is Morita-equivalent to $\Lambda$. The theorem below shows that our definition of $\rdim\Lambda$ given in \XDA\ is surely one-dimensional analogy of the representation dimension of artin algebras defined in [A1].

\vskip.5em{\bf Theorem }{\it
Let $\Lambda$ be an order in a semisimple algebra. Then ${\bf A}(\Lambda)=\{\endm_\Lambda(M)\ |\ M\in\lat\Lambda,\ \Lambda\oplus D\Lambda\in\add M\}$ holds. Thus $\rdim\Lambda=\inf\{\gl\Gamma\ |\ \Gamma\in{\bf A}(\Lambda)\}$ holds.}

\vskip.5em{\sc Proof }
For $M\in\lat\Lambda$ with $\Lambda\oplus D\Lambda\in\add M$, we will show $\Gamma:=\endm_\Lambda(M)\in{\bf A}(\Lambda)$. Put $I_i:=I_i(M)$. Then $\add I_0=\rin\Lambda$ holds by $D\Lambda\in\add M$. We have an exact sequence $0\to\Gamma\to\hom_\Lambda(M,I_0)\stackrel{f}{\to}\hom_\Lambda(M,I_1)$ of $\Gamma$-modules such that $f$ is in $J_{\lat\Gamma}$. Since $\hom_\Lambda(M,D\Lambda)=DM=D\hom_\Lambda(\Lambda,M)\in\rin\Gamma$ holds by $\Lambda\in\add M$, we obtain $I_0(\Gamma)=\hom_\Lambda(M,I_0)$. Since $\hom_\Lambda(M,I_0)\in\pr\Gamma$ holds by $D\Lambda\in\add M$, $\Gamma$ is $2$-Gorenstein. Moreover, $\endm_\Gamma(I_0(\Gamma))=\endm_\Lambda(I_0)$ is Morita-equivalent to $\endm_\Lambda(D\Lambda)=\Lambda^{op}$.

Conversely, take $\Gamma\in{\bf A}(\Lambda)$. We can assume $\Lambda=\endm_\Gamma(I_0(\Gamma))^{op}=\endm_{\Gamma^{op}}(DI_0(\Gamma))$. Then $M:=DI_0(\Gamma)\in\lat\Lambda$ is a projective-injective $\Gamma^{op}$-module. Thus $D\Lambda=D\endm_{\Gamma^{op}}(M)\in\add_{\Lambda}D\hom_{\Gamma^{op}}(M,D\Gamma)=\add_{\Lambda}M$ and $\Lambda=\endm_{\Gamma^{op}}(M)\in\add_\Lambda\hom_{\Gamma^{op}}(\Gamma,M)=\add_{\Lambda}M$ holds. We obtain $\Lambda\oplus D\Lambda\in\add_\Lambda M$. Put $P:=\hom_{\Gamma^{op}}(M,\Gamma)\in\pr\Gamma$. It is easily checked that any $X\in\mod\Gamma$ with $\hom_\Gamma(P,X)=0$ satisfies $\ext^i_\Gamma(X,\Gamma)=0$ for $i=0,1$ by [I2;III.2.2.1]. Thus the functor $\hom_\Gamma(P,\ ):\pr\Gamma\to\mod\Lambda$ is full faithful by [I6;2.2.1], and $\Gamma=\endm_\Gamma(\Gamma)=\endm_\Lambda(\hom_\Gamma(P,\Gamma))=\endm_\Lambda(M)$ holds.\rule{5pt}{10pt}

\vskip.5em{\bf\XE\ Auslander-Reiten theory and Finite equivalence}

Throughtout this section, let $\Lambda$ be an artin algebra or order in a semisimple algebra, and $\dn{M}_\Lambda$ the category in \XZA.

\vskip.5em{\bf\XEA\ }
To study the factor categories of $\dn{M}_\Lambda$, it is convenient to introduce the concept of $\tau$-categories [I2]. Recall that a Krull-Schmidt category $\cc$ is called a {\it $\tau$-category} if the following conditions are satisfied.

\strut\kern1em(i) Any $X\in\ind\cc$ has a complex $(X]_{\cc}:=(\tau^+_{\cc}X\to\theta^+_{\cc}X\to X)$ in $\cc$ such that $\cc(\ ,\tau^+_{\cc}X)\to\cc(\ ,\theta^+_{\cc}X)\to J_{\cc}(\ ,X)\to0$ gives a first two terms of a minimal projective resolution in $\mod\cc$. If $\tau^+_{\cc}X\neq0$, then $\cc(X,\ )\to\cc(\theta^+_{\cc}X,\ )\to J_{\cc}(\tau^+_{\cc}X,\ )\to0$ is exact.

\strut\kern1em(ii) Any $X\in\ind\cc$ has a complex $[X)_{\cc}:=(X\to\theta^-_{\cc}X\to\tau^-_{\cc}X)$ in $\cc$ such that $\cc(\tau^-_{\cc}X,\ )\to\cc(\theta^-_{\cc}X,\ )\to J_{\cc}(X,\ )\to0$ gives a first two terms of a minimal projective resolution in $\mod\cc^{op}$. If $\tau^-_{\cc}X\neq0$, then $\cc(\ ,X)\to\cc(\ ,\theta^-_{\cc}X)\to J_{\cc}(\ ,\tau^-_{\cc}X)\to0$ is exact.

Put $\ind^+_1\cc:=\{ X\in\ind\cc\ |\ \tau^+_{\cc}X=0\}$ and $\ind^-_1\cc:=\{ X\in\ind\cc\ |\ \tau^-_{\cc}X=0\}$. Then $\tau^+_{\cc}$ gives a bijection between $\ind\cc\backslash\ind^+_1\cc\to\ind\cc\backslash\ind^-_1\cc$ with the inverse $\tau^-_{\cc}$ [I2;I.2.3]. For any $X\in\ind\cc\backslash\ind^+_1\cc$, the complexes $(X]_{\cc}$ and $[\tau^+_{\cc}X)_{\cc}$ are isomorphic.

\vskip.5em{\bf\XEAA\ Example }[I2;I.2.2](1) Let $\Lambda$ be an artin algebra or order in a semisimple algebra. Then $\dn{M}_\Lambda$ forms a $\tau$-category by Auslander-Reiten theory (\XEBB).

(2) Let $Q$ be a translation quiver. Then the mesh category of $Q$ forms a $\tau$-category.

\vskip.5em{\bf\XEAB\ }
Let $\cc$ and $\cc^\prime$ be $\tau$-categories and $\fff:\cc\to\cc^\prime$ an equivalence of categories. Then $\fff$ induces bijections $\ind^+_1\cc\to\ind^+_1\cc^\prime$ and $\ind^-_1\cc\to\ind^-_1\cc^\prime$. Moreover, $\tau^+_{\cc^\prime}\circ\fff X=\fff\circ\tau^+_{\cc}X$ holds for any $X\in\ind\cc\backslash\ind^+_1\cc$.

\vskip.5em{\sc Proof }
Since $\ind^+_1$ and $\ind^-_1$ are defined categoircally, the former assertion follows. Similarly, the complex $\tau^+_{\cc}X\to\theta^+_{\cc}X\to X$ is preserved by $\fff$.\rule{5pt}{10pt}

\vskip.5em{\bf\XEAC\ Proposition }[I2;II.1.4] {\it
Let $\cc$ be a $\tau$-category and $\cc^\prime$ a subcategory of $\cc$. Then $\overline{\cc}:=\cc/[\cc^\prime]$ forms a $\tau$-category again. We regard $\ind\cc$ as a disjoint union of $\ind\cc^\prime$ and $\ind(\cc/[\cc^\prime])$. We denote by $\overline{(\ )}:\cc\to\overline{\cc}$ the natural functor. Let $X\in\ind\cc\backslash\ind\cc^\prime$. 

(1) $(X]_{\overline{\cc}}=(0\to0\to X)$ if $\theta^+_{\cc}X\in\cc^\prime$, and $(X]_{\overline{\cc}}=\overline{(X]}_{\cc}$ othewise. Dually, $[X)_{\overline{\cc}}=(X\to0\to 0)$ if $\theta^-_{\cc}X\in\cc^\prime$, and $[X)_{\overline{\cc}}=\overline{[X)}_{\cc}$ otherwise.

(2) $X$ is contained in $\ind^-_1\overline{\cc}$ if and only if either $\tau^-_{\cc}X\in\cc^\prime$ or $\theta^-_{\cc}X\in\cc^\prime$ holds. Hence $\ind^-_1\overline{\cc}$ is a disjoint union of three subsets $\ind^-_1\cc\backslash\ind\cc^\prime$, $\tau^+_{\cc}(\ind\cc^\prime\backslash\ind^+_1\cc)\backslash\ind\cc^\prime$ and $s_{\cc^\prime}(\cc):=\{X\in\ind\cc\backslash\ind\cc^\prime\ |\ \theta^-_{\cc}X\in\cc^\prime,\ \tau^-_{\cc}X\notin\cc^\prime\}$.

(3) $X$ is contained in $\ind^+_1\overline{\cc}$ if and only if either $\tau^+_{\cc}X\in\cc^\prime$ or $\theta^+_{\cc}X\in\cc^\prime$ holds. Hence $\ind^+_1\overline{\cc}$ is a disjoint union of three subsets $\ind^+_1\cc\backslash\ind\cc^\prime$, $\tau^-_{\cc}(\ind\cc^\prime\backslash\ind^-_1\cc)\backslash\ind\cc^\prime$ and $s^{\cc^\prime}(\cc):=\{X\in\ind\cc\backslash\ind\cc^\prime\ |\ \theta^+_{\cc}X\in\cc^\prime,\ \tau^+_{\cc}X\notin\cc^\prime\}$.}

\vskip.5em{\bf\XEB\ Auslander-Reiten theory }

We recall well-known results in Auslander-Reiten theory [ARS]. We call $\underline{\dn{M}}_\Lambda:=\dn{M}_\Lambda/[\pr\Lambda]$ the {\it stable category} of $\Lambda$, and $\overline{\dn{M}}_\Lambda:=\dn{M}_\Lambda/[\rin\Lambda]$ the {\it costable category} of $\Lambda$. In both cases of artin algebras and orders, we have an equivalence $\tau^+_\Lambda:\underline{\dn{M}}_\Lambda\to\overline{\dn{M}}_\Lambda$ called the {\it Auslander translate}. It induces a bijection $\ind\dn{M}_\Lambda\backslash\ind(\pr\Lambda)\to\ind\dn{M}_\Lambda\backslash\ind(\rin\Lambda)$. The existence of a functorial isomorphism $D\underline{\dn{M}}_\Lambda(X,Y)\to\ext^1_\Lambda(Y,\tau^+_\Lambda X)$ for any $X,Y\in\dn{M}_\Lambda$ called {\it Auslander-Reiten isomorphism} is one of the most important theorem in the representation theory. This immediately implies the following existence theorem of almost split sequences.

\vskip.5em{\bf\XEBA\ Theorem }{\it
Let $\Lambda$ be an artin algebra or order in a semisimple algebra. For any $X\in\ind\dn{M}_\Lambda\backslash\ind(\pr\Lambda)$, there exists an exact sequence $0\to\tau^+_\Lambda X\to\theta^+_\Lambda X\to X\to0$ such that $0\to\dn{M}_\Lambda(\ ,\tau^+_\Lambda X)\to\dn{M}_\Lambda(\ ,\theta^+_\Lambda X)\to J_{\dn{M}_\Lambda}(\ ,X)\to0$ and $0\to\dn{M}_\Lambda(X,\ )\to\dn{M}_\Lambda(\theta^+_\Lambda X,\ )\to J_{\dn{M}_\Lambda}(\tau^+_\Lambda X,\ )\to0$ are exact.}

\vskip.5em{\bf\XEBB\ }
Consequently, $\dn{M}_\Lambda$ forms a $\tau$-category with $\ind^+_1\dn{M}_\Lambda=\ind(\pr\Lambda)$ and $\ind^-_1\dn{M}_\Lambda=\ind(\rin\Lambda)$. In fact, for any $P\in\ind(\pr\Lambda)$ and $I\in\ind(\rin\Lambda)$, the sequence $0\to J_\Lambda X\to X$ and $X\to D((DX)J_\Lambda)\to 0$ induces isomorphisms $\dn{M}_\Lambda(\ ,J_\Lambda X)\to J_{\dn{M}_\Lambda}(\ ,X)$ and $\dn{M}_\Lambda(D((DX)J_\Lambda),\ )\to J_{\dn{M}_\Lambda}(X,\ )$ on $\dn{M}_\Lambda$.

\vskip.5em{\bf\XEC\ Relative homology }

Let us recall relative homology theory of Auslander-Solberg [ASo]. There exists a bijection between functorially finite subcategories $\xx$ of $\dn{M}_\Lambda$ such that $\Lambda\in\xx$ and functorially finite subcategories $\yy$ of $\dn{M}_\Lambda$ such that $D\Lambda\in\yy$, where the correspondense is given by $\ind\yy\backslash\ind(\rin\Lambda)=\tau^+_\Lambda(\ind\xx\backslash\ind(\pr\Lambda))$. In the rest, fix such a pair $(\xx,\yy)$ and put $\underline{\underline{\dn{M}}}_\Lambda:=\dn{M}_\Lambda/[\xx]$ and $\overline{\overline{\dn{M}}}_\Lambda:=\dn{M}_\Lambda/[\yy]$. Then the Auslander translate $\tau^+_\Lambda$ still induces equivalences $\underline{\underline{\dn{M}}}_\Lambda\to\overline{\overline{\dn{M}}}_\Lambda$ and $\xx/[\pr\Lambda]\to\yy/[\rin\Lambda]$.

We denote by $\underline{\underline{\dn{M}}}_\Lambda(X,Y)$ the set of morphisms in $\underline{\underline{\dn{M}}}_\Lambda$. For $X,Z\in\dn{M}_\Lambda$, we denote by $F_{\xx}(X,Z)$ (resp. $F^{\yy}(X,Z)$) the subset of $\ext^1_\Lambda(X,Z)$ consisting of $0\to Z\to Y\to X\to0$ such that $\dn{M}_\Lambda(\ ,Y)\to\dn{M}_\Lambda(\ ,X)\to0$ is exact on $\xx$ (resp. $\dn{M}_\Lambda(Y,\ )\to\dn{M}_\Lambda(Z,\ )\to0$ is exact on $\yy$). Then the Auslander-Reiten isomorphism implies that $F_{\xx}=F^{\yy}$ holds, and $F_{\xx}$ gives an additive subfunctor of $\ext^1_\Lambda(\ ,\ )$. Put $\ext^i_{F_{\xx}}(X,Z):=F_{\xx}(\Omega_{\xx}^{i-1}X,Z)$ , where $\Omega_{\xx}$ is defined in \XAA. A main theorem in [ASo] is the relative version below of Auslander-Reiten isomorphism. This immediately implies that indecomposable projective objects in $\mod\underline{\underline{\dn{M}}}_\Lambda$ are $\underline{\underline{\dn{M}}}_\Lambda(\ ,X)$ for $X\in\ind\dn{M}_\Lambda\backslash\ind\xx$, and indecomposable injective objects in $\mod\underline{\underline{\dn{M}}}_\Lambda$ are $F_{\xx}(\ ,X)$ for $X\in\ind\dn{M}_\Lambda\backslash\ind\yy$.

\vskip.5em{\bf\XECA\ Theorem }{\it
There exists a functorial isomorphism $D\underline{\underline{\dn{M}}}_\Lambda(X,Y)\to F_{\xx}(Y,\tau^+_\Lambda X)$ for any $X,Y\in\dn{M}_\Lambda$.}

\vskip.5em{\bf\XECB\ Theorem }{\it
Let $M\in\mod\underline{\underline{\dn{M}}}_\Lambda$. We denote by $0\to\dn{M}_\Lambda(\ ,Z)\stackrel{\cdot g}{\to}\dn{M}_\Lambda(\ ,Y)\stackrel{\cdot f}{\to}\dn{M}_\Lambda(\ ,X)\to M\to0$ the minimal projective resolution of $M$ in $\mod\dn{M}_\Lambda$.

(1) A projective resolution of $M$ in $\mod\underline{\underline{\dn{M}}}_\Lambda$ is given by the sequence below, where the first two terms $\underline{\underline{\dn{M}}}_\Lambda(\ ,Y)\stackrel{\cdot f}{\to}\underline{\underline{\dn{M}}}_\Lambda(\ ,X)\to M\to0$ is minimal.
\begin{eqnarray*}
\cdots\to\underline{\underline{\dn{M}}}_\Lambda(\ ,\Omega_{\xx}^2X)\to\underline{\underline{\dn{M}}}_\Lambda(\ ,\Omega_{\xx} Z)\to\underline{\underline{\dn{M}}}_\Lambda(\ ,\Omega_{\xx} Y)\to\underline{\underline{\dn{M}}}_\Lambda(\ ,\Omega_{\xx} X)\to\ \ \ \ \ \ \ \ \ \ &&\\
\underline{\underline{\dn{M}}}_\Lambda(\ ,Z)\stackrel{\cdot g}{\to}\underline{\underline{\dn{M}}}_\Lambda(\ ,Y)\stackrel{\cdot f}{\to}\underline{\underline{\dn{M}}}_\Lambda(\ ,X)\to M\to0&&\end{eqnarray*}

(2) An injective resolution of $M$ in $\mod\underline{\underline{\dn{M}}}_\Lambda$ is given by sequence below, where the first two terms $0\to M\to\ext^1_{F_{\xx}}(\ ,Z)\to\ext^1_{F_{\xx}}(\ ,Y)$ is minimal.
\begin{eqnarray*}
&&0\to M\to\ext^1_{F_{\xx}}(\ ,Z)\to\ext^1_{F_{\xx}}(\ ,Y)\to\ext^1_{F_{\xx}}(\ ,X)\\
&&\ \ \ \ \ \ \ \ \ \ \to\ext^2_{F_{\xx}}(\ ,Z)\to\ext^2_{F_{\xx}}(\ ,Y)\to\ext^2_{F_{\xx}}(\ ,X)\to\ext^3_{F_{\xx}}(\ ,Z)\to\cdots\end{eqnarray*}}

\vskip-.5em{\sc Proof }
(1) was proved in \XAAC\ except the minimality, which follows from the fact that $f$ and $g$ are in $J_{\dn{M}_\Lambda}$. (2) is the Matlis dual of (1) by \XECA.\rule{5pt}{10pt}

\vskip.5em{\bf\XED\ }
We again fix a pair $(\xx,\yy)$ of subcategories of $\dn{M}_\Lambda$ as in \XEC.

(1) $s_{\xx}(\dn{M}_\Lambda)=\{X\in\ind\dn{M}_\Lambda\backslash(\ind\xx\cup\ind\yy)\ |\ \theta_\Lambda^-X\in\xx\}=s^{\yy}(\dn{M}_\Lambda)=\{X\in\ind\dn{M}_\Lambda\backslash(\ind\xx\cup\ind\yy)\ |\ \theta_\Lambda^+X\in\yy\}$ holds, where we use the notations in \XEAC.

(2) $\ind^-_1\underline{\underline{\dn{M}}}_\Lambda=(\ind\yy\backslash\ind\xx)\cup s_{\xx}(\dn{M}_\Lambda)$ and $\ind^+_1\overline{\overline{\dn{M}}}_\Lambda=(\ind\xx\backslash\ind\yy)\cup s_{\xx}(\dn{M}_\Lambda)$ hold by \XEAC. In particular, putting $\widehat{s}_{\xx}(\dn{M}_\Lambda):=\ind\xx\cup\ind\yy\cup s_{\xx}(\dn{M}_\Lambda)$, we obtain $\underline{\underline{\dn{M}}}_\Lambda/[\ind^-_1\underline{\underline{\dn{M}}}_\Lambda]=\dn{M}_\Lambda/[\widehat{s}_{\xx}(\dn{M}_\Lambda)]=\overline{\overline{\dn{M}}}_\Lambda/[\ind^+_1\overline{\overline{\dn{M}}}_\Lambda]$.

\vskip.5em{\bf\XEDA\ Proposition }{\it
Let $\cc$ and $\cc^\prime$ be finite subcategories of $\dn{M}_\Lambda$ such that $\xx\subseteq\cc^\prime\subseteq\cc\subseteq\cc^\prime+\add s_{\xx}(\dn{M}_\Lambda)$. Then $\resdim{\cc}{\dn{M}_\Lambda}\le\resdim{\cc^\prime}{\dn{M}_\Lambda}$.}

\vskip.5em{\sc Proof }We shall use the induction on $\resdim{\cc}{X}$. We can assume $X\in\ind\dn{M}_\Lambda\backslash\ind\cc$. Since $\Lambda\in\xx$, a right $\cc^\prime$-approximation of $X$ gives an exact sequence $0\to\Omega_{\cc^\prime}X\to Y\stackrel{f}{\to}X\to0$. Since any $a\in\dn{M}_\Lambda(Z,X)$ with $Z\in\ind\cc\backslash\ind\cc^\prime$ factors throught $g:Z\to\theta_\Lambda^-Z$ by $X\notin\cc$, there exists $a^\prime$ such that $a=ga^\prime$. Since $\theta_\Lambda^-Z\in\xx$ holds by $Z\in s_{\xx}(\dn{M}_\Lambda)$, there exists $a^{\prime\prime}$ such that $a^\prime=a^{\prime\prime}f$. Thus $a=(ga^{\prime\prime})f$ holds. This implies that $f$ is a right $\cc$-approximation of $X$. Thus the inductive assumption shows $\resdim{\cc}{X}=\resdim{\cc}{\Omega_{\cc^\prime}X}+1\le\resdim{\cc^\prime}{\Omega_{\cc^\prime}X}+1=\resdim{\cc^\prime}{X}$.\rule{5pt}{10pt}

\vskip.5em{\bf\XEE\ }Let us state the main theorem in this section. Let $\Lambda$ and $\Gamma$ be artin algebras or orders in semisimple algebras, $(\xx,\yy)$ and $(\xx^\prime,\yy^\prime)$ pairs of finite subcategories of $\dn{M}_\Lambda$ and $\dn{M}_\Gamma$ respectively satisfying the same condition in \XEC. Put $\widehat{s}_{\xx}(\dn{M}_\Lambda):=\ind\xx\cup\ind\yy\cup s_{\xx}(\dn{M}_\Lambda)$ and $\widehat{s}_{\xx^\prime}(\dn{M}_\Gamma):=\ind\xx^\prime\cup\ind\yy^\prime\cup s_{\xx^\prime}(\dn{M}_\Gamma)$, which are finite sets.

\vskip.5em{\bf Theorem }{\it
Let $\fff:\underline{\underline{\dn{M}}}_\Lambda:=\dn{M}_\Lambda/[\xx]\to\underline{\underline{\dn{M}}}_\Gamma:=\dn{M}_\Gamma/[\xx^\prime]$ be an equivalence of categories.

(1) $\fff$ induces a bijection $\ind\dn{M}_\Lambda\backslash \widehat{s}_{\xx}(\dn{M}_\Lambda)\to\ind\dn{M}_\Gamma\backslash\widehat{s}_{\xx^\prime}(\dn{M}_\Gamma)$. 

(2) Let $\cc$ and $\cc^\prime$ be finite subcategories of $\dn{M}_\Lambda$ and $\dn{M}_\Gamma$ respectively such that $\xx\subset\cc$, $\xx^\prime\subset\cc^\prime$ and $\fff(\ind\cc\backslash\ind\xx)=\ind\cc^\prime\backslash\ind\xx^\prime$. Then $\widehat{s}_{\xx}(\dn{M}_\Lambda)\subset\cc$ if and only if $\widehat{s}_{\xx^\prime}(\dn{M}_\Gamma)\subset\cc^\prime$. In this case, $g_\Lambda(\cc)=g_\Gamma(\cc^\prime)$ and $r_\Lambda(\cc)=r_\Gamma(\cc^\prime)$ hold.}

\vskip.5em{\bf\XEEA\ }$\ggg:=\tau^+_{\Gamma}\circ\fff\circ\tau^-_\Lambda$ gives an equivalence $\overline{\overline{\dn{M}}}_\Lambda:=\dn{M}_\Lambda/[\yy]\to\overline{\overline{\dn{M}}}_\Gamma:=\dn{M}_\Gamma/[\yy^\prime]$. Thus we have bijections $\fff:\ind\dn{M}_\Lambda\backslash\ind\xx\to\ind\dn{M}_\Gamma\backslash\ind\xx^\prime$ and $\ggg:\ind\dn{M}_\Lambda\backslash\ind\yy\to\ind\dn{M}_\Gamma\backslash\ind\yy^\prime$. Define a subfunctor $F_{\xx^\prime}$ of $\ext^1_\Gamma(\ ,\ )$ by \XEC.

(1) The induced equivalence $\fff:\mod\underline{\underline{\dn{M}}}_\Lambda\to\mod\underline{\underline{\dn{M}}}_\Gamma$ satisfies $\fff\underline{\underline{\dn{M}}}_\Lambda(\ ,X)=\underline{\underline{\dn{M}}}_\Gamma(\ ,\fff X)$ and $\fff F_{\xx}(\ ,X)=F_{\xx^\prime}(\ ,\ggg X)$.

(2) The following diagram is commutative, where the horizontal equivalences are induced by \XEAB\ and the vertical identities are given by \XED(2).
\[\begin{diag}
\underline{\underline{\dn{M}}}_\Lambda/[\ind^-_1\underline{\underline{\dn{M}}}_\Lambda]&\RA{\fff}&\underline{\underline{\dn{M}}}_\Gamma/[\ind^-_1\underline{\underline{\dn{M}}}_\Gamma]\\
\parallel&&\parallel\\
\dn{M}_\Lambda/[\widehat{s}_{\xx}(\dn{M}_\Lambda)]&&\dn{M}_\Gamma/[\widehat{s}_{\xx^\prime}(\dn{M}_\Gamma)]\\
\parallel&&\parallel\\
\overline{\overline{\dn{M}}}_\Lambda/[\ind^+_1\overline{\overline{\dn{M}}}_\Lambda]&\RA{\ggg}&\overline{\overline{\dn{M}}}_\Gamma/[\ind^+_1\overline{\overline{\dn{M}}}_\Gamma]
\end{diag}\]

\vskip.5em{\sc Proof }
(1) $\fff F_{\xx}(\ ,X)=\fff(D\underline{\underline{\dn{M}}}_\Lambda(\tau_\Lambda^-X,\ ))=D\underline{\underline{\dn{M}}}_\Gamma(\fff\circ\tau_\Lambda^-X,\ )=F_{\xx^\prime}(\ ,\ggg X)$ by \XECA.

(2) Since $\ggg X=\tau^+_\Gamma\circ\fff\circ\tau^-_\Lambda X=\fff X$ holds for any $X\in\ind\underline{\underline{\dn{M}}}_\Lambda\backslash\ind^-_1\underline{\underline{\dn{M}}}_\Lambda$ by \XEAB\ and \XEAC(1), the commutativity follows.\rule{5pt}{10pt}

\vskip.5em{\bf\XEEB\ }We denote by $0\to\dn{M}_\Lambda(\ ,Z)\to\dn{M}_\Lambda(\ ,Y)\to\dn{M}_\Lambda(\ ,X)\to M\to0$ the minimal projective resolution of $M\in\mod\underline{\underline{\dn{M}}}_\Lambda$ in $\mod\dn{M}_\Lambda$, and by $0\to\dn{M}_\Gamma(\ ,Z^\prime)\to\dn{M}_\Gamma(\ ,Y^\prime)\to\dn{M}_\Gamma(\ ,X^\prime)\to\fff M\to0$ the minimal projective resolution of $\fff M\in\mod\underline{\underline{\dn{M}}}_\Gamma$ in $\mod\dn{M}_\Gamma$.

(1) $X$ and $X^\prime$ has no direct summand in $\ind\xx$ and $\ind\xx^\prime$ respectively, and $\fff X=X^\prime$.

(2) $\ggg\overline{Y}=\overline{Y}^\prime$ holds in $\overline{\overline{\dn{M}}}_\Gamma$, and $\fff\underline{Y}=\underline{Y}^\prime$ holds in $\underline{\underline{\dn{M}}}_\Gamma$.

(3) $Z$ and $Z^\prime$ has no direct summand in $\ind\yy$ and $\ind\yy^\prime$ respectively, and $\ggg Z=Z^\prime$.

\vskip.5em{\sc Proof }
By \XECB, a minimal projective resolution and a minimal injective resolution of $M$ in $\mod\underline{\underline{\dn{M}}}_\Lambda$ are given by $\underline{\underline{\dn{M}}}_\Lambda(\ ,Y)\to\underline{\underline{\dn{M}}}_\Lambda(\ ,Z)\to M\to0$ and $0\to M\to F_{\xx}(\ ,X)\to F_{\xx}(\ ,Y)$. Similarly, those of $\fff M$ in $\mod\underline{\underline{\dn{M}}}_\Gamma$ is given by $\underline{\underline{\dn{M}}}_\Gamma(\ ,Y^\prime)\to\underline{\underline{\dn{M}}}_\Gamma(\ ,Z^\prime)\to \fff M\to0$ and $0\to \fff M\to F_{\xx^\prime}(\ ,X^\prime)\to F_{\xx^\prime}(\ ,Y^\prime)$. Thus the assertion follows from \XEEA(1).\rule{5pt}{10pt}

\vskip.5em{\bf\XEEC\ Proof of \XEE\ }
(1) and the former assertion of (2) follow from \XEEA(2). We will show the latter assertion of (2). By \XDAB, we only have to show $\resdim{\cc}{X}\ge\resdim{\cc^\prime}{\fff X}$ inductively. We can assume $X\in\ind\dn{M}_\Lambda\backslash\ind\cc$. Take an exact sequence $0\to Z\to Y\stackrel{f}{\to}X\to0$ such that $f$ is a minimal right $\cc$-approximation of $X$. Then $0\to\dn{M}_\Lambda(\ ,Z)\to\dn{M}_\Lambda(\ ,Y)\to\dn{M}_\Lambda(\ ,X)\to M\to0$ gives a minimal projective resolution of $M\in\mod\underline{\underline{\dn{M}}}_\Lambda$ in $\mod\dn{M}_\Lambda$. Now take a minimal projective resolution $0\to\dn{M}_\Gamma(\ ,Z^\prime)\to\dn{M}_\Gamma(\ ,Y^\prime)\stackrel{\cdot f^\prime}{\to}\dn{M}_\Gamma(\ ,X^\prime)\to\fff M\to0$ of $\fff M\in\mod\underline{\underline{\dn{M}}}_\Gamma$ in $\mod\dn{M}_\Gamma$, which is induced by an exact sequence $0\to Z^\prime\to Y^\prime\stackrel{f^\prime}{\to}X^\prime\to0$. Then $X^\prime=\fff X$, $\underline{Y^\prime}=\fff\underline{Y}$ and $Z^\prime=\ggg Z$ hold by \XEEB. Thus $Y^\prime\in\cc^\prime$ holds. Since $\fff M=0$ holds on $\ind\cc^\prime=\fff(\ind\cc\backslash\ind\xx)\cup\ind\xx^\prime$, $f^\prime$ is a right $\cc^\prime$-approximation of $X^\prime=\fff X$. By \XEEA(2), $Z^\prime=\ggg Z$ coincides with $\fff Z$ up to a direct sum of modules in $\widehat{s}_{\xx^\prime}(\dn{M}_\Gamma)\subset\cc^\prime$. Thus the inductive assumption shows $\resdim{\cc}{X}=\resdim{\cc}{Z}+1\ge\resdim{\cc^\prime}{\fff Z}+1=\resdim{\cc^\prime}{Z^\prime}+1\ge\resdim{\cc^\prime}{\fff X}$.\rule{5pt}{10pt}

\vskip.5em{\bf\XEF\ Proof of \XDB\ }
(1) Obviously we can assume $\Lambda\in\xx$ and $\Gamma\in\xx^\prime$. Since $\widehat{s}_{\xx}(\dn{M}_\Lambda)$ and $\widehat{s}_{\xx^\prime}(\dn{M}_\Gamma)$ are finite sets, $|r_\Lambda|=|r_\Gamma|$ follows immediately from \XEE(2).

(2) We shall apply \XEE\ to the case $\xx=\pr\Lambda$, $\yy=\rin\Lambda$, $\xx^\prime=\pr\Gamma$ and $\yy^\prime=\rin\Gamma$. Thus $\widehat{s}_{\xx}(\dn{M}_\Lambda)=\ind(\pr\Lambda)\cup\ind(\rin\Lambda)\cup s_{\xx}(\dn{M}_\Lambda)$ and $\widehat{s}_{\xx^\prime}(\dn{M}_\Gamma)=\ind(\pr\Gamma)\cup\ind(\rin\Gamma)\cup s_{\xx^\prime}(\dn{M}_\Gamma)$. Take a finite subcategory $\cc$ of $\dn{M}_\Lambda$ such that $\Lambda\oplus D\Lambda\in\cc$ and $\rdim\Lambda=g_\Lambda(\cc)$. By \XEDA, we can assume $\widehat{s}_{\xx}(\dn{M}_\Lambda)\subset\cc$. Define the corresponding finite subcategory $\cc^\prime$ of $\dn{M}_\Gamma$ by \XEE. Since $\Gamma\oplus D\Gamma\in\add\widehat{s}_{\xx^\prime}(\dn{M}_\Gamma)\subseteq\cc^\prime$ holds, we obtain $\rdim\Gamma\le g_\Gamma(\cc^\prime)=g_\Lambda(\cc)=\rdim\Lambda$. Exchanging $\Lambda$ and $\Gamma$, we obtain $\rdim\Lambda=\rdim\Gamma$.\rule{5pt}{10pt}

\vskip.5em{\footnotesize
\begin{center}
{\bf References}
\end{center}

[ADW] I. Agoston, V. Dlab, T. Wakamatsu: Neat algebras. Comm. Algebra 19 (1991), no. 2, 433--442.

[As] I. Assem: Tilting theory---an introduction. Topics in algebra, Part 1 (Warsaw, 1988), 127--180, Banach Center Publ., 26, Part 1, PWN, Warsaw, 1990.

[A1] M. Auslander: Representation dimension of Artin algebras. Lecture notes, Queen Mary College, London, 1971.

[A2] M. Auslander: Functors and morphisms determined by objects. Representation theory of algebras (Proc. Conf., Temple Univ., Philadelphia, Pa., 1976), pp. 1--244. Lecture Notes in Pure Appl. Math., Vol. 37, Dekker, New York, 1978. 

[AB] M. Auslander, M. Bridger: Stable module theory. Memoirs of the American Mathematical Society, No. 94 American Mathematical Society, Providence, R.I. 1969 146 pp.

[ABu] M. Auslander, R. O. Buchweitz: The homological theory of maximal Cohen-Macaulay approximations, Soc. Math. de France, Mem. no. 38 (1989) 5--37.

[APT] M. Auslander, M. I. Platzeck, G. Todorov: Homological theory of idempotent ideals. Trans. Amer. Math. Soc. 332 (1992), no. 2, 667--692.

[AR1] M. Auslander, I. Reiten: Applications of contravariantly finite subcategories. Adv. Math. 86 (1991), no. 1, 111--152.

[AR2] M. Auslander, I. Reiten: Homologically finite subcategories. Representations of algebras and related topics (Kyoto, 1990), 1--42, London Math. Soc. Lecture Note Ser., 168, Cambridge Univ. Press, Cambridge, 1992.

[AR3] M. Auslander, I. Reiten: $k$-Gorenstein algebras and syzygy modules. J. Pure Appl. Algebra 92 (1994), no. 1, 1--27.

[AR4] M. Auslander, I. Reiten: Syzygy modules for Noetherian rings. J. Algebra 183 (1996), no. 1, 167--185.

[ARS] M. Auslander, I. Reiten, S. O. Smalo: Representation theory of Artin algebras. Cambridge Studies in Advanced Mathematics, 36. Cambridge University Press, Cambridge, 1995. 

[AS1] M. Auslander, S. O. Smalo: Almost split sequences in subcategories. J. Algebra 69 (1981), no. 2, 426--454. 

[AS2] M. Auslander, S. O. Smalo: Preprojective modules over Artin algebras. J. Algebra 66 (1980), no. 1, 61--122. 

[AS3] M. Auslander, S. O. Smalo: Preprojective lattices over classical orders. Integral representations and applications (Oberwolfach, 1980), pp. 326--344, Lecture Notes in Math., 882, Springer, Berlin-New York, 1981.

[ASo] M. Auslander, O. Solberg: Relative homology and representation theory. I. Relative homology and homologically finite subcategories. Comm. Algebra 21 (1993), no. 9, 2995--3031.

[B] H. Bass: Finitistic dimension and a homological generalization of semiprimary rings, Trans. Amer. Math. Soc. 95 (1960) 466-488.

[CP] F. U. Coelho, M. I. Platzeck : On the representation dimension of some classes of algebras, preprint.

[CB] W. W. Crawley-Boevey: Functorial filtrations. II. Clans and the Gelfand problem. J. London Math. Soc. (2) 40 (1989), no. 1, 9--30.

[CPS1] E. Cline, B. Parshall, L. Scott: Finite-dimensional algebras and highest weight categories. J. Reine Angew. Math. 391 (1988), 85--99.

[CPS2] E. Cline, B. Parshall, L. Scott: Algebraic stratification in representation categories. J. Algebra 117 (1988), no. 2, 504--521.

[CR] C. W. Curtis, I. Reiner: Methods of representation theory. Vol. I. With applications to finite groups and orders. A Wiley-Interscience Publication. John Wiley \& Sons, Inc., New York, 1990.

[DK] Y. A. Drozd, V. V. Kiri\v cenko: The quasi-Bass orders. (Russian) Izv. Akad. Nauk SSSR Ser.
Mat. 36 (1972), 328--370.

[DKR] Y. A. Drozd, V. V. Kiri\v cenko, A. V. Ro\u\i ter: Hereditary and Bass orders. (Russian) Izv. Akad. Nauk SSSR Ser. Mat. 31 1967 1415--1436.

[DR1] V. Dlab, C. M. Ringel: Indecomposable representations of graphs and algebras. Mem. Amer. Math. Soc. 6 (1976), no. 173.

[DR2] V. Dlab, C. M. Ringel: Quasi-hereditary algebras. Illinois J. Math. 33 (1989), no. 2, 280--291. 

[DR3] V. Dlab, C. M. Ringel: Every semiprimary ring is the endomorphism ring of a projective module over a quasihereditary ring. Proc. Amer. Math. Soc.
107 (1989), no. 1, 1--5.

[DR4] V. Dlab, C. M. Ringel: Auslander algebras as quasi-hereditary algebras. J. London Math. Soc. (2) 39 (1989), no. 3, 457--466.

[EHIS] K. Erdmann, T. Holm, O. Iyama, J. Schr\"oer: Radical embeddings and representation dimension, to appear in Advances in Mathematics.

[FGR] R. M. Fossum, P. Griffith, I. Reiten: Trivial extensions of abelian categories. Homological algebra of trivial extensions of abelian categories with applications to ring theory. Lecture Notes in Mathematics, Vol. 456. Springer-Verlag, Berlin-New York, 1975.

[F] H. Fujita: Minimal injective resolution of an order over a local Dedekind domain. Comm. Algebra 26 (1998), no. 2, 447--451.

[G] P. Gabriel: Indecomposable representations. II. Symposia Mathematica, Vol. XI (Convegno di Algebra Commutativa, INDAM, Rome, 1971), pp. 81--104. Academic Press, London, 1973. 

[GD] P. Gabriel, J. A. de la Pena: Quotients of representation-finite algebras. Comm. Algebra 15 (1987), no. 1-2, 279--307. 


[Ha] D. Happel: Triangulated categories in the representation theory of finite-dimensional algebras. London Mathematical Society Lecture Note Series, 119. Cambridge University Press, Cambridge, 1988.

[Han] Y. Han: Controlled wild algebras. Proc. London Math. Soc. (3) 83 (2001), no. 2, 279--298.

[HS] P. J. Hilton, U. Stammbach: A course in homological algebra. Graduate Texts in Mathematics, 4. Springer-Verlag, New York, 1997.

[H] T. Holm: Representation dimension of some tame blocks of finite groups, to appear in Algebra Colloquium.

[HN] H. Hijikata, K. Nishida: Bass orders in nonsemisimple algebras. J. Math. Kyoto Univ. 34 (1994), no. 4, 797--837.

[I1] O. Iyama: A generalization of rejection lemma of Drozd-Kirichenko. J. Math. Soc. Japan 50 (1998), no. 3, 697--718.

[I2] O. Iyama: $\tau$-categories I,II,III, to appear in Algebras and Representation Theory.

[I3] O. Iyama: A proof of Solomon's second conjecture on local zeta functions of orders, J. Algebra 259 (2003), no. 1, 119-126.

[I4] O. Iyama: Finiteness of Representation dimension, Proc. Amer. Math. Soc. 131 (2003), no. 4, 1011-1014.

[I5] O. Iyama: Representation dimension and Solomon zeta function, to appear in Fields Institute Communications.

[I6] O. Iyama: The relationship between homological properties and representation theoretic realization of artin algebras, to appear in Transactions of the American Mathematical Society.

[I7] O. Iyama: Quadratic bimodules and Quadratic orders, preprint.

[I8] O. Iyama: On stable categories of orders, preprint.

[IT] K. Igusa, G. Todorov: On the finitistic global dimension conjecture, preprint.

[J] V. A. Jategaonkar: Global dimension of triangular orders over a discrete valuation ring. Proc. Amer. Math. Soc. 38 (1973), 8--14. 

[K] S. K\"onig: Every order is the endomorphism ring of a projective module over a quasi-hereditary order. Comm. Algebra 19 (1991), no. 8, 2395--2401.

[KW] S. K\"onig, A. Wiedemann: Global dimension two orders are quasi-hereditary. Manuscripta Math. 66 (1989), no. 1, 17--23.

[M] Y. Miyashita: Tilting modules of finite projective dimension. Math. Z. 193 (1986), no. 1, 113--146.

[N] H. Nagase: Non-strictly wild algebras. J. London Math. Soc. (2) 67 (2003), no. 1, 57--72.

[Re] I. Reiner: Maximal orders. London Mathematical Society Monographs. New Series, 28. The Clarendon Press, Oxford University Press, Oxford, 2003.

[R1] C. M. Ringel: Representations of $K$-species and bimodules. J. Algebra 41 (1976), no. 2, 269--302.

[R2] C. M. Ringel: Tame algebras and integral quadratic forms. Lecture Notes in Mathematics, 1099. Springer-Verlag, Berlin, 1984.

[R3] C. M. Ringel: The Gabriel-Roiter measure, preprint.

[RR] C. M. Ringel, K. W. Roggenkamp: Diagrammatic methods in the representation theory of orders. J. Algebra 60 (1979), no. 1, 11--42. 

[Ro] K. W. Roggenkamp: Lattices over orders. II. Lecture Notes in Mathematics, Vol. 142 Springer-Verlag, Berlin-New York 1970.

[Roi] A. V. Ro\u\i ter: Unboundedness of the dimensions of the indecomposable representations of an algebra which has infinitely many indecomposable representations. (Russian) Izv. Akad. Nauk SSSR Ser. Mat. 32 1968 1275--1282. 

[Rou] R. Rouquier, Dimensions of triangulated categories, preprint.

[Ru1] W. Rump: Lattice-finite rings, preprint.

[Ru2] W. Rump: The category of lattices over a lattice-finite ring, preprint.

[Si] S. Sikko: Resolutions with $n$th syzygies. Comm. Algebra 23 (1995), no. 10, 3729--3739.

[S1] L. Solomon: Zeta functions and integral representation theory. Advances in Math. 26 (1977), no. 3, 306--326.

[S2] L. Solomon: Partially ordered sets with colors. Relations between combinatorics and other parts of mathematics (Proc. Sympos. Pure Math., Ohio State Univ., Columbus, Ohio, 1978), pp. 309--329, Proc. Sympos. Pure Math., XXXIV, Amer. Math. Soc., Providence, R.I., 1979.

[St] B. Stenstr\"om: Rings of quotients Die Grundlehren der Mathematischen Wissenschaften, Band 217. An introduction to methods of ring theory. Springer-Verlag, New York-Heidelberg, 1975.

[Xi1] C. C. Xi: Representation dimension and quasi-hereditary algebras. Adv. Math. 168 (2002), no. 2, 193--212.

[Xi2] C. C. Xi: On the finitistic dimension conjecture I: related to representation-finite algebras, preprint.

[Xi3] C. C. Xi: On the finitistic dimension conjecture II: related to finite global dimension, preprint.

[X] G. Xiangqian: Representation dimension: an invariant under stable equivalence, preprint.

[Y] Y. Yoshino: Cohen-Macaulay modules over Cohen-Macaulay rings. London Mathematical Society Lecture Note Series, 146. Cambridge University Press, Cambridge, 1990.

[Z] H. B. Zimmermann: The finitistic dimension conjectures---a tale of $3.5$ decades. Abelian groups and modules (Padova, 1994), 501--517, Math. Appl., 343, Kluwer Acad. Publ., Dordrecht, 1995.}

\vskip.5em{\footnotesize
{\sc Department of Mathematics, Himeji Institute of Technology, Himeji, 671-2201, Japan}

{\it iyama@sci.himeji-tech.ac.jp}}

\end{document}